\documentclass[10pt,paper=a4]{article}
 \usepackage{indentfirst, latexsym,bm}
 \usepackage{amsmath}
 \usepackage{pifont}
 \usepackage{amsfonts}
 \usepackage{mathrsfs}
 \usepackage{array}
 \usepackage{multirow} 
 \usepackage{graphicx}
 \usepackage{subfigure}
 \usepackage{picinpar}
 \RequirePackage[colorlinks,citecolor=blue,urlcolor=blue,linkcolor=blue]{hyperref}
 
 \RequirePackage[numbers]{natbib}
 \usepackage{enumitem} 

 \usepackage{booktabs}
 \usepackage{authblk}
 
 \usepackage{threeparttable}
 \usepackage{mathrsfs,amsfonts,amsmath}
 
 \usepackage{color}

 \makeatletter
 \def\namedlabel#1#2{\begingroup
	#2%
	\def\@currentlabel{#2}%
	\phantomsection\label{#1}\endgroup
 }
 \makeatother

 \newcommand\email[1]{\rm\href{mailto:#1}{ \nolinkurl{#1}}}

 \renewcommand{\theequation}{\arabic{section}.\arabic{equation}}

 \newtheorem{theorem}{Theorem}[section]
 \newtheorem{definition}[theorem]{Definition}
 \newtheorem{lemma}[theorem]{Lemma}
 \newtheorem{corollary}[theorem]{Corollary}
 \newtheorem{proposition}[theorem]{Proposition}
 \newtheorem{remark}[theorem]{Remark}
 \newtheorem{condition}[theorem]{Condition}
 \newtheorem{example}{Example}[section]

 \def\blemma{\begin{lemma}}\def\elemma{\end{lemma}}
 \def\bproposition{\begin{proposition}}\def\eproposition{\end{proposition}}
 \def\ttheorem{\begin{theorem}}\def\etheorem{\end{theorem}}
 \def\bcorollary{\begin{corollary}}\def\ecorollary{\end{corollary}}
 \def\bremark{\begin{remark}}\def\eremark{\end{remark}}
 \def\bcondition{\begin{condition}}\def\econdition{\end{condition}}

 \def\benumerate{\begin{enumerate}}\def\eenumerate{\end{enumerate}}
 \def\bitemize{\begin{itemize}}\def\eitemize{\end{itemize}}

 \def\beqlb{\begin{eqnarray}}\def\eeqlb{\end{eqnarray}}
 \def\beqnn{\begin{eqnarray*}}\def\eeqnn{\end{eqnarray*}}
 \def\ar{\!\!\!&}

 \def\proof{\noindent{\it Proof.~~}}\def\qed{\hfill$\Box$\medskip}

\begin{document} 
 	 \title{\bf Second-Order Regular Variation and Second-Order Approximation of Hawkes Processes}
 	
 	\author{Ulrich Horst\footnote{Department of Mathematics, and School of Business and Economics,  Humboldt-Universit\"at zu Berlin, Unter den Linden 6, 10099 Berlin; email: horst@math.hu-berlin.de. Horst gratefully acknowledges financial support from the DFG CRC/TRR 388 ``Rough Analysis, Stochastic Dynamics and Related Fields", Project B02.} 	\quad\  and \quad  Wei Xu\footnote{School of Mathematics and Statistics, Beijing Institute of Technology, No. 5, South Street, Zhongguancun, Haidian District, 100081 Beijing; email: xuwei.math@gmail.com}
 	}
 	\maketitle

 \begin{abstract}
 This paper provides and extends second-order versions of several fundamental theorems on first-order regularly varying functions such as Karamata's theorem/representation and Tauberian's theorem. 
 Our results are used to establish second-order approximations for the mean and variance of Hawkes processes with general kernels. Our approximations provide novel insights into the asymptotic behavior of Hawkes processes. They are also of key importance when establishing functional limit theorems for Hawkes processes. 
 \medskip
 \smallskip
    
 \textit{MSC2020 subject classifications.} 
  Primary 
   26A12,   
   40E05;  
  secondary  
   60G55,  
   60K05. 	
     		
  \smallskip
     		
 
 \end{abstract}
  	 
   \section{Introduction}
  \setcounter{equation}{0}
  
  A real-valued function $F$ is called \textit{regularly varying} (at infinity) if it asymptotically behaves like a power function, that is, if  
  \beqnn
  \frac{F(tx)}{F(t)}\to x^\alpha, \quad x>0,
  \eeqnn 
  as $t\to\infty$ for some index of regular variation $\alpha\in\mathbb{R}$. Originally introduced by Jovan Karamata in \cite{Karamata1930}, regular variation has long established itself as powerful mathematical theory for analyzing heavy-tailed phenomena, long-range dependencies and domains of attraction. We refer to \cite{DeHaanFerreira2006,Resnick2007,Samorodnitsky2007,Samorodnitsky2016} and references therein for a comprehensive introduction into the theory of regular variation and its many applications.  
  
  To study  the rate of convergence of a regularly varying function to the limiting power function, de Haan, Resnick and Stadtm\"uller \cite{DeHaanResnick1996,DeHaanStadtmuller1996} introduced the notion of \textit{second-order regularly varying functions}. A regular varying function $F$  with first-order index $\alpha$ is called second-order regularly varying with second-order index $\rho\leq 0$ and auxiliary function $A$ if, as $t\to\infty$, 
  \beqnn
  \frac{F(tx)/F(t)- x^\alpha}{A(t)} \to x^\alpha \int_1^x u^{\rho-1}du,\quad x>0.
  \eeqnn
  
  The subject of second-order regular variation has found numerous applications in statistics and probability theory. It provides an efficient tool for measuring the rate of convergence of the distribution of extreme order statistics  \cite{DeHaanPeng1997,DeHaanPeng1999,DeHaanResnick1996}, for characterizing the asymptotic normality of Hill estimators  \cite{GelukdeHaanResnickStarica1997,ResnickStarica1997}, and for establishing the second-order expansion of tail probabilities of sums of random variables \cite{GelukdeHaanResnickStarica1997,MaoHu2013,MaoNg2015}. Second-order regular variation has also been successfully used to model extreme environmental events \cite{DeHaanDeRonde1998}, to assess tail risks in financial markets \cite{ColesHeffernanTawn1999,EmbrechtsKluppelbergMikosch2003}, to estimate losses from rare catastrophic events \cite{NecirMeraghniMeddi2007} and to establish second-order approximations for risk and concentration measures \cite{DegenLambriggerSegers2010,HuaJoe2011}.
  
  Motivated by its many applications, second-order regular variation has long received considerable attention in the more theoretical probability and statistics literature.  
  For instance, Geluk et al. \cite{GelukdeHaanResnickStarica1997} first proved that the second-order regular variation of two i.i.d random variables carries over to  their maxima and sum. Their result has been extended in \cite{DegenLambriggerSegers2010,LengHu2014,LiuMaoHu2017,MaoHu2013,MaoNg2015} to finite and random sums. They also showed the equivalence between the second-order regular variation and asymptotic normality of the Hill estimator. It has been shown in  \cite{LvMaoHu2012,PanLengHu2013} that second-order regular variation of functions is preserved under integration and composition of functions. 
  Mao and Hua \cite{MaoHua2016} proved the equivalence of the second-order regular variation of a tail-distribution function and its Laplace-Stieltjes transform. A refinement of second-order regular variation, known as \textit{second-order extended regular variation}, was introduced and studied extensively in \cite{DeHaanFerreira2006,DeHaanStadtmuller1996}. Its differences and connections to second-order regular variation have been systematically examined in \cite{Mao2013,Neves2009}. We refer to Appendix B in \cite{DeHaanFerreira2006} for details. 
  
  In the first part of this work, motivated by applications of second-order regular variation to renewal theory and self-exciting stochastic systems, we complement and extend the literature on second-order versions of key theorems for first-order regularly varying functions, including Karamata’s theorem and Karamata's Tauberian theorem. Some of our results have been partially considered in the aforementioned references, yet only under various sufficient conditions or merely for distribution functions. In particular,  the second-order version of Karamata’s theorem in \cite{PanLengHu2013} and the second-order Karamata Tauberian theorem in \cite{MaoHua2016}  have only been established for second-order regularly varying distributions with negative second-order index. As a result, the existing literature on second-order regular variation cannot be applied to establish higher-order asymptotic properties of renewal or self-exciting processes, such as Hawkes processes. 
  
  Our main motivation for analyzing second-order regular varying functions is to establish second-order approximations for the mean and variance of Hawkes processes\,\footnote{In our accompanying paper \cite{HorstXu2023} these approximations play a key role in establishing functional limit theorems for Hawkes processes with long-range dependencies.}.  A Hawkes process $N:=\{N(t):t\geq 0\}$ is a random point process that models self-exciting arrivals of random events. Its \textit{intensity} $\Lambda:=\{\Lambda(t):t\geq 0\}$ is usually of the form
  \beqnn
  \Lambda(t):= \mu(t)+ \sum_{0<\tau_i<t} \phi(t-\tau_i) = \mu(t) + \int_{(0,t)} \phi(t-s)N(ds), 
  \eeqnn
  for some \textit{immigration density} $\mu \in    L^1_{\rm loc} (\mathbb{R}_+;\mathbb{R}_+)$ that captures the immigration of exogenous events, and some \textit{fertility/activation function}
  $\phi\in    L^1 (\mathbb{R}_+;\mathbb{R}_+)$ that captures the self-exciting impact of past events on the arrivals of future events. 
  The random variable $\tau_i$ denotes the arrival time of the  $i$-th event, for each $i \in \mathbb N$. 
  
  First introduced by Hawkes in \cite{Hawkes1971a,Hawkes1971b} to model cross-dependencies between earthquakes and their aftershocks, Hawkes processes and their generalizations have become a powerful tool to model a variety of phenomena in biology and neuroscience \cite{HorstXu2021,Pernice2011}, sociology and criminology \cite{Crane2008,Mohler2011} and, in particular, finance. Applications in finance range from intraday transaction dynamics \cite{Bauwens2009,Bowsher2007} to asset price processes \cite{BacryDelattreHoffmannMuzy2013} and volatility modeling  \cite{ElEuchFukasawaRosenbaum2018,HorstXu2022,JaissonRosenbaum2015,JaissonRosenbaum2016}, and from limit order book modelling  \cite{HorstXu2019} to financial contagion \cite{SahaliaCacho-DiazLaeven2015,Giesecke2011,JorionZhang2009}. We refer to  \cite{BacryMastromatteoMuzy2015} for a review of Hawkes processes and their applications in finance. 
  
  %
  %
  Under suitable first/second-order regular variation conditions on the kernel function, we establish second-order approximations for both the mean $\mathbf{E}[N(t)]$ and the variance ${\rm Var}(N(t))$ of a {\sl non-stationary} Hawkes process; some asymptotics for {\sl stationary} Hawkes processes are established in, e.g.~\cite[Section 12.4]{Bremaud2020}.
  Our results identify the average number of child events triggered by a mother event as well as the release of self-excitation, and the dispersion of child events as the key determinants of the asymptotic behavior of Hawkes processes. 
  Specifically, we prove that both the mean and the variance of a Hawkes process can always be approximated by a polynomial function modified by a regularly varying function, i.e., as $t\to\infty$,
  \beqnn
  \mathbf{E}[N(t)] \sim C_1 \cdot t^{\alpha_1} + \varepsilon_1(t)
  \quad\mbox{and}\quad 
  {\rm Var}(N(t)) \sim C_2 \cdot t^{\alpha_2} + \varepsilon_2(t). 
  \eeqnn
  Here, $\alpha_1,\alpha_2$ are two positive constants determined by the average number of child events respectively, the release of self-excitation, and $ \varepsilon_1, \varepsilon_2$ are two regularly varying functions that describe the dispersion of child events. Depending on the tail behaviour of the kernel function, the functions  $\varepsilon_i$ may be ``very close to the first-order approximation'' in the sense that $\varepsilon_i(t) = t^{\alpha_i - \eta} \ell_i(t)$ for some $\eta \in [0,\alpha_i]$ and some slowly varying function $\ell_i$ that converges (possibly slowly) to zero as $t \to \infty$. We emphasize that $\eta$  may be zero. In other words, the first-order approximations $C_i t^{\alpha_i}$ of the mean and variance may not  be good approximations, except in the very long run, in which case the second-order correction terms $\varepsilon_i$ need to be accounted for when approximating $\mathbf{E}[N(t)]$ and ${\rm Var}(N(t))$.    
  
  As a concrete example, we provide the exact second-order approximation for Hawkes processes with Mittag-Leffler kernel, which are also known as a \textit{fractional Hawkes process}. Fractional Hawkes processes have been widely used to model the complex systems that enjoy both self-exciting property and long-range dependence; see \cite{BondiPulidoScotti2024,DavisBaeumerWang2024,KetelbutersHainaut2022}.

  \medskip
  {\it \textbf{Organization of this paper.}} In Section~\ref{Preliminaries}, we recall some notation, definitions and elementary results from the theory of regular variation. 
  The second-order versions of Karamata's theorem and representation as well as some byproducts are formulated in Section~\ref{2RVKaraThmRep}. 
  In Section~\ref{Sec.KaramataTauberianThm}, we establish the second-order Karamata Tauberian theorem and the second-order Wiener-Tauberian theorem in Pitt's form. 
  The first- and second-order approximations are established in Section~\ref{Sec.HawkesP} for Hawkes processes under regular variation conditions.

  {\bf Notation.}
  Let $\mathbb{N}=\{0,1,2,\cdots  \}$, $\Gamma(z_1)$ be the gamma function with parameter $z_1 \neq 0$ and $\mathrm{B}(z_2,z_3)$ be the beta function with parameter $z_2,z_3 >0 $:
  \beqnn
  \mathrm{B}(z_2,z_3) = \int_0^1 t^{z_2-1}(1-t)^{z_3 - 1}dt = \frac{\Gamma(z_2) \Gamma(z_3) }{\Gamma(z_2+z_3)} . 
  \eeqnn 
  We make the convention that for any $t_1\leq t_2$,
  \beqnn
  \int_{t_1}^{t_2} = \int_{(t_1,t_2]}
  \quad \mbox{and}\quad
  \int_{t_1}^\infty = \int_{(t_1,\infty)}.
  \eeqnn
  For each $k\geq 1$, let $\mathcal{I}^k_f$ be the $k$-th repeated integral of $f$ with base point $0$ defined by
  \beqnn
  \mathcal{I}^k_f(t):= \int_0^t\int_0^{t_1}\cdots\int_0^{t_{k-1}}f(t_k)dt_k\cdots dt_2 dt_1 
   = \frac{1}{(k-1)!} \int_0^t (t-s)^{k-1}f(s)ds, \quad t\geq 0,
  \eeqnn
  and $\mathcal{I}_f=\mathcal{I}^1_f$ for convention. 
  Here the repeated integral is over the interval $0\leq t_k\leq \cdots\leq t_1\leq t$.
  Let $\hat{f}$ be the Laplace-Stieltjes transform of function  $f \in L^1_{\rm loc}(\mathbb{R}_+;\mathbb{R})$
  \beqnn
  \hat{f}(\lambda):= \int_0^\infty \lambda e^{-\lambda t}f(t)dt,\quad \lambda >0.
  \eeqnn 
 For two functions $f,g$ on $\mathbb{R}_+$, we use the notation $f\sim g$ if the ratio $f(t)/g(t) $ tends to $1$ as $t\to \infty$ [{\it or} as $t\to 0$]; we write $f \asymp g$ if there exists positive constants $C_0$ and $C_1$ such that, for large $t > 0$, 
  \beqnn
  C_0\leq \left| \frac{f(t)}{g(t)} \right | \leq  C_1.
  \eeqnn 
  Moreover, we write $f(1/\cdot)$ for the function $f(1/t)$ with $t>0$.
 When the function are locally integrable, we denote by $f*g$ their convolution, i.e.,
  \beqnn
  f*g(t):=\int_0^t f(t-s)g(s)ds,
  \quad
  t\geq 0.
  \eeqnn 
  Throughout this paper, we assume the generic constant $C$ may vary from line to line.


 \section{Preliminaries}\label{Preliminaries}
 \setcounter{equation}{0}
 
 In this section, we introduce some additional notation and elementary properties of second-order regularly varying functions; additional properties that will frequently be used are summarized in Appendix~\ref{Appendix.RV}.  Interested readers are referred to the standard references \cite{BinghamGoldieTeugels1987, GelukdeHaan1987, Resnick2007} for regular variation, and \cite{DeHaanFerreira2006, DeHaanStadtmuller1996, Mao2013} for second-order regular variation.
 
 Let us begin by recalling the definitions of regular variation, second-order regular variation and $\Pi$-variation.
 
 \begin{definition}[Regular variation] \label{Def.RV}
 	A measurable function $F: \mathbb{R}_+\mapsto \mathbb{R}$ is  said to be {\rm regularly varying at infinity} with index $\alpha \in\mathbb{R}$, denoted by $F\in {\rm RV}^\infty_\alpha$, if it has constant sign near infinity and 
 	\beqlb\label{Def.RV.eqn}
 	\lim_{t\to\infty}\frac{F(tx)}{F(t)}= x^\alpha,\quad x>0,
 	\eeqlb
 	In particular, $F$ is also said to be  {\rm slowly varying at infinity}  if $\alpha=0$. 
 \end{definition}

 \begin{definition}[Second-order regular variation]\label{Def.2RV}
 	For a function $F\in {\rm RV}^\infty_\alpha$ with $\alpha\in\mathbb{R}$, if there exist a constant $\rho\in\mathbb{R}$  and an eventually positive or negative function $A$ on $\mathbb{R}_+$ such that  
 	\beqlb\label{Def.2RV.eqn}
 	\lim_{t\to\infty} \frac{F(tx)/F(t)-x^\alpha}{A(t)}  =  x^\alpha \int_1^x u^{\rho-1}du, \quad x>0,
 	\eeqlb
 	then $F$ is said to be of {\rm second-order regular variation at infinity} with first-order index $\alpha$, second-order index $\rho$ and auxiliary function  $A$.  
 	The class of all such functions is denoted ${\rm 2RV}^\infty_{\alpha,\rho}(A)$. 
 \end{definition}
 
 We will frequently use the following result; its proof follows directly from Definition~\ref{Def.2RV}. 
 
 \begin{proposition}\label{Prop.Fr}
 For $\alpha\in\mathbb{R}$, $\rho\leq 0$ and $A\in \mathscr{A}_\rho^\infty$, we have that $F\in 2{\rm RV}_{\alpha,\rho}^\infty(A)$ if and only if $t^\theta \cdot F(t)\in 2{\rm RV}_{\alpha+\theta,\rho}^\infty(A)$ for some and hence all $\theta \in \mathbb{R}$. 
 \end{proposition}
 
 According to Theorem~B.1.3 in \cite[p.362]{DeHaanFerreira2006}, if the limit on the left-hand side of (\ref{Def.RV}) exists, then it must be equal to the right-hand side of (\ref{Def.RV}). If the limit on the left side of (\ref{Def.2RV}) exists and is non-zero, then $A(t)\to 0$ as $t\to \infty$. 
 Moreover, by Theorem~B.2.1 in \cite[p.372]{DeHaanFerreira2006} - and adjusting $A$ proportionally if necessary -  the limit must be equal to the right side of (\ref{Def.2RV}), and $A\in {\rm RV}_\rho^\infty$ with $\rho\leq 0$. In what follows we denote by 
 \begin{equation}
 	\mathscr{A}_\rho^\infty :=\Big \{A \in {\rm RV}^\infty_\rho : \lim_{t\to \infty} A(t) = 0 \Big\}
 \end{equation}
 the space of all regularly varying functions functions with index $\rho$ that converge to zero as $t \to \infty$ and
 adopt the convention that $F\in {\rm 2RV}_{\alpha,-\infty}^\infty(0)$ if $F(t)=C\cdot t^\alpha$ for large $t$ and some constant $C\neq 0$. 
 
 \begin{definition}[$\Pi$-variation]\label{Def.PiV}
 	A measurable function $F: \mathbb{R}_+\mapsto \mathbb{R}$ is said to belong to {\rm the class $\Pi$}, if there exists  an eventually positive or negative function  $A$ on  $\mathbb{R}_+$ such that 
 	\beqlb\label{Def.PiV.eqn}
 	\lim_{t\to\infty} \frac{F(tx)-F(t)}{A(t)} =  \log x,\quad x>0.
 	\eeqlb
 	The class of all such functions is denoted $F\in \Pi^\infty(A)$. We refer to   $A$ as an {\rm auxiliary function} for $F$. 
 \end{definition}
 
 By Theorem~B.2.7 in \cite[p.362]{DeHaanFerreira2006}, the auxiliary function $A$ in Definition~\ref{Def.PiV} is slowly varying at infinity.  The next proposition demonstrates the equivalence of the function spaces ${\rm 2RV}^\infty_{0,0}(\cdot)$ and $\Pi^\infty(\cdot)$.
 
 \begin{proposition}\label{Prop.03}
 	For any $A \in \mathscr{A}^\infty_0$,  we have $F \in {\rm 2RV}^\infty_{0,0}(A)$ if and only if $F \in \Pi^\infty(F\cdot A)$. 
 \end{proposition}
 \proof 
 If $F \in {\rm 2RV}^\infty_{0,0}(A)$, by (\ref{Def.2RV.eqn})  with $\alpha=\rho=0$ we have as $t\to\infty$, 
 \beqnn
 \frac{F(tx) - F(t)}{F(t)A(t)}  = \frac{F(tx)/F(t)-x^\alpha}{A(t)}  \to \log x, \quad x>0, 
 \eeqnn 
 and hence $F \in \Pi^\infty(F\cdot A)$. 
 For the converse,  if $F \in \Pi^\infty(F\cdot A)$, by (\ref{Def.PiV.eqn}) with $B=F\cdot A$ we have
 \beqnn
 \frac{F(tx)/F(t)-1}{A(t)}\to \log x,\quad x>0.
 \eeqnn
 as $t\to\infty$, which yields that $F \in {\rm 2RV}^\infty_{0,0}(A)$.
 \qed 
 
 We occasionally refer to regular variation as {\it first-order regular variation} to avoid confusion. 
 The first-order regular variation, second-order regular variation, and $\Pi$-variation of the function $F$ at zero is defined as in (\ref{Def.RV})-(\ref{Def.PiV.eqn}) with the corresponding limits holding as $t\to 0+$. The corresponding function spaces are denoted ${\rm RV}^{0+}_\alpha$, ${\rm 2RV}^{0+}_{\alpha,\rho}(A)$, and $\Pi^{0+}(A)$, respectively. In this case, we have $A \in {\rm RV}^{0+}_{\rho}$ with $\rho\geq 0$. 
 
 The next proposition is a direct consequence of the Definitions~\ref{Def.RV}-\ref{Def.PiV}; its proof is hence omitted. 
 We recall the notation $f(1/\cdot)$ for the function $f(1/t)$ with $t>0$.
 
 \begin{proposition}\label{Prop.01}
 	For $\alpha\in\mathbb{R}$, $\rho\leq 0$, $A_1\in \mathscr{A}_\rho^\infty$ and $A_2\in {\rm RV}^{\infty}_0$, the following hold.
 	\begin{enumerate}
 		\item[(1)] $F\in {\rm RV}^{\infty}_\alpha$ if and only if $F(1/\cdot) \in {\rm RV}^{0+}_{-\alpha}$.
 		
 		\item[(2)] $F\in {\rm 2RV}^{\infty}_{\alpha,\rho}(A_1)$ if and only if $F(1/\cdot) \in {\rm 2RV}^{0+}_{-\alpha,-\rho} \big(-A_1(1/\cdot)\big)$. 
 		
 		\item[(3)] $F \in \Pi^\infty(A_2)$ if and only if $F(1/\cdot)\in  \Pi^{0+} \big(-A_2(1/\cdot) \big)$.  
 	\end{enumerate}
 \end{proposition}
 
 Our focus will be on second-order regular variation at infinity; analogous results for the second-order regular variation at zero can be obtained by using the preceding proposition. We emphasize that only the asymptotics of $F$ at infinity is considered in what follows. To simplify the subsequent statements and proofs, we hence assume without loss of generality that
 \begin{enumerate} 	
 	\item[\namedlabel{H}{\bf(H)}] 
 	\centering{\bf $F$ is positive and locally bounded on $\mathbb{R}_+$} . 
 \end{enumerate}


 \section{Second-order Karamata theorem and representation}\label{2RVKaraThmRep} 
 \setcounter{equation}{0}
 
 In this section we present a general version of the second-order Karamata theorem along with a representation result for second-order regularly varying functions. 
 
 Karamata's theorem (Proposition~\ref{Thm.Karamata}) examines the tail behavior of integrals of regularly varying functions. As many applications call for a more precise analysis of the tail behavior, several authors have analyzed the speed of convergence of the integral functions at infinity. The case of a tail-distribution function $\bar F$ on $\mathbb{R}$ was first considered by Geluk et al. in \cite{GelukdeHaanResnickStarica1997}. A tail-distribution function $\bar F$ is regularly varying at infinity with index $\alpha<0$ if and only if  
 \beqlb\label{enq.Kara.01}
 \lim_{t\to\infty} \frac{ \bar F(t)}{\int_t^\infty s^{-1}\bar F(s)ds} = - \alpha .
 \eeqlb 
 According to Theorem~4.3 in  \cite{GelukdeHaanResnickStarica1997}, for some constant $\rho \leq 0$ and function $A\in\mathscr{A}^\infty_\rho$, the function $\bar F$ belongs to ${\rm 2RV}^\infty_{\alpha,\rho}(A)$ if and only if there exists a function $A_1\in \mathscr{A}^\infty_\rho$ such that
 \beqlb
 \lim_{t\to\infty} \frac{1}{A_1(t)} \Big(  \frac{\bar F(t)}{\int_t^\infty s^{-1}\bar F(s)ds}+ \alpha \Big)  =  C \neq 0. 
 \eeqlb
 An analogous result for general functions has recently been established in \cite[Theorem~3.1]{PanLengHu2013}. As a preparation for the analysis that follows, we first present a minor extension of these results.

 \begin{proposition}[Second-order Karamata theorem]\label{2KaraThm}
 	For $\alpha\in\mathbb{R}$, $\rho\leq 0$, and $A\in \mathscr{A}_\rho^\infty$, the following three statements are equivalent.
 	\begin{enumerate}
 		\item[(1)] $ \displaystyle{F\in 2{\rm RV}_{\alpha,\rho}^\infty(A)}$. 
 		
 		\item[(2)]  For some (and hence all) $\theta>-\alpha-\rho$ and some $t_0\geq 0$ with $\int_{t_0}^t s^{\theta-1}F (s)ds<\infty$ for all $t\geq t_0$, the following limit holds:
 		\beqlb\label{2KaramataThm01}
 		\lim_{t\to\infty} \frac{1}{A(t)} \Big( \frac{t^{\theta}F (t) }{ \int_{t_0}^t s^{\theta-1}F (s)ds} -(\alpha+\theta) \Big) = \frac{\alpha+\theta }{\alpha+\theta+\rho }. 
 		\eeqlb
 		In this case, we also have $\displaystyle{ \int_{t_0}^t s^{\theta -1}F (s)ds\in {\rm 2RV}^\infty_{\alpha+\theta,\rho}  \Big( \frac{\alpha+\theta}{\alpha+\theta+\rho} \cdot A \Big)}$.
 		
 		\item[(3)] For some (and hence all) $\theta<-\alpha$, the following limit holds:
 		\beqlb\label{2KaramataThm02}
 		\lim_{t\to\infty}   \frac{1}{A(t)} \Big( \frac{t^{\theta} F (t) }{ \int_t^\infty s^{\theta-1}F (s)ds} + (\alpha+\theta) \Big)
 		=  -\frac{\alpha+\theta }{\alpha+\theta+\rho }.
 		\eeqlb
 		In this case, we also have  $\displaystyle{ \int_t^\infty s^{\theta -1}F  (s)ds\in {\rm 2RV}^\infty_{\alpha+\theta,\rho}\Big( \frac{\alpha+\theta}{\alpha+\theta+\rho} \cdot A \Big)  } $.
 	\end{enumerate}
 \end{proposition}
 \proof 
 We only prove the equivalence between claims (1) and (2). The equivalence between claims (1) and (3) can be proved in a similar way.  
 
 In view of Proposition~\ref{Prop.Fr}, it suffices to show that $t^{\theta-1}F(t)\in 2{\rm RV}_{\alpha+\theta-1,\rho}^\infty(A)$ if and only if claim (2) holds. The pre-limit on the left-hand side of (\ref{2KaramataThm01}) can be written as
 \beqnn
 \frac{1}{A(t)} \Big(\frac{1}{\alpha+\theta} - \frac{ \int_{t_0}^t s^{\theta-1} F (s)ds}{t^{\theta}F(t) } \Big) \cdot \frac{(\alpha+\theta) \cdot t^{\theta}F(t) }{ \int_{t_0}^t s^{\theta-1} F(s)ds} . 
 \eeqnn 
 By Theorem~3.1(i) in \cite{PanLengHu2013} and Proposition \ref{Thm.Karamata}(1) the above function converges to a non-zero limit as $t \to \infty$ if and only if $t^{\theta-1}F(t)\in 2{\rm RV}_{\alpha+\theta-1,\rho}^\infty(A)$. Moreover, in this case, the limit must be equal to the right-hand side of (\ref{2KaramataThm01}).  
 The second-order regular variation of the integral $\int_{t_0}^t s^{\theta-1}F (s)ds$ at infinity follows from Corollary~3.2 in  \cite{PanLengHu2013}. 
 \qed 
 
 Karamata's theorem yields the well-known Karamata representation of a regularly varying function. For example, a tail-distribution function $\bar F$ belongs to ${\rm RV}^\infty_\alpha$ with $\alpha<0$ if and only if it admits the following representation:
 \beqnn
 \bar F(t)= c(t)\cdot \exp\Big\{ \int_1^t s^{-1} \epsilon(s)ds \Big\} \cdot t^\alpha  , \quad t>0,
 \eeqnn
 where the functions  $c(t): \mathbb{R}_+\mapsto \mathbb{R}_+$ and $ \epsilon(t): \mathbb{R}_+\mapsto \mathbb{R}_+$ satisfy $c(t)\to c\in(0,\infty)$ and $ \epsilon(t)\to0$ as $t\to\infty$; see e.g. Corollary~2.1 in \cite[p.29]{Resnick2007}. 
 This representation was first established by Karamata \cite{Karamata1930} in the continuous case and by Korevaar et al. \cite{KorevaarAardenne-EhrenfestDeBruijn1949} in the measurable case. 
 
 The following result establishes a representation result for second-order regular varying functions, that is, a second-order version of Karamata's representation theorem. 
 
 \begin{theorem}[Second-order Karamata representation theorem]
 	\label{Thm.KaraRep}
 	For $\alpha\in\mathbb{R}$, $\rho\leq 0$ and $A\in\mathscr{A}_\rho^\infty$, we have $F\in 2{\rm RV}_{\alpha,\rho}^\infty(A)$ if and only if there exist two constants $\zeta_1,\zeta_2\neq 0$ and a function $A_1 \in \mathscr{A}^\infty_\rho$ such that   
 	\beqlb \label{eqn.3001}
 	\frac{A(t)}{A_1(t)} \to \rho \zeta_2+1 >0 \quad \mbox{as} \quad  t\to\infty 
 	\eeqlb
 	and 
 	\beqlb\label{eqn.KaraRep}
 	F(t)= \zeta_1\cdot \big(1+ \zeta_2\cdot A_1(t)\big)\cdot \exp\Big\{ \int_1^t s^{-1}A_1(s) ds\Big\}  \cdot t^\alpha  ,\quad t>0. 
 	\eeqlb 
 	Moreover, for any $\vartheta<-\alpha$, the two constants $\zeta_1$ and $\zeta_2$ can be chosen as 
 	\beqlb\label{eqn.1004}
 	\zeta_1 = -(\alpha+\vartheta) \int_1^\infty s^{\vartheta-1} F(s)ds 
 	\quad \mbox{and} \quad 
 	\zeta_2 = \frac{1}{\alpha+\vartheta}. 
 	\eeqlb
 \end{theorem}
 \proof 
 Let us first assume that a function $F$ admits the representation (\ref{eqn.KaraRep}) and prove that it satisfies (\ref{Def.2RV}). For each $t,x>0$, we have
 \beqlb\label{eqn.KaraRep.01}
 \frac{F(tx)}{F(t)}-x^\alpha 
 \ar=\ar  x^\alpha \cdot \frac{1+ \zeta_2\cdot A_1(tx)}{1+ \zeta_2\cdot A_1(t)}\cdot \exp\Big\{ \int_t^{tx} s^{-1}A_1(s) ds\Big\}  -x^\alpha\cr
 \ar=\ar x^\alpha \Big(\Big(\frac{1+ \zeta_2\cdot A_1(tx)}{1+ \zeta_2\cdot A_1(t)}-1\Big)  \exp\Big\{ \int_t^{tx} s^{-1}A_1(s) ds\Big\}   + \exp\Big\{ \int_t^{tx} s^{-1}A_1(s) ds\Big\} -1  \Big)\cr
 \ar =\ar x^\alpha \Big(  \zeta_2\cdot \frac{A_1(tx)-A_1(t)}{1+ \zeta_2\cdot A_1(t)}\exp\Big\{ \int_t^{tx} s^{-1}A_1(s) ds\Big\} + \exp\Big\{ \int_t^{tx} s^{-1}A_1(s) ds\Big\}-1 \Big).
 \eeqlb
 Since $A_1 \in \mathscr{A}^\infty_\rho$, it follows from Proposition~\ref{Thm.UniConver} that as $t\to\infty$,
 \beqnn
 \sup_{s\in[1\wedge x,1\vee x]}\big|A_1(ts)/A_1(t)-s^\rho \big|\to0
 \quad\mbox{and}\quad
 \int_t^{tx} s^{-1}A_1(s) ds =\int_1^{x} s^{-1}A_1(ts) ds \to0.
 \eeqnn
 Plugging these limits back into the right-hand side of the last equality in (\ref{eqn.KaraRep.01}), we have as $t\to\infty$,
 \beqnn
 \frac{ F(tx)/F(t)-x^\alpha }{A(t)} 
 \ar=\ar x^\alpha \Big( \zeta_2 \cdot   \frac{ \exp\big\{ \int_t^{tx} s^{-1}A_1(s) ds\big\} }{1+ \zeta_2\cdot A_1(t)}\cdot \frac{A_1(tx)-A_1(t)}{ A(t)}+\frac{1}{A(t)}\Big( \exp\Big\{ \int_t^{tx} s^{-1}A_1(s) ds\Big\}-1\Big) \Big) \cr
 \ar\sim\ar x^\alpha \Big(\zeta_2\cdot \frac{A_1(tx)-A_1(t)}{A(t)}
 + \frac{1}{A(t)} \int_t^{tx} \frac{A_1(s)}{s} ds  \Big)\cr
 \ar=\ar  x^\alpha \Big(\zeta_2\cdot \frac{A_1(tx)-A_1(t)}{A(t)}
 + \int_1^{x} \frac{A_1(ts)}{sA(t)} ds  \Big) \cr
 \ar\to \ar \frac{x^\alpha}{ 1+\rho \zeta_2} \cdot \Big(  \zeta_2(x^\rho -1) +  \int_1^x s^{\rho-1}ds\Big) \cr 
 \ar = \ar x^\alpha \int_1^x s^{\rho-1}ds,
 \eeqnn
 which shows that $F\in 2{\rm RV}^\infty_{\alpha,\rho}(A)$. 
 
 Let us now assume that $F\in 2{\rm RV}^\infty_{\alpha,\rho}(A)$. To establish the desired representation,  we choose a constant $\vartheta < -\alpha$ and define two constants $\zeta_1,\zeta_2$ as in (\ref{eqn.1004}). 
 In view of Proposition~\ref{Prop.Fr}, we have 
 \beqnn
 t^{\vartheta-1} F(t)\in 2{\rm RV}_{\alpha+\vartheta-1,\rho}^\infty(A) \quad \mbox{and}\quad 
 \int_t^\infty s^{\vartheta-1} F(s)ds<\infty, 
 \eeqnn
 for any $t>0$. By Proposition~\ref{2KaraThm}(3), 
 \beqlb\label{FunG}
 \frac{t^{\vartheta}F (t) }{ \int_t^\infty s^{\vartheta-1} F (s)ds} = -\frac{1}{\zeta_2}  -  A_1(t) , \quad t>0,
 \eeqlb
 for some function $A_1\in\mathscr{A}^\infty_\rho$ satisfying (\ref{eqn.3001}).  
 Since $\int_1^t s^{-1} A_1(s)ds<\infty$ for any $t>0$, it follows from (\ref{FunG}) that
 \beqnn
 -\int_1^t s^{-1}A_1(s)ds
 \ar =\ar \int_1^t  \frac{s^{\vartheta-1} F (s) ds}{ \int_s^\infty r^{\vartheta-1} F(r)dr} + \int_1^t \frac{s^{-1}}{\zeta_2}ds
 = \log  \int_1^\infty s^{\vartheta-1} F (s)ds- \log  \int_t^\infty s^{\vartheta-1} F (s)ds +\frac{\log t}{\zeta_2},
 \eeqnn
 and so (using the definition of $\zeta_2$) 
 \beqnn
 \int_t^\infty s^{\vartheta-1} F(s)ds = \int_1^\infty r^{\vartheta-1} F(r)dr \cdot \exp\Big\{ \int_1^t s^{-1}A_1(s) ds  \Big\} \cdot  t^{\alpha +\vartheta}.
 \eeqnn
 Taking this back into (\ref{FunG}) and then dividing both sides by $t^{\vartheta}$ yield the desired representation (\ref{eqn.KaraRep}).
 \qed
 
 The above results allow us to characterize the second-order regular variation of a function $F$ in terms of the second-order regular variation of the following two functions. These functions will be used below to establish a second-order Karamata Tauberian theorem. Specifically, for $t>0$, we set 
 \beqlb
 \mathcal{I}_{F,\theta}^{t_0,\uparrow}(t)\ar:=\ar  t^\theta F(t)-\theta\int_{t_0}^t s^{\theta-1} F(s)ds,\label{eqn.401}\\
 \mathcal{I}_{F,\theta}^{\infty,\downarrow}(t)\ar:=\ar   t^\theta F(t)+\theta\int_t^\infty s^{\theta-1} F(s)ds , \label{eqn.402}
 \eeqlb 
 where $\theta \in \mathbb{R}$ and $t_0\geq 0$ are two constants that ensure that the above integrals are well-defined. 
 If $F$ has locally bounded variation and $t^\theta F(t)\to 0$ as $t\to\infty$, then an integration by parts argument yields 
 \beqnn
 \mathcal{I}_{F,\theta}^{t_0,\uparrow}(t)\ar=\ar t_0^\theta\cdot F(t_0)+    \int_{t_0}^t s^{\theta} dF(s)
 \quad \mbox{and}\quad 
 \mathcal{I}_{F,\theta}^{\infty,\downarrow}(t)=  \int_t^\infty s^{\theta} dF(s),\quad t>0.
 \eeqnn

 \begin{theorem}[Extended second-order Karamata theorem]\label{Thm.Gen2Kara}
 	For $\alpha\neq 0$, $\rho\leq 0$ and $A\in \mathscr{A}_\rho^\infty$, the following three statements are equivalent.  
 	\begin{enumerate}
 		\item[(1)] $  \displaystyle{ F \in {\rm 2RV}_{\alpha,\rho}^\infty(A)}$.
 		
 		\item[(2)]  $   \displaystyle{ \mathcal{I}_{F,\theta}^{t_0,\uparrow} \in {\rm 2RV}_{\alpha+\theta,\rho}^\infty\left( \frac{(\alpha+\rho)(\alpha+\theta) }{\alpha(\alpha+\theta+\rho )} \cdot A\right)}$ for some and hence all   $\theta>-\alpha-\rho$.
 		
 		\item[(3)]  $  \displaystyle{\mathcal{I}_{F,\theta}^{\infty,\downarrow} \in {\rm 2RV}_{\alpha+\theta,\rho}^\infty \left( \frac{(\alpha + \rho)(\alpha+\theta) }{\alpha(\alpha+\theta+\rho )} \cdot A \right) }$ for some and hence all $\theta<-\alpha$.
 		
 	\end{enumerate}
 	
 \end{theorem}
 \proof 
 We assume w.l.o.g.~that $\theta\neq 0$. To prove that   (1) implies (2) we rewrite the function  $\mathcal{I}_{F,\theta}^{t_0,\uparrow} $ as
 \beqlb\label{eqn2.6.01}
 \mathcal{I}_{F,\theta}^{t_0,\uparrow} (t)\ar=\ar  \Big(\frac{ t^\theta F (t)}{ \int_{t_0}^t s^{\theta-1}F (s) ds} -\theta\Big) \cdot   \int_{t_0}^t s^{\theta-1} F(s) ds,\quad t\geq t_0. 
 \eeqlb
 If $F \in {\rm 2RV}_{\alpha,\rho}^\infty(A)$, then it follows from Proposition~\ref{Prop.Fr} that $ t^{\theta-1} F(t)\in 2{\rm RV}_{\alpha+\theta-1,\rho}^\infty(A)$, and from  
 Proposition \ref{2KaraThm} that  
 \beqnn
 \frac{ t^\theta F (t)}{ \int_{t_0}^t s^{\theta-1}F (s) ds} = \alpha+\theta+ A_1(t) ,\quad t\geq t_0,
 \eeqnn
 for some $A_1\in\mathscr{A}_\rho^\infty$ with $A_1(t)\sim \frac{\alpha+\theta}{\alpha+\theta+\rho} \cdot A(t)$ as $t\to\infty$. Using the same arguments as in the proof of Theorem~\ref{Thm.KaraRep}, 
 \beqnn
 \int_{t_0}^t s^{\theta-1}F (s) ds = \int_{t_0}^1 s^{\theta-1}F (s) ds \cdot \exp\Big\{ \int_1^t \frac{A_1(s)}{s} ds\Big\}\cdot t^{\alpha+\theta}.
 \eeqnn
 Taking this back into (\ref{eqn2.6.01}) shows that
 \beqnn
 \mathcal{I}_{F,\theta}^{t_0,\uparrow} (t)= \alpha\int_{t_0}^1 s^{\theta-1}F (s)ds \cdot \big(1+ \alpha^{-1}\cdot A_1(t)\big)\cdot \exp\Big\{ \int_1^t \frac{A_1(s)}{s} ds\Big\}\cdot t^{\alpha+\theta}.
 \eeqnn
 Hence,  an application of Theorem~\ref{Thm.KaraRep} with $\zeta_2=1/\alpha$ shows that $\mathcal{I}_{F,\theta}^{t_0,\uparrow}$ is second-order regularly varying at infinity with a first-order index of $\alpha+\theta$, second-order index of $\rho$, and auxiliary function 
 \beqnn
 \Big( \frac{\rho}{\alpha}+1\Big) A_1(t) \sim \frac{(\alpha+\rho)(\alpha+\theta) }{\alpha(\alpha+\theta+\rho )} \cdot A(t) ,
 \eeqnn
 as $t\to\infty$. This proves (2). The proof that (1) implies (3) is similar.

 We now prove that $(2)$ implies $(1)$ if $\alpha<0$. By Fubini's theorem, 
 \beqnn
 \theta \int_t^\infty s^{-\theta -1} \int_{t_0}^s r^{\theta-1}F(r)dr ds 
 \ar=\ar \theta \int_t^\infty s^{-\theta -1}  ds\cdot \int_{t_0}^t r^{\theta-1}F(r)dr 
 +  \theta \int_t^\infty s^{-\theta -1} \int_t^s r^{\theta-1}F(r)dr ds \cr
 \ar=\ar  t^{-\theta }\cdot \int_{t_0}^t r^{\theta-1}F(r)dr  + \theta \int_t^\infty s^{\theta-1}F(s)  \int_s^\infty r^{-\theta -1}dr ds\cr
 \ar=\ar  t^{-\theta }\cdot \int_{t_0}^t r^{\theta-1}F(r)dr  +  \int_t^\infty s^{-1}F(s) ds,
 \eeqnn
 for any $t\geq t_0$. 
 Hence,
 \beqnn
 \int_t^\infty s^{-\theta -1} \mathcal{I}_{F,\theta}^{t_0,\uparrow}(s)ds
 \ar=\ar \int_t^\infty s^{-1} F(s)ds -\theta \int_t^\infty s^{-\theta -1} \int_0^s r^{\theta-1}F(r)dr ds  
 = -   t^{-\theta }\cdot \int_{t_0}^t r^{\theta-1}F(r)dr  .
 \eeqnn
 Taking this back into (\ref{eqn.401}) shows that
 \beqnn
 F(t)\ar=\ar  t^{-\theta}   \mathcal{I}_{F,\theta}^{t_0,\uparrow}(t) +\theta t^{-\theta} \int_{t_0}^t s^{\theta-1} F(s)ds 
 = t^{-\theta}  \mathcal{I}_{F,\theta}^{t_0,\uparrow}(t) -\theta \int_t^\infty s^{-\theta-1} \mathcal{I}_{F,\theta}^{\uparrow}(s)ds =  \mathcal{I}^{\infty,\downarrow}_{ \mathcal{I}_{F,\theta}^{t_0,\uparrow},-\theta}(t). 
 \eeqnn
 Using the fact that $(1)$ implies $(3)$, we see that
 \[
 F=  \mathcal{I}^{\infty,\downarrow}_{ \mathcal{I}_{F,\theta}^{t_0,\uparrow},-\theta} \in {\rm 2RV}_{\alpha,\rho}^\infty(A_2) 
 \]
 with
 \[ 
 A_2 := \frac{(\rho + \theta+\alpha)(\theta+\alpha-\theta) }{(\theta+\alpha)(\theta+\alpha-\theta+\rho )} \cdot   \frac{(\alpha+\rho)(\alpha+\theta) }{\alpha(\alpha+\theta+\rho )} \cdot A  = A
 \]
 Similarly, if $\alpha> 0$, then it follows again from Fubini's theorem that 
 \beqnn
 F(t)= t^{-\theta}  \mathcal{I}_{F,\theta}^{t_0,\uparrow}(t) -\theta \int_{t_0}^t s^{-\theta-1} \mathcal{I}_{F,\theta}^{\uparrow}(s)ds =  \mathcal{I}^{t_0,\uparrow}_{ \mathcal{I}_{F,\theta}^{t_0,\uparrow},-\theta}(t). 
 \eeqnn
 Using the fact that $(1)$ implies $(2)$ we obtain that $F =\mathcal{I}^{t_0,\uparrow}_{ \mathcal{I}_{F,\theta}^{t_0,\uparrow},-\theta}\in {\rm 2RV}^\infty_{\alpha,\rho}(A) $. This shows that $(2)$ implies  $(1)$.  The proof that $(3)$ implies $(1)$ is similar. 
 \qed

 The assumption that the function $A$ is eventually positive or negative guarantees that the integral $ \int_1^\infty s^{-1}A(s)ds$ is well-defined and finite when $\rho<0$. In this case, the following corollary provides an alternative asymptotic representation of a second-order regularly varying function by comparing the asymptotic behavior of the auxiliary function and its integral at infinity. This result was first given in \cite[Lemma~3]{HuaJoe2011}. We provide an alternative proof.
 \begin{corollary}\label{Coro.KaraRep}
 	For $\alpha\in\mathbb{R}$, $\rho< 0$ and $A\in \mathscr{A}_\rho^\infty$, we have $F\in  {\rm 2RV}_{\alpha,\rho}^\infty (A)$ if and only if there exists a constant $C_F\neq 0$ such that as $t\to\infty$, 
 	\beqlb\label{eqn.KaraRep.04}
 	F(t)= C_F\cdot t^\alpha  \cdot \Big(1  + \frac{A(t)}{\rho}+o\big(A(t)\big) \Big).
 	\eeqlb
 \end{corollary}
 \proof 
 Let us first assume that $F$ admits the representation (\ref{eqn.KaraRep.04}). In this case it holds for any $x,t>0$ that  
 \beqnn
 \frac{F(tx)}{F(t)} -x^\alpha 
 \ar=\ar  x^\alpha \cdot  \frac{ (A(tx)-A(t))/\rho+ o(A(tx))-o(A(t))}{1+A(t)/\rho+ o(A(t))} . 
 \eeqnn
 Since $A(tx)/A(t) \to x^\rho$ as $t \to \infty$ it is not difficult to see that
 \beqnn
 \lim_{t\to\infty}  \frac{F(tx)/F(t)-x^\alpha}{A(t)} 
 = x^\alpha \cdot \frac{1}{\rho} \cdot  \lim_{t\to\infty} \frac{A(tx)-A(t)}{A(t)} = x^\alpha \cdot \frac{1}{\rho} \big( x^\rho -1 \big) =  x^\alpha   \int_1^x u^{\rho-1}du,
 \eeqnn
 from which we deduce that $F\in   {\rm 2RV}_{\alpha,\rho}^\infty (A)$. 
 
 Conversely, if $F\in  {\rm 2RV}_{\alpha,\rho}^\infty (A)$, then the representation (\ref{eqn.KaraRep}) holds. In terms of the constant 
 \beqnn
 C_F := \zeta_1\cdot  \exp \Big\{ \int_1^\infty s^{-1}A_1(s) ds \Big\} \neq 0,
 \eeqnn
 we hence get that 
 \beqnn
 F(t)= C_F  \cdot t^\alpha\cdot \big(1+ \zeta_2\cdot A_1(t)\big)\cdot  \exp\Big\{ -\int_t^\infty s^{-1}A_1(s) ds\Big\}  ,\quad t>1.
 \eeqnn 
 A Taylor expansion of the exponential function along with the fact that $\int_t^\infty s^{-1}A_1(s) ds \to 0$ as $t\to\infty$  shows that
 \beqlb\label{eqn.301}
 F(t)
 \ar=\ar C_F \cdot t^\alpha\cdot \big(1+ \zeta_2\cdot A_1(t)\big)\cdot  \Big( 1-  \int_t^\infty s^{-1}A_1(s) ds + o\Big(  \Big| \int_t^\infty s^{-1}A_1(s) ds\Big| \Big)  \Big) . 
 \eeqlb
 Applying Karamata's theorem (Proposition~\ref{Thm.Karamata}) to the function $t^{-1}A_1(t)$, we see that, as $t \to \infty$, 
 \beqnn
 -  \int_t^\infty s^{-1}A_1(s) ds \sim \frac{A_1(t)}{\rho} \sim   \frac{A(t)/\rho}{\zeta_2\rho +1}.
 \eeqnn 
 Plugging these asymptotic results back into the presentation (\ref{eqn.301}), we obtain that
 \beqnn
 F(t)
 \ar=\ar  C_F  \cdot t^\alpha\cdot \Big(1+ (\zeta_2+1/\rho)\cdot A_1(t) +o\big(A_1(t)\big) \Big) = C_F \cdot t^\alpha\cdot \big(1+ A(t)/\rho+ o(A(t)) \big).
 \eeqnn
 \qed 
 
 For any $F\in  {\rm 2RV}_{\alpha,\rho}^\infty (A)$ with $\rho <0$ or $\rho=0$ and $\int_t^\infty s^{-1}A(s)ds<\infty$ for some $t\geq 0$, the representations (\ref{eqn.KaraRep.04}) and (\ref{eqn.KaraRep.03}) imply that $t^{-\alpha}\cdot F(t) \to C_F \neq 0$ as $t\to \infty$. 
 In this case we may conveniently choose 
 $$C_F:=\lim_{t\to\infty}t^{-\alpha}\cdot F(t) \mbox{ whenever it  exits.}$$

 \begin{remark}
 For a function $A\in \mathscr{A}^\infty_0$, the asymptotic behavior of $A$ and the tail-integral function $\int_t^\infty s^{-1} A(s)ds$ at infinity cannot accurately be  compared. 
 This makes it difficult to establish an  analogue of Corollary~\ref{Coro.KaraRep} for second-order regularly varying functions with zero second-order index. 
 	We can, however, still give a representation result. In fact, let $F\in  {\rm 2RV}_{\alpha,0}^\infty (A)$ for some $\alpha\in\mathbb{R} $ and $A\in\mathscr{A}_0^\infty$  with $\int^\infty s^{-1} A(s)ds <\infty $. 
 	Then the representation (\ref{eqn.301}) still holds. Moreover, by Karamata's theorem,
 	\beqnn 
 	A(t)^{-1}\int_t^\infty s^{-1} A(s)ds \to \infty, 
 	\eeqnn
 	as $t\to\infty$. Plugging this back into (\ref{eqn.301}) shows that  
 	\beqlb\label{eqn.KaraRep.03}
 	F(t)= C_F\cdot t^\alpha  \cdot \Big(1 -\int_t^\infty s^{-1}A_1(s)ds \cdot \big(1+o(1) \big)   \Big).
 	\eeqlb 
 \end{remark}
 
 
 The next corollary shows that powers of second-order regularly varying function are second-order regularly varying. 
 
 \begin{corollary}
 	For $\alpha\in\mathbb{R}$, $\rho< 0$ and $A\in \mathscr{A}_\rho^\infty$, if $F\in  {\rm 2RV}_{\alpha,\rho}^\infty (A)$, then $|F|^\theta\in   {\rm 2RV}_{\theta\alpha,\rho}^\infty (\theta\cdot A) $ for any $\theta\neq 0$. 
 \end{corollary}
 \proof 
 By using the Taylor expansion of the function $(1+x)^\theta$, we have as $t\to\infty$, 
 \beqnn
 |F(t)|^\theta = |C_F|^\theta\cdot t^{\theta\alpha}\cdot \big(1+ \theta A(t)/\rho + o(A(t)) \big).
 \eeqnn
 Hence, the result follows from Corollary~\ref{Coro.KaraRep}.  
 \qed 
 
 Corollary~\ref{Coro.KaraRep} also allows us to show that the convolution 
 \beqnn
 F_1*F_2(t): = \int_0^t F_1(t-s) F_2(s)ds ,\quad t\geq 0.
 \eeqnn	      
 of two second-order regularly varying functions $F_1$ and $F_2$ is also second-order regularly varying. The second-order regular variation of the convolution tail,  
 \beqnn
 1- \int_0^t F_1(t-s) dF_2(s),\quad t\geq 0,
 \eeqnn
 of two probability distribution functions $F_1$ and $F_2$ on $\mathbb{R}_+$ has previously been analyzed for the benchmark case $F_1=F_2$ in \cite{GelukdeHaanResnickStarica1997} and for the general case $F_1 \neq F_2$ in \cite{BarbeMCormick2005,LiuMaoHu2017}. 
 
 For two functions $f,g$ on $\mathbb{R}_+$, we write $f \asymp g$ if there exists positive constants $C_0$ and $C_1$ such that, for large $t > 0$, 
 \beqnn
 C_0\leq \left| \frac{f(t)}{g(t)} \right | \leq  C_1.
 \eeqnn
 
 \begin{corollary}\label{Coro.2.6}
 	For $i\in\{1,2\}$, let $F_i\in {\rm 2RV}^\infty_{\alpha_i,\rho_i} (A_i)$ with $\alpha_i>-1$, $\rho_i\in(-\alpha_i-1,0)$ and $A_i\in\mathscr{A}^\infty_{\rho_i}$. Let
 	\beqlb\label{eqn.3005}
 	A_0(t) := 	\frac{\mathrm{B}(\alpha_1+\rho_1+1,\alpha_2+1)}{\mathrm{B}(\alpha_1+1,\alpha_2+1)} \cdot A_1(t)+ \frac{\mathrm{B}(\alpha_1+1,\alpha_2+\rho_2+1)}{\mathrm{B}(\alpha_1+1,\alpha_2+1)}  \cdot A_2(t)  ,\quad t>0,
 	\eeqlb   
 	where $\mathrm{B}(\cdot,\cdot)$ denotes the beta function. 
 	If $ A_0\asymp  |A_1|+|A_2|$, then $F_1*F_2\in {\rm 2RV}^\infty_{\alpha_1+\alpha_2+1,\rho_1\vee \rho_2} (A_0)$.
 \end{corollary} 
 \proof 
 For each $i\in\{1,2\}$, it follows from Corollary~\ref{Coro.KaraRep} that
 \beqnn
 F_i(t) = C_{F_i} \big( t^{\alpha_i} + G_i(t) + G^\circ_i(t) \big) 
 \eeqnn
 where $C_{F_i}\neq 0$, $G_i(t)  = t^{\alpha_i}  A_i(t)/\rho_i  \in {\rm RV}^\infty_{\alpha_i+\rho_i}$ and $G_i^\circ$ is a locally bounded function that satisfies $G_i^\circ(t)=o(G_i(t))$ as $t\to \infty$.   
 The convolution of $F_1$ and $F_2$ can thus be expressed as $F_1*F_2(t) =   C_{F_1}C_{F_2}\cdot \sum_{i=1}^9 H_i(t)$ with
 \beqnn
 H_1(t)\ar:=\ar \int_0^t (t-s)^{\alpha_1} s^{\alpha_2}ds,\qquad \ \,
 H_2(t) :=  \int_0^t G_1(t-s) s^{\alpha_2}ds,\quad\ \ \ \,
 H_3(t):= \int_0^t (t-s)^{\alpha_1} G_2(s)ds,\quad \\
 H_4(t)\ar:=\ar \int_0^t G_1(t-s)G_2(s)  ds,\quad 
 H_5(t):=  \int_0^t G_1^\circ(t-s) s^{\alpha_2}ds,\quad \ \ \ \,
 H_6(t):= \int_0^t (t-s)^{\alpha_1} G_2^\circ(s)ds,\quad\\
 H_7(t)\ar:=\ar  \int_0^t G_1^\circ(t-s) G_2(s)ds,\quad 
 H_8(t):= \int_0^t G_1(t-s) G_2^\circ(s)ds,\quad 
 H_9(t):= \int_0^t G_1^\circ(t-s) G_2^\circ(s)ds.
 \eeqnn
 A change of variables yields that
 \beqnn
 H_1(t) = t^{\alpha_1+\alpha_2+1} \int_0^1 (1-s)^{\alpha_1}  s^{\alpha_2} ds = \mathrm{B}(\alpha_1+1,\alpha_2+1)\cdot  t^{\alpha_1+\alpha_2+1} . 
 \eeqnn
 Applying Corollary~\ref{Appendix.Coro01} to $H_2$, $H_3$ and $H_4$, we have as $t\to\infty$,
 \beqnn
 H_2(t) \ar\sim\ar \mathrm{B}(\alpha_1+\rho_1+1,\alpha_2+1) \cdot t^{\alpha_2+1}\cdot G_1(t) 
 = \mathrm{B}(\alpha_1+\rho_1+1,\alpha_2+1) \cdot t^{\alpha_1+\alpha_2+1}\cdot \frac{A_1(t)}{\rho_1}, \\
 H_3(t) \ar\sim\ar \mathrm{B}(\alpha_1+1,\alpha_2+\rho_2+1) \cdot t^{\alpha_1+1}\cdot G_2(t) 
 = \mathrm{B}(\alpha_1+1,\alpha_2+\rho_2+1) \cdot t^{\alpha_1+\alpha_2+1}\cdot \frac{A_2(t)}{\rho_2}, \\
 H_4(t)\ar\sim\ar  \mathrm{B}(\alpha_1+\rho_1+1,\alpha_2+\rho_2+1) \cdot t \cdot G_1(t) G_2(t)
 = \mathrm{B}(\alpha_1+\rho_1+1,\alpha_2+\rho_2+1) \cdot t^{\alpha_1+\alpha_2+1}\cdot \frac{A_1(t)A_2(t)}{\rho_1\rho_2}.
 \eeqnn
 
 Let us now consider the function $H_5$. For any $\epsilon\in(0,1)$, there exists a constant $t_0>0$ such that $$|G_1^\circ(t)|\leq \epsilon \cdot |G_1(t)|$$ for any $t\geq t_0$. 
 By a change of variables and the local boundedness of $G_1^\circ$, we have for any $t\geq 2t_0$, 
 \beqnn
 H_5(t) = \int_0^t G_1^\circ(s) (t-s)^{\alpha_2} ds  
 \ar \leq \ar  \epsilon \int_{t_0}^t G_1 (s) (t-s)^{\alpha_2} ds + \sup_{r\in[0,t_0]}|G_1^\circ(r)|  \int_0^{t_0}  (t-s)^{\alpha_2} ds 
 \leq   \epsilon \cdot H_2(t)+ C\cdot t^{\alpha_2},
 \eeqnn
 for some constant $C>0$ independent of $\epsilon$. 
 The preceding result and the fact that $\alpha_1+\rho_1+1>0$ shows that 
 \beqnn 
 \limsup_{t\to \infty} \Big|\frac{H_5(t)}{t^{\alpha_1+\alpha_2+1} A_1(t)}\Big| \leq  \mathrm{B}(\alpha_1+\rho_1 +1, \alpha_2+1) \cdot  \epsilon
 \eeqnn
 and hence that 
 \beqnn
 H_5(t)=o(t^{\alpha_1+\alpha_2+1} A_1(t)) \quad \mbox{as} \quad t \to \infty.
 \eeqnn	
 
 Similarly, $G_6(t)= o(t^{\alpha_1+\alpha_2+1} A_2(t))$, while the functions $G_7(t)$, $G_8(t)$, $G_9(t)$ are all of order $o(t^{\alpha_1+\alpha_2+1} A_1(t)A_2(t))$ as $t\to\infty$. As a result, 
 \beqlb\label{eqn.0005}
 \lefteqn{F_1*F(t) =  C_{F_1}C_{F_2}\cdot t^{\alpha_1+\alpha_2+1} \cdot \mathrm{B}(\alpha_1 +1, \alpha_2 +1)} \quad \ar\ar \cr
 \ar\ar\cr
 \ar\ar \times \Big( 1+ \frac{ \mathrm{B}(\alpha_1+\rho_1 +1, \alpha_2+1)}{\mathrm{B}(\alpha_1 +1, \alpha_2 +1)} \cdot \frac{A_1(t)}{\rho_1} + \frac{ \mathrm{B}(\alpha_1 +1, \alpha_2+\rho_2+1)}{\mathrm{B}(\alpha_1 +1, \alpha_2 +1)} \cdot \frac{A_2(t)}{\rho_2} + o\big( |A_1(t)| + |A_2(t)| \big) \Big) .
 \eeqlb
 If $\alpha_1\neq \alpha_2$ or $\rho_1\neq \rho_2$, then $ A_0\asymp  |A_1|+|A_2|$ and, as $t\to\infty$, 
 \beqnn
 A_0(t) \sim  
 \begin{cases}
 	\displaystyle{	\frac{\mathrm{B}(\alpha_1+\rho_1+1,\alpha_2+1)}{\mathrm{B}(\alpha_1+1,\alpha_2+1)} \cdot A_1(t),} & \mbox{  if }\rho_1<\rho_2; \vspace{5pt} \\
 	\displaystyle{	\frac{\mathrm{B}(\alpha_1+1,\alpha_2+\rho_2+1)}{\mathrm{B}(\alpha_1+1,\alpha_2+1)} \cdot A_2(t),} & \mbox{  if }\rho_1>\rho_2.
 \end{cases} 
 \eeqnn
 If $\alpha_1 = \alpha_2  $ and $\rho_1=\rho_2$, then
 $$A_0(t)=\frac{\mathrm{B}(\alpha_1+\rho_1+1,\alpha_2+1)}{\mathrm{B}(\alpha_1+1,\alpha_2+1)} \cdot \big(A_1(t)+A_2(t)\big), \quad t>0 $$
 and $ A_0\asymp  |A_1|+|A_2| $ if and only if $A_1+A_2\asymp |A_1|+|A_2|$. 
 In any case,  if $ A_0\asymp  |A_1|+|A_2| $, by (\ref{eqn.0005}) we have 
 \beqnn
 F_1*F(t) =   C_{F_1}C_{F_2}\cdot t^{\alpha_1+\alpha_2+1} \cdot \mathrm{B}(\alpha_1 +1, \alpha_2 +1) 
 \Big( 1 + \frac{A_0(t)}{\rho_1\vee \rho_2} + 
 o\big( |A_0(t)|\big) \Big),
 \eeqnn
 then it follows from Corollary~\ref{Coro.KaraRep}  that 
 \beqnn
 F_1*F_2 \in {\rm 2RV}^\infty_{\alpha_1+\alpha_2+1,\rho_1\vee \rho_2}(A_0).
 \eeqnn
 \qed

  
    \section{Second-order Karamata Tauberian theorem}\label{Sec.KaramataTauberianThm}
  \setcounter{equation}{0}
  
  The Laplace-Stieltjes transform of a locally integrable function $F$ is defined as 
  \beqnn
  \hat{F} (\lambda ) :=  \int_0^\infty \lambda e^{-\lambda t}F(t)dt,\quad \lambda \geq 0.
  \eeqnn
  If $F$ has locally bounded variation, then it is convenient to consider its Laplace transform 
  \beqnn
  \int_0^\infty  e^{-\lambda t}d F(t),\quad \lambda \geq 0. 
  \eeqnn  
  
  The connection between the first-order regular variation of $F$ and that of its Laplace-Stieltjes transform $\hat{F}$ has been extensively investigated in the literature on Abelian and Tauberian theorems. For instance, \textit{Karamata's Tauberian theorem} (Proposition~\ref{Thm.KaramataTauberian}) states that if $\alpha>-1$ and $F$ is eventually monotone, then $F\in {\rm RV}_\alpha^\infty$ if and only if $\hat{F}(1/\cdot) \in {\rm RV}_\alpha^\infty$. In this case, $\hat{F}(1/\cdot)\sim \Gamma(1+\alpha) \cdot F(\cdot)$ at infinity.  
  
  In this section we establish a second-order version of Karamata's Tauberian theorem. 
  In particular, we prove that (under mild conditions) $F$ is second-order regularly varying if and only if $\hat{F}(1/\cdot)$ is second-order regularly varying. 
  To this end, we introduce, for any  $\alpha\in \mathbb{R}$, the set 
  \begin{equation}
  	\begin{split}
  		\mathscr{M}^\infty_{\alpha,-} & := \left\{ F:\mathbb{R}_+ \to \mathbb{R} : \mbox{ $F$ is locally bounded, eventually positive or negative,} \right. \\
  		& \qquad \qquad ~~~~ \left. \mbox{$ t^{-\alpha}\cdot F(t) \to C_F \neq 0$ as $t \to \infty$, and $  F(t)- C_F \cdot t^\alpha$ is eventually monotone.} \right\}
  	\end{split}	
  \end{equation}	
  
  
  We  start by establishing a second-order Karamata Tauberian theorem for second-order regularly varying functions with negative second-oder index.

  \begin{theorem}[Second-order Karamata Tauberian theorem: $\rho<0$] \label{Thm.2KaraTaub01}
  	For $\alpha>-1$, $\rho\in (-1-\alpha, 0)$ and $A\in\mathscr{A}^\infty_\rho$,   the following hold.
  	\begin{enumerate}
  		\item[(1)] If $F\in {\rm 2RV}_{\alpha,\rho}^\infty (A)$, then $ \displaystyle{\hat{F}(1/\cdot) \in {\rm 2RV}_{\alpha,\rho}^{\infty}\Big(\frac{\Gamma(\alpha+\rho +1)}{\Gamma(\alpha+1)}\cdot A\Big)}$.
  		
  		\item[(2)] If $F\in \mathscr{M}^\infty_{\alpha,-}$ and  $ \displaystyle{\hat{F}(1/\cdot) \in {\rm 2RV}_{\alpha,\rho}^{\infty}\Big(\frac{\Gamma(\alpha+\rho +1)}{\Gamma(\alpha+1)}\cdot A\Big) }$, then $F\in {\rm 2RV}_{\alpha,\rho}^\infty (A)$.
  	\end{enumerate}   
  	
  \end{theorem} 
  \proof  
  If $F\in {\rm 2RV}_{\alpha,\rho}^\infty (A)$, then it follows from Corollary~\ref{Coro.KaraRep} that  $t^{-\alpha} F(t)\to C_F\neq 0$ and 
  \beqnn
  G(t):=F(t)-C_F \cdot t^{\alpha}= \frac{C_F}{\rho} \cdot t^{\alpha} \cdot A(t)\big(1+o(1)\big) \in {\rm RV}_{\alpha+\rho}^\infty. 
  \eeqnn 
  Since $G $ is locally integrable, it follows from Karamata's Tauberian theorem (Proposition~\ref{Thm.KaramataTauberian}) that, as $\lambda \to \infty$, 
  \beqnn
  \hat{G}(1/\lambda) = \int_0^\infty  \frac{1}{\lambda}  e^{-x/\lambda }G(t)dt
  \ar=\ar \frac{C_F}{\rho}\Gamma(\alpha+\rho +1) 
  \cdot \lambda^{\alpha}A(\lambda)  \cdot \big(1+o(1)\big)
  \in {\rm RV}^{\infty}_{\alpha+\rho} . 
  \eeqnn
  Additionally, by using integration by parts,  
  \beqlb\label{eqn.002}
  \int_0^\infty  \frac{1}{\lambda} e^{-t/\lambda }  t^\alpha dt = \lambda^{\alpha} 	\int_0^\infty     t^\alpha e^{- t}dt =  \Gamma(\alpha+1)\lambda^{\alpha},
  \quad \lambda >0.
  \eeqlb 
  Putting these two results together, we see that the Laplace-Stieltjes transform of $F$ satisfies, as $\lambda \to \infty$,
  \beqnn
  \hat{F}(1/\lambda)= C_F  \Gamma(\alpha+1)  \cdot \lambda^{\alpha}  \cdot \Big(1+ \frac{\Gamma(\alpha+\rho +1)}{\Gamma(\alpha+1)} \cdot \frac{A(\lambda)}{\rho} \cdot  \big(1+o(1)\big)   \Big) .
  \eeqnn
  Using Corollary~\ref{Coro.KaraRep} again, we conclude that $\hat{F}(1/\cdot) \in {\rm 2RV}_{\alpha,\rho}^{\infty}\big(\frac{\Gamma(\alpha+\rho +1)}{\Gamma(\alpha+1)}A\big)$. 
  
  To prove claim (2) we assume that $F \in \mathscr{M}^\infty_{\alpha,-}$, from which we deduce that $F(t)\sim C_F\cdot t^\alpha$ as $t\to\infty$. 
  Using Karamata's Tauberian theorem (Proposition~\ref{Thm.KaramataTauberian}) again, we get that
  \beqnn
  \hat{F} (1/\lambda )\sim \Gamma(1+\alpha)F(\lambda) \sim \Gamma(1+\alpha)\cdot C_F\cdot\lambda^{\alpha} = C_F\int_0^\infty \frac{1}{\lambda}e^{-t/\lambda} t^\alpha dt,
  \eeqnn
  as $\lambda\to\infty$. 
  In view of the equality (\ref{eqn.002}) and Corollary~\ref{Coro.KaraRep}, the Laplace-Stieltjes transform of the function $F(t) - C_F t^\alpha$ satisfies 
  \beqnn
  \int_0^\infty \frac{1}{\lambda}e^{-t/\lambda} \big(F(t)-C_F \cdot t^\alpha \big)dt  \ar=\ar \hat{F}(1/\lambda)-\Gamma(1+\alpha) C_F \cdot \lambda^{\alpha} \cr
  \ar=\ar
  \frac{ C_F }{\rho} \Gamma(\alpha+\rho +1) \cdot \lambda^\alpha A(\lambda) \cdot \big(1+o(1)\big) \in {\rm RV}_{\alpha+\rho}^\infty .
  \eeqnn
  Since the function $F(t)-C_F\cdot t^\alpha$ is eventually monotone, an application of Proposition~\ref{Thm.KaramataTauberian} shows that   
  \beqnn
  F(t)-C_F t^\alpha = \frac{C_F}{\rho} \cdot t^\alpha A(t) \cdot \big( 1+ o(1) \big) \in {\rm RV}_{\alpha+\rho}^\infty, 
  \eeqnn
  and  hence it follows from Corollary~\ref{Coro.KaraRep} that $F\in {\rm 2RV}_{\alpha,\rho}^\infty(A)$. 
  \qed

  Establishing a second-order version of Karamata's Tauberian theorem for second-order regularly varying functions with zero second-order index is more challenging. As a preparation we first extend the Wiener-Tauberian theorem in Pitt's form to eventually monotone functions. The case of non-decreasing functions has been established in \cite[Theorem 2.34]{GelukdeHaan1987}. 
  
  Before presenting the theorem, we recall several important auxiliary definitions. 
  A function $F$ is called \textit{slowly decreasing at infinity} if
  \beqnn
  \lim_{u\to 1+}\liminf_{t\to\infty} \inf_{x\in[1,u]} \big\{F(tx) - F(t) \big\} \geq 0,
  \eeqnn
  and  \textit{slowly increasing at infinity} if $-F$ is slowly decreasing at infinity.
  For a  locally integrable function $\mathrm{K}$ on $\mathbb{R}_+$, the {\it Mellin convolution} of $F$ associated with $\mathrm{K}$ is the function $\mathrm{K}\overset{\rm M}*F$ on $(0,\infty)$ given by
  \beqlb
  \mathrm{K}\overset{\rm M}*F(\lambda) := \int_0^\infty \mathrm{K}(x)F(\lambda x)dx = \int_0^\infty \frac{\mathrm{K}(x/\lambda)}{\lambda}F(x)dx ,\quad \lambda >0.
  \eeqlb
  Here, $\mathrm{K}$ is usually referred to as the \textit{kernel}.  
  If $\mathrm{K}(x)=e^{-x}$ is the exponential kernel,  then $ \mathrm{K}\overset{\rm M}*F =\hat{F}(1/\cdot)$. 
  
  For $\alpha\in\mathbb{R}$, we denote by $\mathscr{W}_\alpha$ the space of all  non-negative  kernels $\mathrm{K}$ on $\mathbb{R}_+$ satisfying the \textit{Wiener condition}, i.e., for some $\delta>0$ and any $z\in\mathbb{R}$,
  \beqlb\label{eqn.App01}
  \int_0^\infty (t^{\alpha+\delta}\vee t^{\alpha-\delta} )\mathrm{K}(t)dt<\infty
  \quad \mbox{and}\quad
  \mathtt{k}_\alpha(z):=\int_0^\infty  t^{\alpha-\mathtt{i}z}\mathrm{K}(t)dt\notin \mathbb{C}\setminus\{ 0 \}. 
  \eeqlb 
  Here $\mathtt{i}$ is the imaginary unit.  If $\mathrm{K}(x)= e^{-x}$, then $\mathrm{K} \in \mathscr{W}_\alpha$ for any $\alpha>-1$ with 
  $\mathtt{k}_\alpha(0)= \Gamma(\alpha+1)$. Moreover,  for any $\kappa\in\mathbb{R}$, we put
  \begin{equation}
  	\begin{split}
  		\mathscr{M}^\infty_{\kappa,0}& := \left\{ F:\mathbb{R}_+ \to \mathbb{R} : \mbox{ $F$ is locally bounded, eventually positive or negative } \right. \\
  		& \qquad \qquad  \left. \mbox{and $  t^{-\kappa}\cdot F(t) $  is eventually monotone.} \right\} 
  	\end{split}
  \end{equation}	
  
  
  \begin{lemma}[The Wiener-Tauberian theorem in Pitt's form]\label{Thm.WienerPittTauberian}
  	For $\alpha\in\mathbb{R}$ and $\mathrm{K}\in \mathscr{W}_\alpha$,   the following hold. 
  	\begin{enumerate}
  		\item[(1)] If $F \in {\rm RV}_\alpha^\infty$, then $\mathrm{K}\overset{\rm M}*F \in {\rm RV}_\alpha^\infty$ and $\mathrm{K}\overset{\rm M}*F \sim \mathtt{k}_\alpha(0) \cdot F $ at infinity.
  		
  		\item[(2)] If $\mathrm{K}\overset{\rm M}*F \in {\rm RV}_\alpha^\infty$ and $F\in \mathscr{M}^\infty_{\kappa,0}$ for some $\kappa\in \mathbb{R}$, then  $ F \in {\rm RV}_\alpha^\infty$. 
  	\end{enumerate} 
  \end{lemma}
  \proof 
  To establish claim (1),  it suffices to prove that $ \mathrm{K}\overset{\rm M}*F(t) \sim   \mathtt{k}_\alpha(0) \cdot F(t) $ as $t\to\infty$. 
  Indeed, in this case, for each $x>0$, 
  \beqnn
  \lim_{t\to\infty} \frac{\mathrm{K}\overset{\rm M}*F(tx)}{\mathrm{K}\overset{\rm M}*F(t)} =    \lim_{t\to\infty} \frac{\mathrm{K}\overset{\rm M}*F(tx)}{F(tx)} \cdot  \lim_{t\to\infty} \frac{F(t)}{\mathrm{K}\overset{\rm M}*F(t)}  \cdot \lim_{t\to\infty} \frac{F(tx)}{F(t)} = x^\alpha.
  \eeqnn
  To prove that $\mathrm{K}\overset{\rm M}*F \sim \mathtt{k}_\alpha(0) \cdot F $, it suffices to show that as $t\to\infty$, 
  \beqlb\label{eqn.003}
  \Big|\frac{\mathrm{K}\overset{\rm M}*F(t)}{ F(t)}- \mathtt{k}_\alpha(0)  \Big| 
  = \Big|\int_0^\infty \mathrm{K}(y) \frac{F(ty)}{F(t)}dy-\int_0^\infty y^\alpha  \mathrm{K}(y)dy\Big| 
  \leq \int_0^\infty \mathrm{K}(y) \Big|\frac{F(ty)}{F(t)}- y^\alpha  \Big| dy \to 0.
  \eeqlb
  
  To this end, we fix $\varepsilon>0$. By Potter's theorem (Proposition~\ref{Thm.PotterThm}), there exists a constant $T_\varepsilon>0$ such that for any $t,ty\geq T_\varepsilon$, 
  \beqnn
  \Big|  \frac{F(ty)}{F(t)} -y^\alpha \Big| \leq \varepsilon (y^{\alpha+\delta}\vee y^{\alpha-\delta}). 
  \eeqnn 
  As a result,  as $\varepsilon \to 0$,
  \beqnn
  \int_{T_\varepsilon/t}^\infty \mathrm{K}(y)\Big| \frac{F(ty)}{F(t)}-y^\alpha\Big|  dy \leq  \varepsilon \cdot 	\int_0^\infty (z^{\alpha+\delta}\vee z^{\alpha-\delta} ) \mathrm{K}(z)dz \to 0, 
  \eeqnn 
  due to  (\ref{eqn.App01}). On the other hand, since $F$ is locally bounded, there exists a constant $C>0$ such that,
  \beqnn
  \int_0^{T_\varepsilon/t} \mathrm{K}(y)\Big| \frac{F(ty)}{F(t)}-y^\alpha\Big|  dy   \leq  \frac{C}{F(t)}  \int_0^{T_\varepsilon/t} \mathrm{K}(y) dy  + \int_0^{T_\varepsilon/t} y^\alpha \mathrm{K}(y) dy  . 
  \eeqnn
  
  Using (\ref{eqn.App01}) again, we see that the second term on the right of the above inequality tends to $0$ as $t\to\infty$. 
  To see that the first term vanishes as well, we distinguish two cases. If  $\alpha>0$, then $F(t)\to\infty$ as $t\to\infty$. Since $\mathrm{K}$ is locally integrable, we conclude that, as $t\to\infty$, 
  \beqnn
  \frac{C}{F(t)}  \int_0^{T_\varepsilon/t} \mathrm{K}(y) dy  \to 0.
  \eeqnn
  If $\alpha\leq 0$, then 
  \beqnn
  \int_0^{T_\varepsilon/t}  \mathrm{K}(y) dy  \leq   \int_0^{T_\varepsilon/t}  \big|ty/T_\varepsilon\big|^{\alpha-\delta}\cdot \mathrm{K}(y) dy \leq \big|t/T_\varepsilon\big|^{\alpha-\delta} \int_0^\infty y^{\alpha-\delta}\mathrm{K}(y)dy , 
  \eeqnn
  and hence 
  \beqnn
  \frac{C}{F(t)}  \int_0^{T_\varepsilon/t} \mathrm{K}(y) dy \leq C \cdot \frac{t^{\alpha-\delta}}{F(t)} ,
  \eeqnn
  which goes to $0$ as $t\to\infty$. 
  Putting these estimates together, we see that (\ref{eqn.003}) and hence (1) hold. 
  
  To prove claim (2) we assume w.l.o.g.~that $t^{-\kappa} F(t)$ is   non-increasing on $[t_0,\infty)$ for some $t_0\geq 0$. Since $\mathtt{k}_\alpha(0)>0$ we may also assume w.l.o.g.~that\footnote{Else, we may reformulate the theorem with the kernel $\mathrm{K}$ replaced by the rescaled kernel $\mathrm{K}(ct)$ for some constant $c>0$.} 
  \beqnn
  \int_{1/2}^1 s^{\kappa} \mathrm{K}(s)ds>0. 
  \eeqnn
  Moreover, we introduce the function $f_0(t):=f(t)\cdot \mathbf{1}_{\{t>2t_0\}}$ on $\mathbb{R}_+$ where 
  \[
  f(t) :=  \frac{F(t)}{\mathrm{K}\overset{\rm M}*F (t)}.
  \]
  If we can prove that $f_0$ is bounded and slowly increasing at infinity and that 
  \beqlb\label{eqn.302}
  \lim_{t\to\infty} \int_0^\infty s^\alpha \mathrm{K}(s)f_0(ts)ds = 1,
  \eeqlb
  then it follows from the Wiener-Pitt theorem (Proposition \ref{Wiener-Pitt.Thm}) with $k_0(t)=t^\alpha \mathrm{K}(t)$ and $g(t)= f_0(t)$ that
  \beqnn
  f_0(t)= f(t)=  \frac{F(t)}{\mathrm{K}\overset{\rm M}*F (t)} \to \frac{1}{\mathtt{k}_\alpha(0)},
  \eeqnn
  as $t\to \infty$ and hence that claim (2) holds. We proceed in two steps.
  
  {\bf Step 1.} {\sl The function $f_0$.}   To prove that $f$ is bounded on $[2t_0,\infty)$ and hence that $f_0$ is bounded, let $s\in[1/2,1]$ and $t\geq 2t_0$. Our monotonicity and positivity condition on $r^{-\kappa}F(r)$ imply that
  \beqnn
  (ts)^{-\kappa} F(ts)\geq t^{-\kappa} F(t) \qquad \mbox{and} \qquad F(ts)/F(t) > 0. 
  \eeqnn
  As a result, it holds  uniformly in $t\geq 2t_0$ that
  \beqnn
  \infty>\frac{1}{f(t)}:=  \frac{\mathrm{K}\overset{\rm M}*F (t)}{F(t)} \geq \int_{1/2}^1 \mathrm{K}(s) \frac{F(ts)}{F(t)}ds= \int_{1/2}^1 s^{\kappa}\mathrm{K}(s)   \frac{(ts)^{-\kappa} F(ts)}{t^{-\kappa} F(t)}ds \geq \int_{1/2}^1 s^{\kappa}\mathrm{K}(s) ds>0.
  \eeqnn
  
  To  prove that $f$ and hence $f_0$ is slowly increasing at infinity, we fix $u>1$ and $x\in[1,u]$. Then $(tx)^\kappa F(tx)\leq t^\kappa F(t)$, and we have for all $t\geq t_0$ that
  \beqnn
  f(t)- f(tx)
  \ar=\ar f(t)\cdot   \Big( 1- \frac{(tx)^{-\kappa} F(tx)}{t^{-\kappa} F(t)} \cdot  \frac{x^\kappa  \cdot \mathrm{K}\overset{\rm M}*F (t)}{   \mathrm{K}\overset{\rm M}*F (tx)}   \Big)\cr
  \ar\geq\ar 
  f(t)\Big( 1-  \frac{x^{\kappa} \cdot\mathrm{K}\overset{\rm M}*F (t)}{ \mathrm{K}\overset{\rm M}*F (tx)}  \Big) \cr
  \ar \geq \ar f(t)\Big( 1-  \frac{x^{\kappa} \cdot\mathrm{K}\overset{\rm M}*F (t)}{  \inf_{y\in[1,u]} \mathrm{K}\overset{\rm M}*F (ty)}  \Big).
  \eeqnn
  As a result, 
  \beqnn
  \inf_{x\in[1,u]} \big\{   f(t)- f(tx) \big\} 
  \ar\geq\ar f(t)\Big( 1-  \frac{(1\vee u^{\kappa} )\cdot \mathrm{K}\overset{\rm M}*F (t)}{ \inf_{y\in[1,u]} \mathrm{K}\overset{\rm M}*F (ty)}  \Big) \cr
  \ar  = \ar f(t)\Big( 1- (1\vee u^{\kappa} )\cdot \sup_{y\in[1,u]} \frac{\mathrm{K}\overset{\rm M}*F (t)}{  \mathrm{K}\overset{\rm M}*F (ty)}  \Big) . 
  \eeqnn
  The uniform convergence theorem for regularly varying functions given in Proposition~\ref{Thm.UniConver}, yields that
  \beqnn
  \lim_{t\to\infty }\Big( 1-  (1\vee u^{\kappa} )\cdot \sup_{y\in[1,u]} \frac{\mathrm{K}\overset{\rm M}*F (t)}{  \mathrm{K}\overset{\rm M}*F (ty)}  \Big)  =  1- (1\vee u^{\kappa} )\cdot(1\vee u^{-\alpha}),
  \eeqnn
  which goes to $0$ as $u\to 1+$. From this we deduce that 
  $$ \lim_{u\to 1+}\liminf_{t\to\infty} \inf_{x\in[1,u]} \{f(t) - f(tx)\} \geq 0$$ and
  hence that $f$ is slowly increasing at infinity.
  
  {\bf Step 2.}{\sl     The integrability condition.} 
  To prove the integrability condition (\ref{eqn.302}) we first rewrite our integral as
  \beqlb\label{eqn.30021}
  \int_0^\infty s^\alpha \mathrm{K}(s)f_0(ts)ds 
  = \int_0^\infty \mathrm{K}(s)f_0(ts)\Big( s^\alpha- \frac{\mathrm{K}\overset{\rm M}*F (ts)}{\mathrm{K}\overset{\rm M}*F (t)} \Big) ds
  +\int_0^\infty \mathrm{K}(s)f_0(ts) \frac{\mathrm{K}\overset{\rm M}*F (ts)}{\mathrm{K}\overset{\rm M}*F (t)} ds . 
  \eeqlb
  Similarly as in the proof of (\ref{eqn.003}), one can prove that the first term on the right-hand side of the above equality vanishes as $t\to\infty$. Using the definition of $f_0$ the second term can be rewritten as  
  \beqnn
  \int_0^\infty \mathrm{K}(s)f_0(ts) \frac{\mathrm{K}\overset{\rm M}*F (ts)}{\mathrm{K}\overset{\rm M}*F (t)} ds
  \ar=\ar  \int_{2t_0/t}^\infty \mathrm{K}(s)f(ts) \frac{\mathrm{K}\overset{\rm M}*F (ts)}{\mathrm{K}\overset{\rm M}*F (t)} ds 
  = \int_0^\infty  \frac{ \mathrm{K}(s) F(ts)}{\mathrm{K}\overset{\rm M}*F (t)} ds
  -\int_0^{2t_0/t}   \frac{\mathrm{K}(s)F(ts)}{\mathrm{K}\overset{\rm M}*F (t)} ds .
  \eeqnn
  
  The definition of $\mathrm{K}\overset{\rm M}*F $ implies that the first integral on the right side of the second equality equals one. 
  Moreover, by the local boundedness of $F$ there exists a constant $C>0$ such that for all $t\geq t_0$, 
  \beqnn
  \int_0^{2t_0/t}   \frac{\mathrm{K}(s)F(ts)}{\mathrm{K}\overset{\rm M}*F (t)} ds  \leq   \frac{C}{\mathrm{K}\overset{\rm M}*F (t)} \int_0^{2t_0/t}   \mathrm{K}(s) ds  . 
  \eeqnn
  
  If $\alpha>0$, the above quantity vanishes as $t\to\infty$, since $K$ is locally integrable and $\mathrm{K}\overset{\rm M}*F (t) \to \infty$.  
  If $\alpha\leq 0$, then it follows from (\ref{eqn.App01}) that there exists a constant $C>0$ such that for any $t>2t_0$, 
  \beqnn
  \int_0^{2t_0/t}   \mathrm{K}(s) ds  \leq  \int_0^{2t_0/t} \Big| \frac{ts} {2t_0}\Big|^{\alpha-\delta} K(s)ds \leq C \cdot t^{\alpha-\delta} 
  \eeqnn
  in which case
  \beqnn 
  \int_0^{2t_0/t}   \frac{\mathrm{K}(s)F(ts)}{\mathrm{K}\overset{\rm M}*F (t)} ds  \leq  \frac{C\cdot t^{\alpha-\delta} }{\mathrm{K}\overset{\rm M}*F (t)},
  \eeqnn
  which vanishes as $t\to \infty$,  since $\mathrm{K}\overset{\rm M}*F  \in {\rm RV}_\alpha^\infty$. 
  Putting these estimates together, we have  
  \beqnn 
  \int_0^\infty \mathrm{K}(s)f_0(ts) \frac{\mathrm{K}\overset{\rm M}*F (ts)}{\mathrm{K}\overset{\rm M}*F (t)} ds \to 1 .
  \eeqnn 
  \qed
  
  The following theorem extends the previous result to second-order regular varying functions.

  \begin{theorem}[Second-order Wiener-Tauberian theorem in Pitt's form]\label{Prop.App.B2}
  	For  $\alpha\in\mathbb{R}$, $\mathrm{K}\in \mathscr{W}_\alpha$ and $A\in \mathscr{A}_0^\infty$, the following hold.
  	\begin{enumerate}
  		\item[(1)] If $F\in {\rm 2RV}_{\alpha,0}^\infty( A)$, then $\mathrm{K}\overset{\rm M}*F  \in {\rm 2RV}_{\alpha,0}^\infty(A)$.
  		
  		\item[(2)] If $\mathrm{K}\overset{\rm M}*F  \in {\rm 2RV}_{\alpha,0}^\infty(A)$ and $F\in \mathscr{M}^\infty_{\alpha,0}$, then $F\in {\rm 2RV}_{\alpha,0}^\infty( A)$. 
  	\end{enumerate}
  	
  \end{theorem}
  \proof 
  In what follows it will be convenient to use the following functions: for $t>0$,   
  \beqnn
  F_{-\alpha}(t):= t^{-\alpha}\cdot F(t) \quad \mbox{and} \quad 
  \mathrm{K}_{\alpha}(t):= t^\alpha\cdot \mathrm{K}(t).
  \eeqnn
  
  To prove (1), we fix $F\in {\rm 2RV}_{\alpha,0}^\infty( A)$. It follows from Proposition~\ref{Prop.Fr} that $F_{-\alpha} \in  {\rm 2RV}_{0,0}^\infty(A)$ and hence $$F_{-\alpha} \in \Pi^\infty(F_{-\alpha}\cdot A),$$ due to Proposition~\ref{Prop.03}. Using the same arguments as in the proof of (4.11.5) in \cite[p.242]{BinghamGoldieTeugels1987}, one can hence show that 
  \begin{equation} \label{eq:K}
  	\mathrm{K}_{\alpha}\overset{\rm M}*F_{-\alpha} \in \Pi^\infty(\mathtt{k}_\alpha(0)\cdot F_{-\alpha}\cdot A).
  \end{equation}
  Thus, by Lemma~\ref{Thm.WienerPittTauberian}(1),
  \beqnn
  \mathrm{K}_{\alpha}\overset{\rm M}*F_{-\alpha}(t) \sim \mathtt{k}_\alpha(0)\cdot F_{-\alpha}(t),
  \eeqnn
  as $t\to\infty$. Moreover, the definition of the Mellin convolution yields, 
  \beqlb\label{eqn.005}
  \mathrm{K}_\alpha \overset{\rm M}*F_{-\alpha} (t)=  \int_0^\infty s^\alpha \mathrm{K} (s) (ts)^{-\alpha} F(ts)ds = t^{-\alpha}  \cdot \mathrm{K}\overset{\rm M}* F(t),\quad t> 0.
  \eeqlb
  In view of \eqref{Def.PiV.eqn} an application of the above results shows that, as $t\to\infty$, 
  \beqnn
  \frac{\mathrm{K}\overset{\rm M}*F (tx)/\mathrm{K}\overset{\rm M}*F (t)-x^\alpha}{A(t)}
  \ar=\ar x^\alpha \frac{\mathrm{K}_{\alpha}\overset{\rm M}*F_{-\alpha} (tx)-\mathrm{K}_{\alpha}\overset{\rm M}*F_{-\alpha} (t)}{\mathrm{K}_{\alpha}\overset{\rm M}*F_{-\alpha}(t) \cdot A(t)} \cr
  \ar\sim\ar  x^\alpha \frac{\mathrm{K}_{\alpha}\overset{\rm M}*F_{-\alpha} (tx)-\mathrm{K}_{\alpha}\overset{\rm M}*F_{-\alpha} (t)}{ \mathtt{k}_\alpha(0)\cdot F_{-\alpha}(t)\cdot A(t)} \cr
  \ar\ar\cr
  \ar \to \ar x^\alpha \log(x).
  \eeqnn 
  This shows that $\mathrm{K}\overset{\rm M}*F  \in {\rm 2RV}^\infty_{\alpha,0}(A)$.

  To prove (2) we first assume that $F_{-\alpha}$ is locally integrable and apply Proposition \ref{Prop.A.5} twice, first with $\overline F = \mathrm{K}_\alpha \overset{\rm M}*F_{-\alpha}$ and then again with $\overline F = F_{-\alpha}$ to show that $F_{-\alpha} \in {\rm 2RV}_{0,0}^\infty( A)$ and hence $F \in {\rm 2RV}_{\alpha,0}^\infty( A)$. 
  
  By (\ref{eqn.005}) and  Proposition~\ref{Prop.Fr}, we have $\mathrm{K}_\alpha \overset{\rm M}* F_{-\alpha} \in {\rm 2RV}^\infty_{0,0}(A)$. 
  Let $\mathcal{I}^{0,\uparrow}_{F_{-\alpha},1}$ be defined by (\ref{eqn.401}) with $t_0=0$, $F=F_{-\alpha}$ and $\theta=1$, and let
  %
  \beqnn
  G(t):=t^{-1}\cdot \mathcal{I}_{F_{-\alpha},1}^{0,\uparrow}(t), \qquad t>0.
  \eeqnn
  Taking the Mellin convolution of $G$ associated to the kernel $\mathrm{K}_\alpha$, we have that
  \beqnn
  \mathrm{K}_{\alpha}\overset{\rm M}*G(t)
  \ar=\ar \mathrm{K}_{\alpha}\overset{\rm M}*F_{-\alpha}(t) -\int_0^\infty \mathrm{K}_\alpha(s) \frac{1}{ts}\int_0^{ts}  F_{-\alpha}(r)drds\cr
  \ar =\ar \mathrm{K}_{\alpha}\overset{\rm M}*F_{-\alpha}(t)- \frac{1}{t}\int_0^t \int_0^\infty \mathrm{K}_{\alpha}(s) F_{-\alpha}(sr)ds ds \cr
  \ar=\ar  \mathrm{K}_{\alpha}\overset{\rm M}*F_{-\alpha}(t)- \frac{1}{t}\int_0^t 
  \mathrm{K}_{\alpha}\overset{\rm M}* F_{-\alpha}(s)ds \cr
  \ar\ar\cr
  \ar = \ar t^{-1}\cdot \mathcal{I}_{\mathrm{K}_{\alpha}\overset{\rm M}*F_{-\alpha},1}^{0,\uparrow}(t).
  \eeqnn
  Since $\mathrm{K}_{\alpha}\overset{\rm M}*F_{-\alpha} \in \Pi^\infty(\mathrm{K}_{\alpha}\overset{\rm M}*F_{-\alpha}\cdot A)$
  %
  we can deduce from Proposition~\ref{Prop.A.5} that
  \beqnn
  \mathrm{K}_{\alpha}\overset{\rm M}*G \in {\rm RV}_0^\infty
  \eeqnn
  and from \eqref{eq:K} that as $t\to\infty$, 
  \beqnn
  \mathrm{K}_{\alpha}\overset{\rm M}*G (t) \sim  \mathrm{K}_{\alpha}\overset{\rm M}*F_{-\alpha}(t)\cdot A(t) \sim \mathrm{k}_\alpha(0) F_{-\alpha}(t) \cdot A(t). 
  \eeqnn
  The tail behavior of the function $G$ can now be inferred from Lemma~\ref{Thm.WienerPittTauberian}(2). In fact, the function \beqnn
  t\cdot G(t)= \mathcal{I}_{F_{-\alpha},1}^{0,\uparrow}(t)
  \eeqnn
  is locally bounded and eventually positive (or negative, depending on $F$) and, as we show below, eventually monotone. Hence, it follows from Lemma~\ref{Thm.WienerPittTauberian}(2) that  
  \beqnn
  \mathtt{k}_\alpha(0) G(t)  \sim \mathrm{K}_{\alpha}\overset{\rm M}*G (t)   \sim \mathrm{k}_\alpha(0) F_{-\alpha}(t) \cdot A(t) ,
  \eeqnn
  as $t\to\infty$. In particular, 
  \beqnn
  G\sim F_{-\alpha} \cdot A \in {\rm RV}_{0}^\infty.
  \eeqnn
  To see that $\mathcal{I}_{F_{-\alpha},1}^{0,\uparrow}$ is eventually monotone, let $F_{-\alpha}$  be monotone on $[t_0,\infty)$ for some $t_0\geq 0$. For $t>0$, let $F^{t_0}_{-\alpha}(t):= F_{-\alpha}(t\vee t_0)$. 
  Both $F^{t_0}_{-\alpha}$ and $ \mathcal{I}_{F^{t_0}_{-\alpha},1}^{0,\uparrow}$ are monotone.  Moreover, for $t\geq t_0$, 
  \beqnn
  \mathcal{I}_{F_{-\alpha},1}^{0,\uparrow}(t)  
  \ar=\ar t F^{t_0}_{-\alpha}(t) - \int_0^tF^{t_0}_{-\alpha}(s)ds + \int_{0}^{t_0} [F_{-\alpha}(t_0)-F_{-\alpha}(s)]ds   \cr
  \ar=\ar  \mathcal{I}_{F^{t_0}_{-\alpha},1}^{0,\uparrow}(t) + \int_{0}^{t_0} [F_{-\alpha}(t_0)-F_{-\alpha}(s)]ds,
  \eeqnn
  and hence $ \mathcal{I}_{F_{-\alpha},1}^{0,\uparrow}$ is monotone on $[t_0,\infty)$. 
  
  Finally, it follows from (\ref{eqn.401}) that
  \beqnn
  G(t)= F_{-\alpha}(t) - \frac{1}{t} \int_0^t F_{-\alpha}(s)ds. 
  \eeqnn 
  Having shown that $G \in {\rm RV}_{0}^\infty$ we can again apply Proposition~\ref{Prop.A.5} to deduce that $F_{-\alpha}\in {\rm 2RV}^\infty_{0,0}(A)$. Hence, $F \in {\rm 2RV}^\infty_{\alpha,0}(A)$, due to Proposition~\ref{Prop.Fr}. 
  
  If $F_{-\alpha}$ is not locally integrable (i.e., $\alpha\geq 1$),  we consider two functions $\tilde{F}(t):= F(t)\cdot \mathbf{1}_{\{ t\geq t_0+1 \}}$ and   $\tilde{F}_{-\alpha}(t):= t^{-\alpha}\tilde{F}(t)$ on $\mathbb{R}_+$, which are locally integrable.  
  It suffices to prove that $\tilde{F} \in {\rm 2RV}^\infty_{\alpha,0}(A)$. 
  The preceding result shows that this holds if and only if $\mathrm{K}\overset{\rm M}* \tilde{F}\in {\rm 2RV}^\infty_{\alpha,0}(A)$.
  Notice that $  \mathrm{K}\overset{\rm M}* \tilde{F} =  \mathrm{K}\overset{\rm M}* F-\varepsilon$ with
  \beqnn
  \varepsilon(t):= \int_0^{(t_0+1)/t}  K(s)F(ts)ds, \quad t>0.
  \eeqnn
  Since $F$ is locally bounded and $\mathrm{K}$ is locally integrable, we have $ \varepsilon(t) \to0 $ and $ \mathrm{K}\overset{\rm M}* \tilde{F}(t) \sim  \mathrm{K}\overset{\rm M}* F(t) \in {\rm RV}^\infty_{\alpha}$ as $t\to\infty$. 
  Thus, for each $x>0$,
  \beqlb\label{eqn.004}
  \frac{\mathrm{K}\overset{\rm M}* \tilde{F}(tx)/\mathrm{K}\overset{\rm M}* \tilde{F}(t)-x^\alpha}{A(t)} 
  \ar=\ar \frac{\mathrm{K}\overset{\rm M}* \tilde{F}(tx)-x^\alpha\mathrm{K}\overset{\rm M}* \tilde{F}(t)}{\mathrm{K}\overset{\rm M}* \tilde{F}(t)A(t)}\cr
  \ar=\ar
  \frac{\mathrm{K}\overset{\rm M}* F(t)}{\mathrm{K}\overset{\rm M}* \tilde{F}(t)}\cdot \frac{\mathrm{K}\overset{\rm M}* F(tx)-x^\alpha \mathrm{K}\overset{\rm M}* F(t)  }{\mathrm{K}\overset{\rm M}* F(t) A(t)} 
  - \frac{\varepsilon(tx)-x^\alpha \varepsilon(t)}{\mathrm{K}\overset{\rm M}* \tilde{F}(t)A(t)},
  \eeqlb
  which tends to $ x^\alpha \log(x)$ as $t\to\infty$. This shows that $\mathrm{K}\overset{\rm M}* \tilde{F} \in {\rm 2RV}^\infty_{\alpha,0}(A)$. 
  \qed
  
  As a direct consequence of Theorem~\ref{Prop.App.B2}, the second-order version of Karamata's Tauberian theorem with zero second-order index can be obtained immediately; see the next theorem.

  \begin{theorem}[Second-order Karamata Tauberian theorem: $\rho=0$]\label{Thm.2KaraTaub02}
  	For any $\alpha>-1$ and $A \in \mathscr{A}^\infty_0$, the following hold.
  	\begin{enumerate}
  		\item[(1)] If $F\in {\rm 2RV}_{\alpha,0}^\infty (A)$,
  		then $\hat{F}(1/\cdot) \in {\rm 2RV}_{\alpha,0}^{\infty}( A)$.
  		
  		\item[(2)] If $\hat{F}(1/\cdot) \in {\rm 2RV}_{\alpha,0}^{\infty}( A)$ and $F\in \mathscr{M}^\infty_{\alpha,0}$, then $F\in {\rm 2RV}_{\alpha,0}^\infty (A)$.
  	\end{enumerate}   
  \end{theorem}
  
  
  The next corollary shows that the second-order regular variation of a function $F$ carries over to the Laplace-Stieltjes transforms of the functions $\mathcal{I}^{t_0,\uparrow}_{F,\,\cdot }$ and $\mathcal{I}^{\infty,\downarrow}_{F,\,\cdot}  $.   
  
  \begin{corollary}\label{Corollary.2KaramTauber}
  	For $\alpha \in \mathbb{R}$, $\rho\leq 0$ and $A\in\mathscr{A}^\infty_\rho$, the following hold.
  	\begin{enumerate}
  		\item[(1)] If $F\in {\rm 2RV}_{\alpha,\rho}^\infty(A)$, then  for each $t_0\geq 0$, $k_1,k_2\in\mathbb{Z}_+$ with $\alpha + k_1 +\rho > -2$ and $-1-\rho<\alpha + k_2<0$,
  		\begin{enumerate}
  			\item[(1.a)] $ \displaystyle{\widehat{\mathcal{I}^{t_0,\uparrow}_{F,k_1}}(1/\cdot) \in  {\rm 2RV}_{k_1+\alpha+1,\rho}^{\infty}\Big(\frac{\Gamma(k_1+\alpha+\rho +1)}{\Gamma(k_1+\alpha+1)}\frac{\alpha+\rho}{\alpha }\cdot A\Big)}$;
  			
  			\item[(1.b)] $ \displaystyle{\widehat{\mathcal{I}^{\infty,\downarrow}_{F,k_2}}(1/\cdot)  \in  {\rm 2RV}_{k_2+\alpha,\rho}^{\infty}\Big(\frac{\Gamma(k_2+\alpha+\rho )}{\Gamma(k_2+\alpha)}\frac{\alpha+\rho}{\alpha }\cdot A\Big)}$.
  		\end{enumerate}
  		
  		\item[(2)] We have $F\in {\rm 2RV}_{\alpha,\rho}^\infty(A)$  if one of the following holds: for some $t_0\geq 0$  and $k_1,k_2\in\mathbb{Z}_+$ defined in (1),
  		\begin{enumerate}
  			\item[(2.a)] Claim (1.a) holds and $\mathcal{I}^{t_0,\uparrow}_{F,k_1}\in \mathscr{M}^\infty_{k_1+\alpha+1,-}$ when $\rho<0$ or $\mathcal{I}^{t_0,\uparrow}_{F,k_1}\in \mathscr{M}^\infty_{k_1+\alpha+1,0}$ when $\rho=0$.
  			
  			\item[(2.b)] Claim (1.b) holds and $\mathcal{I}^{\infty,\downarrow}_{F,k_2}\in \mathscr{M}^\infty_{k_2+\alpha,-}$ when $\rho<0$ or $\mathcal{I}^{\infty,\downarrow}_{F,k_2}\in \mathscr{M}^\infty_{k_2+\alpha,0}$ when $\rho=0$.
  		\end{enumerate}
  	\end{enumerate}
  \end{corollary}
  \proof By Theorem~\ref{Thm.Gen2Kara}, we know that $F\in {\rm 2RV}_{\alpha,\rho}^\infty(A)$ if and only if  
  \beqnn
  \mathcal{I}^{t_0,\uparrow}_{F,k_1} \in {\rm 2RV}_{\alpha+k_1+1,\rho}^{\infty} \Big( \frac{(\alpha+\rho)(\alpha+k_1+1)}{\alpha(\alpha+\rho+k_1+1)} \cdot A \Big)
  \quad\mbox{or}\quad
  \mathcal{I}^{\infty,\downarrow}_{F,k_2} \in  {\rm 2RV}_{k_2+\alpha,\rho}^{\infty} \Big(\frac{(\alpha+\rho)(\alpha+k_2 )}{\alpha(\alpha+\rho+k_2)} \cdot A\Big).
  \eeqnn
  One can identify that conditions in Theorem~\ref{Thm.2KaraTaub01} and \ref{Thm.2KaraTaub02} are satisfied. 
  Hence the second-order regular variation of $\mathcal{I}^{t_0,\uparrow}_{F,k_1}$ and $ \mathcal{I}^{\infty,\downarrow}_{F,k_2} $ can be inherited by $\widehat{\mathcal{I}^{t_0,\uparrow}_{F,k_1}}(1/\cdot) $ and $\widehat{\mathcal{I}^{\infty,\downarrow}_{F,k_2}}(1/\cdot) $ respectively, and vice versa. 
  \qed

  
  \section{Second-order approximation of Hawkes processes} \label{Sec.HawkesP}
 \setcounter{equation}{0}
 
 In this section, we apply the previous results on second-order regular variation to establish second-order approximations for the mean and variance of Hawkes processes.
 
 Let $(\Omega,\mathscr{F},\mathbf{P})$ be a complete probability space endowed with a filtration $\{\mathscr{F}_t:t\geq 0\}$ satisfying the usual hypotheses.
 A Hawkes process $\{N(t):t\geq 0\}$ defined on this probability space is a random point process whose $(\mathscr{F}_t)$-intensity process $\{ \Lambda(t):t\geq 0 \}$ is of the form
 \beqlb\label{HawkesDensity}
 \Lambda(t)= \mu_0+ \sum_{0<\tau_i<t} \phi(t-\tau_i)=   \mu_0 + \int_0^t \phi(t-s)N(ds), \quad t\geq 0,
 \eeqlb
 where $\mu_0 > 0$ represents the \textit{exogenous density}, $\phi \in L^1 (\mathbb{R}+, \mathbb{R}+)$ is the \textit{kernel} or \textit{fertility/activation function}, and $\tau_i$ is the arrival time of the $i$-th event.  
 Based on the average number of child events
 \beqlb
 m := \|\phi\|_{L^1}= \int_0^\infty \phi(s)ds, 
 \eeqlb
 we distinguish three types of Hawkes processes. We call a Hawkes process \textit{subcritical} if $m < 1$, \textit{critical} if $m = 1$, and \textit{supercritical} if $m > 1$. 
 
 Our goal is to establish second-order approximations of the mean and variance of critical and subcritical Hawkes process. Specifically, we prove that  
 \beqnn
 \mathbf{E}[N(t)] \sim C_1 \cdot t^{\alpha_1} + \varepsilon_1(t)
 \quad\mbox{and}\quad 
 {\rm Var}(N(t)) \sim C_2 \cdot t^{\alpha_2} + \varepsilon_2(t), 
 \eeqnn
 where $\alpha_1,\alpha_2$ are two positive constants determined by the average number of child events and the release of self-excitation, and $ \varepsilon_1, \varepsilon_2$ are two regularly varying functions. We provide explicit representations of the functions $\varepsilon_1, \varepsilon_2$ from which we deduce that 
 \beqnn
 \varepsilon_i(t) \sim t^{\alpha_i - \eta} A_i(t), 
 \eeqnn
 as $t \to \infty$ for some slowly varying function $A_i$ that converges to zero and some coefficient $\eta \in [0,\alpha_i]$. 
 
 Depending on the asymptotics of the kernel $\phi$ this coefficient may be very close to zero or even equal to zero, in which case the first-order approximations  $C_i \cdot t^{\alpha_i}$ may only provide poor approximations of the mean and variance, except in the very long run.

 \subsection{Preliminaries}
 
 Associated with the kernel $\phi$ we define the resolvent $R$ by the unique solution of the \textit{resolvent equation}
 \beqlb\label{Resolvent}
 R(t)= \phi(t)+ R*\phi(t),\quad t\geq 0. 
 \eeqlb
 The function $R$  is also known as the \textit{renewal function} in renewal theory. 
 It combines the direct and indirect impact of an external event on the arrival of future events. 
 Properties of $(N, \Lambda)$ are often formulated in terms of $R$. We will also need the following three functions:
 \beqnn
 \mathcal{I}_\Lambda(t):= \int_0^t \Lambda(s)ds, \quad 
 \mathcal{I}_R(t):= \int_0^t R(s)ds
 \quad\mbox{and}\quad
 \mathcal{I}_R^2(t):= \int_0^t \mathcal{I}_R(s) ds, \quad t\geq 0.
 \eeqnn 
 
 We notice that Hawkes process $N$ has compensator $\mathcal{I}_\Lambda$ and that the compensated point process $\widetilde{N}:=N-\mathcal{I}_\Lambda$ is an $(\mathscr{F}_t)$-martingale.   
 The following martingale representation of $\Lambda$ was first introduced by Bacry et al. \cite{BacryDelattreHoffmannMuzy2013}.  
 
 \begin{lemma}[Martingale representation]\label{MartRep}
 	The intensity process $\Lambda$ is the unique solution to
 	\beqlb\label{SVR}
 	\Lambda(t) \ar=\ar \mu_0  \big(1+ \mathcal{I}_R(t) \big) + \int_0^t R(t-s)\widetilde{N}(ds) ,\quad t\geq 0.
 	\eeqlb  
 \end{lemma}

 The martingale representation theorem provides an exact representation of the mean and variance of the Hawkes process $N$ in term of the resolvent $R$  as shown by the following corollary. 
 
 \begin{corollary}\label{Coro.401}
 	For the Hawkes process $(N, \Lambda)$, we have for $t\geq 0$, 
 	\begin{equation}\label{eqn.3007} 
 		\begin{split}
 			\mathbf{E} \big[N(t) \big] & = \mu_0 \big(t+ \mathcal{I}^2_R(t)\big),  \\
 			{\rm Var} \big(N(t) \big) 
 			&  =\mu_0\big( t+ 3 \cdot \mathcal{I}^2_R (t) + 2 \cdot \mathcal{I}_R*\mathcal{I}_R(t) + |\mathcal{I}_R|^2* \mathcal{I}_R(t) + \mathcal{I}_{|\mathcal{I}_R|^2}(t)  \big). 
 		\end{split}
 	\end{equation}
 \end{corollary}
 \proof Integrating both side of (\ref{SVR}) over $(0,t]$ and then using Fubini's theorem as well as  the stochastic Fubini theorem; see \cite[Theorem~D.2]{Xu2021b} or \cite[Theorem~2.6]{Walsh1986}, yields  
 \beqlb\label{eqn.3006}
 \mathcal{I}_\Lambda(t) = \mu_0  \big(t+ \mathcal{I}_R^2(t) \big) + \int_0^t R(t-s) \widetilde{N}(s)ds, \quad t\geq 0.
 \eeqlb
 Taking expectations on both sides of this equality,   we have 
 $  \mathbf{E}[N(t)]= \mathbf{E}[\mathcal{I}_\Lambda(t)]   = \mu_0  \big(t+ \mathcal{I}_R^2(t) \big) $. 
 By the perfect square trinomial,
 \beqnn
 \mathbf{E} \big[|N(t)|^2 \big] = \mathbf{E} \big[|\widetilde{N}(t)+\mathcal{I}_\Lambda(t)|^2 \big] 
 = \mathbf{E} \big[| \widetilde{N}(t) |^2 \big]   + 2\mathbf{E} \big[\widetilde{N}(t)\mathcal{I}_\Lambda(t)  \big] + \mathbf{E} \big[|\mathcal{I}_\Lambda(t)|^2 \big] .
 \eeqnn
 The martingality of $\widetilde{N}$ yields that 
 \[
 \mathbf{E} \big[\widetilde{N}(t)\widetilde{N}(r) \big]= \mathbf{E} \big[|\widetilde{N}(t\wedge r)|^2 \big] =\mathbf{E} \big[N(t\wedge r) \big] =\mu_0  \big(t\wedge r+ \mathcal{I}_R^2(t\wedge r) \big),
 \]
 for any $t,r\geq 0$. By this and Fubini's theorem, 
 \beqnn
 \mathbf{E} \big[\widetilde{N}(t)\mathcal{I}_\Lambda(t)  \big]  
 \ar=\ar \int_0^t R(t-s) \mathbf{E} \big[\widetilde{N}(t)\widetilde{N}(s) \big]ds  \cr
 \ar =\ar \int_0^t R(t-s) \mathbf{E} \big[|\widetilde{N}(s)|^2 \big]ds 
 =  \mu_0  \int_0^t R(t-s) \big(s+ \mathcal{I}_R^2(s) \big)ds 
 = \mu_0 \big( \mathcal{I}_R^2(t) + \mathcal{I}_R*\mathcal{I}_R(t) \big). 
 \eeqnn
 Squaring both sides of (\ref{eqn.3006}) and then taking expectations,
 \beqnn
 \mathbf{E} \big[|\mathcal{I}_\Lambda(t)|^2 \big] 
 \ar=\ar |\mu_0|^2  \big(t+ \mathcal{I}_R^2(t) \big)^2  
 + \int_0^t ds \int_0^t R(t-s) R(t-r) \mathbf{E}\big[\widetilde{N}(s)\widetilde{N}(r) \big]dr   \cr
 \ar=\ar |\mu_0|^2  \big(t+ \mathcal{I}_R^2(t) \big)^2  
 + \mu_0\int_0^tds\int_0^t R(t-r) R(t-s)  \big(s\wedge r+ \mathcal{I}_R^2(s\wedge r) \big) dr  \cr
 \ar=\ar |\mu_0|^2  \big(t+ \mathcal{I}_R^2(t) \big)^2  
 + 2\mu_0\int_0^t R(t-s)  \big(s+ \mathcal{I}_R^2(s) \big)  ds \int_s^t R(t-r) dr  \cr
 \ar\ar\cr
 \ar=\ar  |\mu_0|^2  \big(t+ \mathcal{I}_R^2(t) \big)^2    
 + \mu_0  \cdot \big( \mathcal{I}_{|\mathcal{I}_R|^2}(t) + |\mathcal{I}_R|^2* \mathcal{I}_R(t) \big) .  
 \eeqnn
 Putting the above estimates together and then plugging them into ${\rm Var} \big(N(t) \big)= \mathbf{E} \big[|N(t)|^2 \big]-  \big|\mathbf{E}[N(t)] \big|^2$,  yields the representation (\ref{eqn.3007}).
 \qed

 The above corollary shows that the long-term behavior of $\mathbf{E}\big[N(t)\big]$ and ${\rm Var}\big(N(t)\big)$ is fully determined by the tail behaviour of the functions $\mathcal{I}_R$ and $\mathcal{I}_R^2$ at infinity. The results that follow show that the tail behaviour of these functions heavily relies on the tail behaviour of the functions
 \beqnn
 {\it\Phi}(t):=\int_t^\infty \phi(s)ds
 \quad \mbox{and}\quad
 {\it \Psi}_\beta(t):= \int_0^t s^\beta \phi(s)ds, \quad
 t,\beta \geq 0.
 \eeqnn 
 Since ${\it\Phi}$ is non-increasing with ${\it\Phi}(0)=m<\infty$ it satisfies our assumption \ref{H}. 
 Their respective Laplace-Stieltjes transforms ${\it \hat\Phi}$ and ${\it \hat\Psi}_\beta$ are given by
 \beqlb
 {\it\hat\Phi}(\lambda)
 \ar:=\ar  \int_0^\infty \lambda e^{-\lambda t} {\it\Phi}(t)dt = \int_0^\infty \big(1-e^{-\lambda t}\big)\phi(t)dt,\label{LapPhi}\\
 {\it\hat\Psi}_\beta(\lambda)
 \ar:=\ar \int_0^\infty \lambda e^{-\lambda t} {\it\Psi}_\beta(t)dt = \int_0^\infty e^{-\lambda t}t^\beta\phi(t)dt,\label{LapPsi}
 \eeqlb
 for $\lambda> 0$.
 In particular, the function $ {\it \Psi}_0$ corresponds to the integrated function of $\phi$, and  $m = {\it\Phi}(0) = {\it\Psi}_0(\infty)$. 
 
 It has been shown in  \cite{HorstXu2023}  that the long-run behavior of critical Hawkes processes  $(m=1)$ strongly depends on the finiteness of the quantity 
 \beqnn
 \sigma := \Psi_1(\infty) =\int_0^\infty  {\it\Phi}(s) ds \in (0, \infty].
 \eeqnn
 Following  \cite{HorstXu2023}, we call a critical Hawkes process $N$ \textit{strongly critical} if $\sigma=\infty$  and \textit{weakly critical} if $\sigma < \infty$. In the weakly critical case, the mean-value theorem implies that as $\lambda \to \infty$, 
 \beqlb\label{sigma1}
 {\it\hat\Phi}(1/\lambda) \sim \frac \sigma \lambda.
 \eeqlb

 \subsection{First-order approximation} 
 
 In this section, we establish first-order approximations of the mean and the variance of a Hawkes process in terms of the first-order approximations of the functions $\mathcal{I}_R$ and $\mathcal{I}_R^2$. Their Laplace-Stieltjes transforms are given by 
 \beqlb\label{eqn.LapIR}
 \hat{\mathcal{I}}_R(\lambda)= \frac{m-{\it\hat\Phi}(\lambda)}{1-m+ {\it\hat\Phi}(\lambda)}
 \quad \mbox{and}\quad
 \hat{\mathcal{I}}_R^2(\lambda)= \frac{\hat{\mathcal{I}}_R(\lambda)}{\lambda} ,\quad \lambda >0.
 \eeqlb  
 
 \begin{proposition}\label{Thm.AsymR}
 	Three regimes arise for the long-term behavior of $\mathcal{I}_R$ and $\mathcal{I}_R^2$.
 	\begin{enumerate}
 		\item[(1)] When $m<1$, we have as $t\to\infty$, 
 		\beqnn
 		\mathcal{I}_R(t)\to \frac{m}{1-m}  
 		\quad \mbox{and}\quad 
 		\mathcal{I}_R^2(t)\sim \frac{m\cdot t}{1-m}.
 		\eeqnn

 		\item[(2)] When $m=1$ and $\sigma <\infty$,  we have as $t\to\infty$, 
 		\beqnn
 		\mathcal{I}_R(t)\to \frac{t}{\sigma}  
 		\quad \mbox{and}\quad 
 		\mathcal{I}_R^2(t)\sim \frac{t^2}{2\sigma}.
 		\eeqnn

 		\item[(3)] When $m=1$ and $\sigma =\infty$, if ${\it \Phi}\in{\rm RV}_{-\alpha}^\infty$ for some\footnote{Note that $\sigma <\infty$ when ${\it \Phi}\in{\rm RV}_{-\alpha}^\infty$ with $\alpha>1$.} $\alpha\in[0,1]$,  we have as $t\to\infty$, 
 		\beqnn
 		\mathcal{I}_R(t)   \sim  
 		\begin{cases}
 			\displaystyle{ \frac{1/{\it\Phi}(t)}{\Gamma(1-\alpha)\Gamma(1+\alpha)}  \in {\rm RV}_{\alpha }^\infty,} & \mbox{if }\alpha\in[0,1);\vspace{5pt} \\
 			\displaystyle{\frac{t}{{\it\Psi}_1(t)}   \in {\rm RV}_{ 1}^\infty,} & \mbox{if }\alpha=1,
 		\end{cases} 
 		\quad \mbox{and}\quad
 		\mathcal{I}_R^2(t) 
 		\sim \frac{t\cdot 	\mathcal{I}_R(t)}{1+\alpha}   \in {\rm RV}_{\alpha +1}^\infty.
 		\eeqnn 
 	\end{enumerate} 
 \end{proposition}
 \proof 
 Using  integration by parts and then the first equality in (\ref{eqn.LapIR}), 
 \beqlb\label{LapResolvent}
 \int_0^\infty e^{-t/\lambda } R(t)dt =\hat{\mathcal{I}}_R (1/\lambda) =\frac{m-{\it\hat\Phi}(1/\lambda)}{1-m+{\it\hat\Phi}(1/\lambda)},  \quad \lambda >0. 
 \eeqlb
 \begin{enumerate}
 	\item[(1)] 
 	Let $m<1$. Since ${\it\hat\Phi}(1/\lambda)\to0$ as  $\lambda \to \infty$, we have  as $t\to\infty$,
 	\beqlb\label{eqn.50001}
 	\mathcal{I}_R(t) \to \mathcal{I}_R(\infty)= \lim_{\lambda \to \infty} \hat{\mathcal{I}}_R (1/\lambda) = \frac{m}{1-m} ,
 	\eeqlb
 	which shows that $\mathcal{I}_R \in {\rm RV}_0^\infty$. By Proposition~\ref{Thm.Karamata}(1), as $t\to\infty$,
 	\beqnn
 	\mathcal{I}_R^2 (t) = \int_0^t \mathcal{I}_R(s)ds \sim  t \cdot \mathcal{I}_R(t) \sim  \frac{m\cdot t}{1-m}. 
 	\eeqnn
 	as $t\to\infty$. Hence claim (1) holds. 
 	
 	\item[(2)] Let $m=1$ and $\sigma < \infty$. In view of (\ref{LapResolvent}) and \eqref{sigma1} we have as $\lambda \to \infty$,
 	\beqlb\label{eqn.30022}
 	\hat{\mathcal{I}}_R (1/\lambda) \sim \frac{1}{ {\it\hat\Phi}(1/\lambda) } \sim \frac{\lambda}{\sigma} \in {\rm RV}^\infty_1.
 	\eeqlb
 	By Karamata's Tauberian theorem (Proposition~\ref{Thm.KaramataTauberian}) this shows that $$\mathcal{I}_R(t) \sim \frac{t}{ \sigma},$$ as $t\to\infty.$ 
 	From this and Proposition~\ref{Thm.Karamata}(1), the assertion follows from
 	\beqnn 
 	\mathcal{I}^2_R(t) =  \int_0^t \mathcal{I}_R(s)ds \sim  \frac{t^2}{2\sigma}.
 	\eeqnn
 	\item[(3)] Let $m=1$ and $\sigma=\infty$. In this case, we claim that  ${\it\hat\Phi}(1/\cdot) \in{\rm RV}_{-\alpha}^\infty$
 	and as $\lambda \to \infty$, 
 	\beqlb\label{eqn.1001}
 	{\it\hat\Phi}(1/\lambda)  \sim  
 	\begin{cases}
 		\Gamma(1-\alpha)  {\it \Phi}(\lambda) \in {\rm RV}^\infty_{-\alpha}, & \mbox{if }\alpha\in[0,1); \vspace{5pt}\\
 		\lambda^{-1}\cdot	{\it\Psi}_1(\lambda) \in {\rm RV}^\infty_{0}, & \mbox{if }\alpha=1. 
 	\end{cases} 
 	\eeqlb 
 	In fact, since ${\it\Phi}$ is a tail function of a probability distribution, equation \eqref{eqn.1001} for $\alpha \in [0,1)$ follows from Proposition~\ref{Prop.A.7} as condition $(2.i)$ therein is satisfied.  
 	If $\alpha = 1$, by using integration by parts and then Proposition~\ref{Thm.Karamata}(1), 
 	\beqlb\label{eqn.50002}
 	{\it \Psi}_1(t) = -t{\it\Phi}(t) + \int_0^t {\it \Phi}(s)ds = \Big( 1- \frac{t{\it\Phi}(t) }{\int_0^t {\it \Phi}(s)ds} \Big) \cdot \int_0^t {\it \Phi}(s)ds \sim  \int_0^t {\it \Phi}(s)ds \in {\rm RV}^\infty_{0}
 	\eeqlb
 	as $t\to\infty$. Thus, equation \eqref{eqn.1001} for $\alpha = 1$ it follows again from Proposition~\ref{Prop.A.7} as condition $(2.ii)$ therein holds. 
 	
 	Taking (\ref{eqn.1001}) back into the first asymptotic equivalence in (\ref{eqn.30022}) shows that 
 	$\hat{\mathcal{I}}_R (1/\cdot)   \in {\rm RV}^\infty_{\alpha}$. 
 	Using Karamata's Tauberian theorem (Proposition~\ref{Thm.KaramataTauberian}) and then Proposition~\ref{Thm.Karamata}(1) this shows that, as $t\to\infty$,
 	\beqnn
 	\mathcal{I}_R (t) \sim  \frac{1/{\it\hat\Phi}(1/t)}{\Gamma(1+\alpha)}
 	\quad \mbox{and}\quad 
 	\mathcal{I}_R^2(t) 
 	\sim \frac{t\cdot \mathcal{I}_R(t)}{1+\alpha} . 
 	\eeqnn 
 	Assertion (3) now follows from (\ref{eqn.1001}). 
 \end{enumerate}
 \qed
 
 The next theorem establishes first-order approximations of $\mathbf{E}[N(t)]$ and ${\rm Var}(N(t))$ by taking the asymptotic results obtained in Proposition~\ref{Thm.AsymR} back into Corollary~\ref{Coro.401}. 
 
 \begin{theorem}
 	The following regimes arise for the long-term behavior of $\mathbf{E}[N(t)]$ and ${\rm Var}(N(t))$.
 	\begin{enumerate}
 		\item[(1)] (Subcritical case) When $m<1$, we have as $t\to\infty$,
 		\beqnn
 		\mathbf{E}[N(t)] \sim \frac{\mu_0\cdot t }{1-m}
 		\quad \mbox{and}\quad 
 		{\rm Var} \big(N(t)\big) \sim  \frac{\mu_0\cdot t }{(1-m)^3} . 
 		\eeqnn
 		
 		\item[(2)] (Weakly critical case) When $m=1$ and $\sigma <\infty$, we have as $t\to\infty$,
 		\beqnn
 		\mathbf{E}[N(t)] \sim \frac{\mu_0 }{2\sigma} \cdot t^2
 		\quad \mbox{and}\quad 
 		{\rm Var} \big(N(t)\big) \sim  \frac{\mu_0 }{12\sigma^3} \cdot t^4 . 
 		\eeqnn
 		
 		\item[(3)] (Strongly critical case) When $m=1$ and $\sigma =\infty$, we have as $t\to\infty$,
 		\begin{enumerate}
 			\item[(3.a)] if ${\it \Phi}\in{\rm RV}_{-\alpha}^\infty$ with $\alpha\in[0,1)$, then
 			\beqnn
 			\mathbf{E}[N(t)]  \sim \frac{\mu_0 }{\Gamma(1-\alpha)\Gamma(2+\alpha)} \cdot \frac{t}{{\it\Phi}(t)} 
 			\quad \mbox{and}\quad 
 			{\rm Var} \big(N(t)\big) \sim 
 			\frac{\mu_0 \mathrm{B}(2\alpha+1,\alpha+1)}{|\Gamma(1-\alpha)\Gamma(1+\alpha)|^3} \cdot \frac{t}{|{\it\Phi}(t)|^3};
 			\eeqnn  
 			
 			\item[(3.b)] if ${\it \Phi}\in{\rm RV}_{-1}^\infty$, then  
 			\beqnn
 			\mathbf{E}[N(t)]  \sim \frac{\mu_0}{2}\cdot \frac{ t^2}{{\it\Psi}_1(t)} 
 			\quad \mbox{and}\quad 
 			{\rm Var} \big(N(t)\big) \sim \frac{\mu_0}{12}\cdot	\frac{t^4}{ |{\it\Psi}_1(t)  |^3}  .
 			\eeqnn
 		\end{enumerate}
 	\end{enumerate}
 	
 \end{theorem}
 \proof 
 We only prove claim (3). The other two claims can be established in the same way. 
 Since $\mathcal{I}_R(t) \to \infty$ and hence $\mathcal{I}^2_R(t) / t \to\infty$ as $t\to\infty$, it follows from Corollary~\ref{Coro.401} that  
 $$\mathbf{E}[N(t)]= \mu_0 \big(t+ \mathcal{I}^2_R(t)\big) \sim \mu_0\cdot 	\mathcal{I}^2_R(t), $$
 and the result for $\mathbf{E}[N(t)]$ follows from Proposition~\ref{Thm.AsymR}(3). To establish the approximation of the variance, we use the fact that $\mathcal{I}_R \in {\rm RV}^\infty_{\alpha}$ and first apply Proposition~\ref{Thm.Karamata}(1) and then Corollary~\ref{Appendix.Coro01} to get that
 \begin{align*}
 	\mathcal{I}_{|\mathcal{I}_R |^2}(t)  & \sim  (2\alpha+1)^{-1}\cdot t\cdot |\mathcal{I}_R(t)|^2,\\
 	\mathcal{I}_R * \mathcal{I}_R(t) &  \sim  \mathrm{B}(\alpha+1,\alpha+1)\cdot t \cdot |\mathcal{I}_R(t)|^2, \\ 
 	|\mathcal{I}_R|^2 * \mathcal{I}_R(t) & \sim  \mathrm{B}(2\alpha+1,\alpha+1)\cdot t \cdot |\mathcal{I}_R(t)|^3. 
 \end{align*}
 Taking these and the fact that $\mathcal{I}_R^2(t) \sim (1+\alpha)^{-1} \cdot t\cdot 	\mathcal{I}_R(t)$ into (\ref{eqn.3007}), we have as $t\to\infty$, 
 \beqnn
 {\rm Var}\big( N(t) \big) \sim \mu_0 \cdot  \mathrm{B}(2\alpha+1,\alpha+1)\cdot t \cdot |\mathcal{I}_R(t)|^3,
 \eeqnn 
 and the desired asymptotic results for ${\rm Var} \big(N(t)\big)$ follow from Proposition~\ref{Thm.AsymR}(3). 
 \qed

 \subsection{Second-order approximation} 
 
 We proceed to prove second-order approximations of $\mathbf{E}\big[ N(t) \big]$ and ${\rm Var}\big( N(t) \big)$ under first- and second-order regular variation conditions on the function ${\it\Phi}$.  
 As preparation, the next three propositions establish second-order approximations of the functions $\mathcal{I}_R$ and $\mathcal{I}_R^2$ for subcrititical, weakly critical and strongly critical Hawkes processes, respectively.  
 
 \begin{proposition}\label{Prop.2Est.01}
 	If $m<1$ and ${\it \Phi}\in {\rm RV}^\infty_{-\alpha}$ with $\alpha \geq 0$, we have   as $t\to\infty$, 
 	\beqlb\label{eqn.40010}
 	\frac{m}{1-m} -\mathcal{I}_R(t) \sim \frac{ {\it\Phi}(t)}{(1-m)^2} 
 	\quad \mbox{and} \quad
 	\frac{m \cdot t}{1-m}-\mathcal{I}_R^2 (t) \sim 
 	\begin{cases}
 		\displaystyle{ \frac{t{\it \Phi}(t)}{(1-m)^2(1-\alpha)} } \in {\rm RV}^\infty_{1-\alpha}, & \mbox{if $\alpha\in[0,1)$}; \vspace{5pt}\\
 		\displaystyle{ \frac{{\it \Psi}_1(t)}{(1-m)^2 } } \in {\rm RV}^\infty_{0}, & \mbox{if $\alpha\geq 1$}.
 	\end{cases} 
 	\eeqlb
 \end{proposition}
 \proof 
 To prove the first result, let $G$ be a probability distribution on $\mathbb{R}_+$ with density function $m^{-1}\cdot \phi$. 
 For each $n\geq 0$, let $G^{*n}$ be the $n$-fold convolution of $G$ defined by
 \beqnn
 G^{*0}(t)\equiv 1 , \quad G^{*1}(t)= G(t) \quad \mbox{and}\quad 
 G^{*n} (t)= \int_0^t G^{*(n-1)} (t-s) dG(s) = m^{-n}\int_0^t \phi^{*n}(s)ds, \quad t\geq 0.
 \eeqnn
 Since the resolvent $R$ admits the Neumann series expansion
 \beqlb\label{eqn.4001}
 R(t)= \sum_{n=1}^\infty \phi^{*n}(t),\quad t\geq 0; 
 \eeqlb
 see for instance (1.4) in \cite[p.37]{GripenbergLondenStaffans1990}.
 we see that the function $\mathcal{I}_R(t) $ admits the represenattion
 \beqnn
 \mathcal{I}_R(t) =  \sum_{n=1}^\infty \int_0^t \phi^{*n}(s)ds = \sum_{n=1}^\infty m^{n} \cdot G^{*n} (t), \quad t\geq 0.
 \eeqnn
 
 Since the tail distribution function $1-G(t) = m^{-1}\cdot {\it\Phi}(t)$ is slowly varying at infinity, it follows from Lemma~1 in \cite{ChoverNeyWainger1973} with $\gamma =m$ and $d=1$ that the function $U_m(t):= 1+  \mathcal{I}_R(t)$ satisfies 
 \beqnn 
 \frac{(1-m)^{-1} -U_m(t) }{1-G(t)} \sim \frac{m}{(1-m)^2},
 \eeqnn
 as $t\to\infty$. As a result,  
 \beqlb\label{eqn.10021}
 \frac{m}{1-m} - \mathcal{I}_R(t) \sim \frac{ {\it\Phi}(t)}{(1-m)^2} \in {\rm RV}^\infty_{-\alpha} .  
 \eeqlb 
 This proves the first result. To prove the second result, we use the representation
 \beqlb\label{eqn.1002}
 \frac{m \cdot t}{1-m}-	\mathcal{I}_R^2 (t) = \int_0^t \Big(   \frac{m}{1-m} - \mathcal{I}_R(s)\Big)ds
 \eeqlb
 and distinguish two cases. 
 
 \begin{itemize}
 	\item If $\alpha \in [0,1)$ or if $\alpha=1$ and $\sigma={\it\Psi}_1(\infty)=\infty$, then it follows from (\ref{eqn.50002}) that 
 	\beqnn
 	\int_0^t {\it\Phi}(s)ds \to \infty, 
 	\eeqnn
 	as $t\to\infty$. Hence, it follows from \eqref{eqn.10021} and\eqref{eqn.1002} that
 	\beqnn
 	\frac{m \cdot t}{1-m}-	\mathcal{I}_R^2 (t) 
 	\sim \frac{1}{(1-m)^2} \cdot  \int_0^t {\it\Phi}(s)ds 
 	\sim \begin{cases}
 		\displaystyle{ \frac{t{\it \Phi}(t)}{(1-m)^2(1-\alpha)} }, & \mbox{if $\alpha\in[0,1)$}; \vspace{5pt}\\
 		\displaystyle{ \frac{{\it \Psi}_1(t)}{(1-m)^2 } }, & \mbox{if $\alpha=1$ and $\sigma=\infty$} 
 	\end{cases}  
 	\eeqnn
 	as $t\to\infty$, where the second asymptotic equivalence follows from Proposition~\ref{Thm.KaramataTauberian}(1) in case $\alpha\in[0,1)$ and from (\ref{eqn.50002}) in case $\alpha=1$ and $\sigma=\infty$. 
 	
 	\item  Let $\alpha>1$ or $\alpha=1$ and $\sigma={\it\Psi}_1(\infty)<\infty$. Since $\mathcal{I}_R(\infty) = m/(1-m)$, integration by parts shows that  
 	\beqlb\label{eqn.1003}
 	\int_0^\infty \Big(   \frac{m}{1-m} - \mathcal{I}_R(t)\Big)dt 
 	\ar=\ar \int_0^\infty \int_t^\infty R(s)ds dt =\int_0^\infty tR(t)dt  . 
 	\eeqlb
 	Differentiating both sides of the first equality in (\ref{eqn.LapIR}) with respect to $\lambda $, and then setting $\lambda =0$ it follows that 
 	\begin{equation}\label{eqn.10004}
 		\begin{split}
 			\int_0^\infty tR(t)dt  =  
 			\frac{d}{d\lambda} \hat{\mathcal{I}}_{R}(0) 
 			=  \frac{-\frac{d}{d\lambda}  {\it\hat\Phi}(0)}{(1-m+  {\it\hat\Phi}(0))^2} 
 			=  \frac{\int_0^\infty s \phi(s)ds}{(1-m)^2} 
 			=  \frac{{\it\Psi}_1(\infty)}{(1-m)^2 }. 
 		\end{split}
 	\end{equation}
 	
 	Since ${\it\Psi}_1(\infty) < \infty$ taking \eqref{eqn.1003} and \eqref{eqn.1004} back into \eqref{eqn.1002} we deduce that
 	\beqnn
 	\frac{m \cdot t}{1-m}-	\mathcal{I}_R^2 (t) 
 	\ar = \ar  \int_0^t \Big(   \frac{m}{1-m} - \mathcal{I}_R(t)\Big)dt = \frac{{\it\Psi}_1(\infty)}{(1-m)^2}- \int_t^\infty \Big(   \frac{m}{1-m} - \mathcal{I}_R(s)\Big)ds.
 	\eeqnn
 	
 	Hence, it remains to prove that the integral term on the right hand side of the second equation vanishes. This follows from
 	%
 	\beqnn
 	\int_t^\infty \Big(   \frac{m}{1-m} - \mathcal{I}_R(s)\Big)ds 
 	\sim  \frac{1}{(1-m)^2} \cdot  \int_t^\infty {\it\Phi}(s)ds 
 	\sim \frac{ t {\it\Phi}(t)}{(1-m)^2} 
 	\eeqnn 
 	and 
 	\beqnn
 	t {\it\Phi}(t) = \int_t^\infty t \phi(s)ds \leq \int_t^\infty s \phi(s)ds = {\it\Psi}_1(\infty)-{\it\Psi}_1(t) \to 0,
 	\eeqnn
 	as $t\to\infty$ to conclude that,  
 	\beqnn
 	\frac{m \cdot t}{1-m}-	\mathcal{I}_R^2 (t)  \sim \frac{{\it\Psi}_1(t)  }{(1-m)^2}, \quad \mbox{as $t \to \infty$.}
 	\eeqnn
 	%
 	%
 \end{itemize}
 \qed

 \begin{proposition}\label{Prop.2Est.02}
 	If $m=1$, $\sigma =  {\it\Psi}_1(\infty)<\infty$ and $ R(t)-1/\sigma$ has constant sign near infinity, then four regimes arise for the second-order approximations of $\mathcal{I}_R$ and $\mathcal{I}_R^2$. 
 	\begin{enumerate}
 		\item[(1)] If ${\it \Phi} \in {\rm RV}^\infty_{-1}$, we have as $t\to\infty$, 
 		\beqnn
 		\mathcal{I}_R(t) - \frac{t}{\sigma} 
 		\sim \frac{1}{\sigma^2} \cdot t \big(\sigma-  {\it\Psi}_1(t) \big)  \in {\rm RV}^\infty_{1}
 		\quad\mbox{and}\quad
 		\mathcal{I}_R^2(t)- \frac{t^2}{2\sigma} 
 		\sim  \frac{1}{2\sigma^2} \cdot t^2 \big(\sigma-  {\it\Psi}_1(t) \big) \in {\rm RV}^\infty_{2}. 
 		\eeqnn
 		
 		\item[(2)] If ${\it \Phi} \in {\rm RV}^\infty_{-\alpha}$ with $\alpha\in(1,2)$, we have as $t\to\infty$, 
 		\beqnn
 		\mathcal{I}_R(t) -\frac{t}{\sigma}
 		\sim - \frac{\Gamma(1-\alpha)}{\Gamma(3-\alpha) \sigma^2} \cdot t^2 {\it\Phi}(t) \in {\rm RV}^\infty_{2-\alpha}
 		\quad\mbox{and}\quad
 		\mathcal{I}_R^2(t) -\frac{t^2}{2\sigma}
 		\sim   - \frac{\Gamma(1-\alpha)}{\Gamma(4-\alpha) \sigma^2} \cdot t^3 {\it\Phi}(t) \in {\rm RV}^\infty_{3-\alpha}. 
 		\eeqnn
 		
 		\item[(3)] If ${\it \Phi} \in {\rm RV}^\infty_{-2}$ and ${\it \Psi}_2(\infty)=\infty$, 
 		we have as $t\to\infty$, 
 		\beqnn
 		\mathcal{I}_R(t) - \frac{t}{\sigma}
 		\sim  \frac{ {\it\Psi}_2(t)}{2\sigma^2}  \in {\rm RV}^\infty_{0}
 		\quad\mbox{and}\quad
 		\mathcal{I}_R^2(t) -\frac{t^2}{2\sigma}
 		\sim   \frac{t{\it\Psi}_2(t)}{2\sigma^2} \in {\rm RV}^\infty_{1}. 
 		\eeqnn
 		
 		\item[(4)] If ${\it \Psi}_2(\infty)<\infty$, we have as $t\to\infty$, 
 		\beqnn
 		\mathcal{I}_R(t)- \frac{t}{\sigma} 
 		\sim  \frac{{\it \Psi}_2(\infty)}{2\sigma^2} -1 \in {\rm RV}^\infty_{0}
 		\quad\mbox{and}\quad
 		\mathcal{I}_R^2(t) - \frac{t^2}{2\sigma}
 		\sim \Big(  \frac{{\it \Psi}_2(\infty)}{2\sigma^2} -1 \Big) \cdot t \in {\rm RV}^\infty_{1}. 
 		\eeqnn
 	\end{enumerate}
 \end{proposition}
 \proof 
 The desired asymptotic results for $\mathcal{I}_R(t) -t/\sigma$ suggest that we can just as well consider the asymptotics of the function $\mathcal{I}_R(t) -t/\sigma +1$.
 
 Since $ R(t)-1/\sigma$ has constant sign near infinity, the function $\mathcal{I}_R(t) -t/\sigma +1$ is eventually monotone. In view of Proposition~\ref{Thm.KaramataTauberian}(2), it hence suffices to study the Laplace-Stieltes transform, i.e~ to prove that 
 as $\lambda\to \infty$,
 \beqlb\label{eqn.50004}
 \int_0^\infty \frac{1}{\lambda} e^{-t/\lambda}  \big( \mathcal{I}_R(t) - \frac{t}{\sigma} +1 \big) dt
 \sim
 \begin{cases}
 	\displaystyle{	 \frac{1}{\sigma^2} \cdot \lambda (\sigma-  {\it\Psi}_1(\lambda)) \in {\rm RV}^\infty_{1}}, & \mbox{if ${\it \Phi} \in {\rm RV}^\infty_{-1}$};  \vspace{5pt}\\
 	\displaystyle{\frac{-\Gamma(1-\alpha)}{\sigma^2} \cdot \lambda^2 {\it\Phi}(\lambda) \in {\rm RV}^\infty_{2-\alpha}}, & \mbox{if ${\it \Phi} \in {\rm RV}^\infty_{-\alpha}$ with $\alpha\in(1,2)$};  \vspace{5pt}\\
 	\displaystyle{\frac{1}{2\sigma^2} \cdot {\it \Psi}_2(\lambda)  \in {\rm RV}^\infty_{0}}, & \mbox{if ${\it \Phi} \in {\rm RV}^\infty_{-2}$  and ${\it \Psi}_2(\infty)= \infty$};  \vspace{5pt}\\
 	\displaystyle{ \frac{ {\it \Psi}_2(\infty)}{2\sigma^2}  \in {\rm RV}^\infty_{0}}, & \mbox{if $ {\it \Psi}_2(\infty)<\infty$.}
 \end{cases} 
 \eeqlb
 
 Using the explicit representations of the Laplace-Stieltjes transforms (\ref{LapPhi}) and (\ref{eqn.LapIR}), an application of the mean-value theorem shows that
 \beqlb\label{eqn.50003}
 \int_0^\infty \frac{1}{\lambda} e^{-t/\lambda}  (\mathcal{I}_R(t) - \frac{t}{\sigma} +1) dt
 \ar=\ar \frac{\int_0^\infty e^{-t/\lambda} \phi(t)dt}{\int_0^\infty (1-e^{-t/\lambda} )\phi(t)dt} -\frac{\lambda}{\sigma} +1 \cr
 \ar=\ar \frac{\lambda}{\sigma}\frac{\int_0^\infty (e^{-t/\lambda}-1+t/\lambda)\phi(t)dt}{\int_0^\infty (1-e^{-t/\lambda} )\phi(t)dt}  \cr
 \ar\sim\ar  \frac{\lambda^2}{\sigma^2}  \int_0^\infty \Big(e^{- t/ \lambda }-1+\frac{t}{\lambda}\Big)\phi(t)dt ,
 \eeqlb
 as $\lambda \to\infty$. 
 Since $\phi$ is a probability density function on $\mathbb{R}_+$ with tail-distribution ${\it\Phi}$ and finite first moment ${\it\Psi}_1(\infty)$, it follows from the implication $(2) \Rightarrow (1)$ in Proposition~\ref{Prop.A.7.1} that as $\lambda \to\infty$,
 \beqlb\label{eqn.5002}
 \int_0^\infty \Big(e^{-t/\lambda }-1+\frac{t}{\lambda}\Big)\phi(t)dt \sim 
 \begin{cases}
 	\displaystyle{
 		\frac{\sigma -{\it\Psi}_1(\lambda) }{\lambda}} \in {\rm RV}^\infty_{-1}, & \mbox{if ${\it \Phi} \in {\rm RV}^\infty_{-1}$}; \vspace{5pt}\\
 	\displaystyle{-\Gamma(1-\alpha)\cdot {\it\Phi}(\lambda) \in {\rm RV}^\infty_{-\alpha},} & \mbox{if  ${\it \Phi} \in {\rm RV}^\infty_{-\alpha}$ with $\alpha\in(1,2)$}; \vspace{5pt}\\
 	\displaystyle{   \frac{ {\it \Psi}_2(\lambda)}{2\lambda^2}} \in {\rm RV}^\infty_{-2}, & \mbox{if ${\it \Phi} \in {\rm RV}^\infty_{-2}$}. 
 	\end{cases}
 	\eeqlb
 	Taking this back into (\ref{eqn.50003}) proves the second-order approximation of the function  $\mathcal{I}_R(t)$. The corresponding result for the function $\mathcal{I}^2_R(t)$
 	can be obtained from Proposition~\ref{Thm.Karamata}(1) and our asymptotic results for $\mathcal{I}_R(t) -t/\sigma$.  
 	\qed
 	
 	\smallskip
 	
 	The next proposition analyzes the second-order regular variation of $\mathcal{I}_R$ and $\mathcal{I}_R^2$ if ${\it\Phi} \in {\rm 2RV}^\infty_{-\alpha,\rho}(A)$ for some $\alpha\in (0,1)$, $\rho\leq 0$ and $A\in \mathscr{A}^\infty_\rho$. 
 	By Proposition~\ref{Thm.AsymR}(3), we have  as $t\to\infty$,
 	\beqlb\label{eqn.70002}
 	t^\alpha {\it\Phi}(t)\to C_{\it\Phi} \neq 0, \quad   
 	\frac{\mathcal{I}_R(t)}{t^{\alpha} }\to C_{\mathcal{I}_R}:= \frac{1}{\Gamma(1-\alpha)\Gamma(1+\alpha)C_{\it\Phi}},
 	\quad\ 
 	\frac{\mathcal{I}_R^2(t)}{t^{\alpha+1}}\to C_{\mathcal{I}_R^2}:= \frac{1}{\Gamma(1-\alpha)\Gamma(2+\alpha)C_{\it\Phi}}.
 	\eeqlb
 	
 	\smallskip
 	
 	\begin{proposition}\label{Prop.2Est.03}
 Let $m=1$ and $\sigma = \infty$.  For $\alpha\in (0,1)$, $\rho \in(\alpha-1,0]$ and $A\in \mathscr{A}^\infty_\rho$, 
 if ${\it\Phi} \in {\rm 2RV}^\infty_{-\alpha,\rho}(A)$ and  $\mathcal{I}_R \in \mathscr{M}_{\alpha, -}^\infty$ when $\rho<0$ or  $\mathcal{I}_R \in \mathscr{M}_{\alpha, 0}^\infty$ when $\rho=0$, then the following hold.
 \begin{enumerate}
 	\item[(1)] If $\rho <-\alpha$, then $\displaystyle{\mathcal{I}_R \in {\rm 2RV}^\infty_{\alpha,-\alpha} \big(\alpha \Gamma(1-\alpha)\Gamma(1+\alpha) \cdot {\it\Phi}  \big)}$ and $\displaystyle{\mathcal{I}_R^2 \in {\rm 2RV}^\infty_{\alpha+1,-\alpha} \big(\alpha \Gamma(1-\alpha)\Gamma(2+\alpha) \cdot {\it\Phi}  \big)}$. 
 	\medskip
 	\item[(2)]  If $\rho >-\alpha$,  then $\displaystyle{\mathcal{I}_R \in {\rm 2RV}^\infty_{\alpha,\rho} \left(- \frac{\Gamma(1+\alpha)\Gamma(1-\alpha+\rho)}{\Gamma(1+\alpha+\rho)\Gamma(1-\alpha)} \cdot A \right) }$ and $\displaystyle{ \mathcal{I}_R^2 \in {\rm 2RV}^\infty_{\alpha+1,\rho} \left(-  \frac{\Gamma(2+\alpha)\Gamma(1-\alpha+\rho)}{\Gamma(2+\alpha+\rho)\Gamma(1-\alpha)}\cdot A  \right)}$. 
 	
 \end{enumerate} 
 \end{proposition}
 \proof 
 We provide a detailed proof of the second-order regular variation of $\mathcal{I}_R^2 $ by using Theorem~\ref{Thm.2KaraTaub01} and~\ref{Thm.2KaraTaub02}.
 The proof for $\mathcal{I}_R$ is similar.  
 
 Let us first consider the second-order regular variation of the Laplace-Stieltjes transform $\hat{\mathcal{I}}^2_R(1/\cdot)$. 
 In view of (\ref{eqn.LapIR}),  
 \beqlb\label{eqn.5001}
 \frac{\hat{\mathcal{I}}^2_R(1/(\lambda x))}{\hat{\mathcal{I}}^2_R(1/\lambda)} -x^{\alpha+1} 
 \ar=\ar   x \cdot \frac{{\it\hat\Phi}(1/\lambda)}{{\it\hat\Phi}(1/(\lambda x))} \cdot \frac{1-{\it\hat\Phi}(1/\lambda)/{\it\hat\Phi}(1/(\lambda x))}{1-{\it\hat\Phi}(1/\lambda)} \cdot {\it\hat\Phi}(1/\lambda)\cr
 \ar\ar  -x^{\alpha+1} \cdot \Big( \frac{{\it\hat\Phi}(1/(\lambda x))}{{\it\hat\Phi}(1/\lambda)} - x^{-\alpha}  \Big) \cdot \frac{{\it\hat\Phi}(1/\lambda)}{{\it\hat\Phi}(1/(\lambda x))}, 
 \eeqlb
 for any $\lambda,x>0$. 
 By Proposition~\ref{Thm.KaramataTauberian} and ${\it \Phi} \in {\rm RV}^\infty_{-\alpha}$ with $\alpha\in(0,1)$, we have ${\it\hat\Phi}(1/\cdot) \in {\rm RV}^\infty_{-\alpha}$ and hence as $\lambda \to\infty$,
 \beqlb \label{eqn.4003}
 {\it\hat\Phi}(1/\lambda) \sim \Gamma(1-\alpha) {\it \Phi} (\lambda) , \quad \frac{{\it\hat\Phi}(1/(\lambda x))} {{\it\hat\Phi}(1/\lambda)} \sim x^{-\alpha},\quad 
 1- \frac{{\it\hat\Phi}(1/(\lambda x))} {{\it\hat\Phi}(1/\lambda)}\sim \alpha \int_1^x u^{-\alpha-1}du . 
 \eeqlb
 An application of Theorem~\ref{Thm.2KaraTaub01}(1) and Theorem \ref{Thm.2KaraTaub02}(1) together with the assumption that ${\it \Phi} \in {\rm 2RV}^\infty_{-\alpha,\rho}(A)$ yields that ${\it \hat \Phi}(1/\cdot) \in {\rm 2RV}_{-\alpha,\rho}^\infty \Big(\frac{\Gamma(1-\alpha+\rho)}{\Gamma(1-\alpha)}  \cdot A \Big)$, and so
 \beqlb \label{eqn.4004}
 \frac{{\it\hat\Phi}(1/(\lambda x))}{{\it\hat\Phi}(1/\lambda)} - x^{-\alpha} \sim 
 \begin{cases}
 \displaystyle{ x^{-\alpha}\cdot \frac{\Gamma(1-\alpha+\rho)}{\Gamma(1-\alpha)} \cdot \frac{A(\lambda)}{\rho} \int_1^x u^{\rho-1}du,} & \mbox{if $\rho<0$}; \vspace{5pt} \cr
 \displaystyle{ x^{-\alpha}\log(x)\cdot A(\lambda),} & \mbox{if $\rho=0$} , 
 \end{cases}
 \eeqlb
 as $\lambda \to\infty$. 
 Plugging all estimates in (\ref{eqn.4003}) and (\ref{eqn.4004}) back into the right-hand side of (\ref{eqn.5001}), we can obtain the following asymptotic results. 
 
 {\bf Case I.} 
 If $\rho<0$, we have as $\lambda \to\infty$,
 \beqlb\label{eqn.400004}
 \frac{ \hat{\mathcal{I}}^2_R (1/(\lambda x))}{\hat{\mathcal{I}}^2_R (1/\lambda) } -x^{\alpha+1}
 \ar\sim\ar  x^{\alpha+1}\int_1^x u^{-\alpha-1}du \cdot \alpha\Gamma(1-\alpha) {\it \Phi}(\lambda)  + x^{\alpha+1} \int_1^x u^{\rho-1}du\cdot \frac{-\Gamma(1-\alpha+\rho)}{\Gamma(1-\alpha)}   A(\lambda) . 
 \eeqlb
 If $\rho<-\alpha$, then $A(\lambda)= o({\it \Phi}(\lambda) )$ as $\lambda\to\infty$ and hence
 \beqnn
 \hat{\mathcal{I}}^2_R (1/\cdot) \in {\rm 2RV}^\infty_{\alpha+1, \rho}   \big(\alpha\Gamma(1-\alpha) \cdot {\it \Phi} \big).
 \eeqnn
 The assumption that $\mathcal{I}_R \in \mathscr{M}_{\alpha,-}^\infty$ yields that $\mathcal{I}_R^2 \in \mathscr{M}_{\alpha+1,-}^\infty$, and so it follows from 
 Theorem~\ref{Thm.2KaraTaub01}(1) that 
 \beqnn
 \mathcal{I}_R^2   \in {\rm 2RV}^\infty_{\alpha+1,-\alpha}\big(\alpha\Gamma(\alpha+2)
 \Gamma(1-\alpha) \cdot {\it\Phi} \big).
 \eeqnn 
 Similarly, when $\rho>-\alpha$, we have ${\it \Phi}(\lambda)= o(A(\lambda) )$ as $\lambda \to\infty$ and hence 
 \beqnn
 \hat{\mathcal{I}}^2_R (1/\cdot) \in {\rm 2RV}^\infty_{\alpha+1, \rho}   \Big(-\frac{\Gamma(1-\alpha+\rho)}{\Gamma(1-\alpha)}   A \Big).
 \eeqnn
 The desired result follows by using Theorem~\ref{Thm.2KaraTaub01}(2). 
 
 {\bf Case II.} If $\rho=0$, we also have as $\lambda \to \infty$,
 \beqnn
 \frac{ \hat{\mathcal{I}}^2_R (1/(\lambda x))}{\hat{\mathcal{I}}^2_R (1/\lambda) } -x^{\alpha+1}
 \ar\sim\ar - x^{\alpha+1}\log(x) \cdot   A(\lambda) ,
 \eeqnn
 which induces that $\hat{\mathcal{I}}^2_R (1/\cdot) \in {\rm 2RV}^\infty_{\alpha+1, 0}\big(-A\big)$. 
 By Theorem~\ref{Thm.2KaraTaub02}(2), we have $ \mathcal{I}^2_R  \in {\rm 2RV}^\infty_{\alpha+1, 0} (-A )$. 
 \qed
 
 %
 %
 \begin{remark}
 \begin{itemize}
 	\item[i)] For $\alpha\in(0,1)$, $\rho\in(-2-\alpha,0)$ and $A\in \mathscr{A}^\infty_\rho$, assume $\mathcal{I}^k_{\it\Phi}\in {\rm 2RV}^\infty_{k-\alpha,\rho}\big(\frac{\Gamma(k-\alpha+\rho+1)}{\Gamma(k-\alpha+1)}\cdot A)$ for some non-negative integer $k>\alpha-\rho-1$,  applying integration by parts $k$-times to $\hat{\mathcal{I}}_{\it\Phi}^k(1/\lambda)$ gives that 
 	\beqnn
 	{\it\hat\Phi}(1/\lambda)= \lambda^{-k} \cdot \hat{\mathcal{I}}_{\it\Phi}^k(1/\lambda).
 	\eeqnn
 	Applications of Theorem~\ref{Thm.2KaraTaub01} and Proposition~\ref{Prop.Fr} show that ${\it \hat \Phi}(1/\cdot) \in {\rm 2RV}_{-\alpha,\rho}^\infty \Big(\frac{\Gamma(1-\alpha+\rho)}{\Gamma(1-\alpha)}  \cdot A \Big)$. Repeating the proof of Proposition~\ref{Prop.2Est.03},  one can see that if $\mathcal{I}_R^2 \in \mathscr{M}_{\alpha+1,-}^\infty$, then results in Proposition~\ref{Prop.2Est.03} still hold for $\mathcal{I}_R^2$. 
 	Additionally, by   Proposition~\ref{Thm.AsymR}(3), we see that limits in (\ref{eqn.70002}) also hold. 
 	
 	\item[ii)] If $\rho \in (-K-\alpha-1,-K-\alpha]$ with $K\geq 3$, under the assumptions that $\mathcal{I}^k_{\it \Phi}  \in {\rm 2RV}^\infty_{k-\alpha,\rho}(A)$ and
 	$\mathcal{I}^{k_i}_{{\it \Psi}_1}  \in {\rm 2RV}^\infty_{k_i+1-\alpha,\rho_i}( A_1)$ with $\rho_i\leq 0$, $k_i>\alpha-\rho_i-i-1$, $A_i\in\mathscr{A}^\infty_{\rho_i}$ for  $i=1,\cdots, K$, similarly as in the proof of Proposition~\ref{Prop.2Est.03}, one can prove that $\mathcal{I}^2_R\in {\rm 2RV}^\infty_{\alpha+1,\varrho_K}(A^*_K)$ with $\varrho_K:=
 	(-\alpha)\vee\rho_1\vee\cdots\vee \rho_K$ and $A^*_K \in\mathscr{A}^\infty_{\varrho_K}$.
 \end{itemize}
 \end{remark}
 
 Using the precedingly establsihed second-order asymptotics for $\mathcal{I}_R$ and $\mathcal{I}_R^2$, we are now ready to establish the second-order approximations for $\mathbf{E} [N(t) ]$ and ${\rm Var} ( N(t) )$ in various settings.

 \begin{theorem}\label{Thm.2RV.HP}
 The following regimes arise for the second-order approximations of $\mathbf{E}[N(t)]$ and ${\rm Var}(N(t))$.
 \begin{enumerate}
 	\item[(1)] (Subcritical case) When $m<1$, we have as $t\to\infty$,
 	
 	\begin{enumerate}
 		\item[(1.a)]  if ${\it\Phi} \in {\rm RV}^\infty_{-\alpha}$ with $\alpha\in[0,1)$, then 
 		\beqnn
 		\mathbf{E}\big[N(t)\big] - \frac{\mu_0\cdot t}{1-m} \sim  - \frac{\mu_0 \cdot t{\it\Phi}(t)}{(1-m)^2(1-\alpha)} 
 		\quad \mbox{and}\quad
 		{\rm Var}\big( N(t)\big) - \frac{\mu_0\cdot t}{(1-m)^3} \sim - 
 		\frac{3\mu_0\cdot t{\it\Phi}(t)}{(1-m)^4(1-\alpha)};  
 		\eeqnn
 		
 		\item[(1.b)] if ${\it\Phi} \in {\rm RV}^\infty_{-1}$ and $\sigma={\it\Psi}_1(\infty)=\infty$, then 
 		\beqnn
 		\mathbf{E}\big[N(t)\big] - \frac{\mu_0\cdot t}{1-m} \sim - \frac{\mu_0 \cdot {\it\Psi}_1(t)}{(1-m)^2} 
 		\quad \mbox{and}\quad
 		{\rm Var}\big( N(t)\big) - \frac{\mu_0\cdot t}{(1-m)^3} \sim - 
 		\frac{3\mu_0\cdot {\it\Psi}_1(t)}{(1-m)^4};  
 		\eeqnn
 		
 		\item[(1.c)] if ${\it\Phi} \in {\rm RV}^\infty_{-\alpha}$ with $\alpha>1$ or if  ${\it\Phi} \in {\rm RV}^\infty_{-1}$ and   $\sigma={\it\Psi}_1(\infty)<\infty$, then 
 		\beqnn
 		\mathbf{E}\big[N(t)\big] - \frac{\mu_0\cdot t}{1-m} \sim - \frac{\mu_0 \sigma }{(1-m)^2} 
 		\quad \mbox{and}\quad
 		{\rm Var}\big( N(t)\big) - \frac{\mu_0\cdot t}{(1-m)^3} \sim - 
 		\frac{3\mu_0  \sigma }{(1-m)^4}. 
 		\eeqnn
 	\end{enumerate}

 	\item[(2)] (Weakly critical case) When $m=1$ and $\sigma={\it\Psi}_1(\infty)<\infty$, we have as $t\to\infty$,
 	
 	\begin{enumerate}
 		\item[(2.a)] if ${\it\Phi} \in {\rm RV}^\infty_{-1}$, then 
 		\beqnn
 		\mathbf{E}\big[N(t)\big] -\frac{\mu_0\cdot t^2 }{2\sigma} \sim \frac{\mu_0}{2\sigma^2} \cdot t^2 \big(\sigma-  {\it\Psi}_1(t)\big)
 		\quad \mbox{and}\quad 
 		{\rm Var}\big( N(t)\big)- \frac{\mu_0\cdot t^4}{12\sigma^3} \sim \frac{\mu_0}{4\sigma^4} \cdot t^4 \big(\sigma-  {\it\Psi}_1(t)\big);
 		\eeqnn
 		
 		\item[(2.b)] if ${\it\Phi} \in {\rm RV}^\infty_{-\alpha}$ with $\alpha\in(1,2)$, then
 		\beqnn
 		\mathbf{E}\big[N(t)\big] - \frac{\mu_0\cdot t^2}{2\sigma}  \sim
 		- \frac{\mu_0 \Gamma(1-\alpha) }{\Gamma(4-\alpha) \sigma^2} \cdot  t^3{\it\Phi}(t)
 		\quad \mbox{and}\quad 
 		{\rm Var}\big( N(t)\big) -\frac{\mu_0\cdot t^4}{12\sigma^3}  \sim 
 		- \frac{2\mu_0\Gamma(1-\alpha) }{\Gamma(4-\alpha)\sigma^4} \cdot \frac{ t^5{\it\Phi}(t)}{5-\alpha} ; 
 		\eeqnn 		
 		
 		\item[(2.c)] if ${\it\Phi} \in {\rm RV}^\infty_{-2}$ and ${\it \Psi}_2(\infty)=\infty$, then 
 		\beqnn
 		\mathbf{E}\big[N(t)\big] - \frac{\mu_0\cdot t^2 }{2\sigma} \sim \frac{\mu_0}{2\sigma^2} \cdot t{\it \Psi}_2(t)
 		\quad \mbox{and}\quad 
 		{\rm Var}\big( N(t)\big) - \frac{\mu_0 \cdot t^4 }{12\sigma^3} \sim  \frac{\mu_0}{3\sigma^4} \cdot t^3 {\it \Psi}_2(t); 
 		\eeqnn
 		
 		\item[(2.d)] if ${\it \Psi}_2(\infty)<\infty$, then
 		\beqnn
 		\mathbf{E}\big[N(t)\big]-  \frac{\mu_0\cdot t^2}{2\sigma} \sim  \frac{\mu_0{\it \Psi}_2(\infty)}{2\sigma^2} \cdot t
 		\quad \mbox{and}\quad 
 		{\rm Var}\big( N(t)\big) - \frac{\mu_0\cdot t^4}{12\sigma^3}  \sim   \frac{\mu_0{\it \Psi}_2(\infty)}{3\sigma^4} \cdot t^3. 
 		\eeqnn
 	\end{enumerate}
 	
 	\item[(3)] (Strongly critical case) When $m=1$ and $\sigma={\it \Psi}_1(\infty)=\infty$, 
 	if ${\it\Phi} \in {\rm 2RV}^\infty_{-\alpha,\rho}(A)$ and  $\mathcal{I}_R \in \mathscr{M}_{\alpha, -}^\infty$ with $\alpha\in (0,1)$, $\rho \in(\alpha-1,0)$ and $A\in \mathscr{A}^\infty_\rho$, we have as $t\to\infty$,
 	
 	\begin{enumerate}
 		\item[(3.a)] if $\rho <-\alpha$, then 
 		\beqnn
 		\mathbf{E}\big[N(t)\big]-  \mu_0 C_{\mathcal{I}^2_R}\cdot t^{\alpha+1}  \sim o (1)
 		\quad \mbox{and}\quad 
 		{\rm Var}\big( N(t)\big) - \mu_0 |C_{\mathcal{I}_R}|^3\cdot \mathrm{B}(2\alpha+1,\alpha+1)\cdot t^{3\alpha+1}
 		\sim   o(t) ;
 		\eeqnn
 		
 		\item[(3.b)]  if $\rho >-\alpha$,  then 
 		\beqnn
 		\mathbf{E}\big[N(t)\big]-   \mu_0 C_{\mathcal{I}^2_R}\cdot t^{\alpha+1} \sim - C^{\tt HP}_1\cdot\frac{\mu_0 C_{\mathcal{I}^2_R}}{\rho}\cdot t^{\alpha+1}A(t)
 		\eeqnn
 		and
 		\beqnn
 		{\rm Var}\big( N(t)\big) - \mu_0 |C_{\mathcal{I}_R}|^3\cdot \mathrm{B}(2\alpha+1,\alpha+1)\cdot t^{3\alpha+1}
 		\sim -C^{\tt HP}_2\cdot \frac{\mu_0|C_{\mathcal{I}_R}|^3}{\rho}\cdot t^{3\alpha+1}A(t),
 		\eeqnn
 		with $C^{\tt HP}_1:=   \frac{\mathrm{B}(2+\alpha,1-\alpha+\rho)}{\mathrm{B}(2+\alpha+\rho,1-\alpha)} $ and $C^{\tt HP}_2:=  \big(2\cdot\mathrm{B}(2\alpha+\rho+1,\alpha+1)+ \mathrm{B}(2\alpha+1,\alpha+\rho+1)\big)
 		\frac{\mathrm{B}(\alpha+1,1-\alpha+\rho)}{\mathrm{B}(\alpha+1+\rho,1-\alpha)}$.
 	\end{enumerate} 
 \end{enumerate}
 
 \end{theorem}
 \proof We provide a detailed proof for claim (1.a) and claim (3.a). 
 All other claims can be proved in the same way.  
 
 In view of (\ref{eqn.3007}), we first need to establish the exact second-order approximations of the functions $\mathcal{I}_{R}$, $\mathcal{I}_{R}^2$, $\mathcal{I}_{|\mathcal{I}_R|^2}$, $\mathcal{I}_{R} * \mathcal{I}_{R}$ and $|\mathcal{I}_{R}|^2 * \mathcal{I}_{R}$. If $m<1$ and ${\it\Phi} \in {\rm RV}^\infty_{-\alpha}$ with $\alpha\in[0,1)$, then it follows from 
 Proposition~\ref{Prop.2Est.01} that
 \beqnn 
 \mathcal{I}_R(t) =\frac{m}{1-m} - \frac{ {\it\Phi}(t)}{(1-m)^2} +o({\it\Phi}(t))
 \quad \mbox{and} \quad
 \mathcal{I}_R^2 (t) =  \frac{m \cdot t}{1-m}-\frac{t{\it \Phi}(t)}{(1-m)^2(1-\alpha)} +o(t{\it\Phi}(t)) .  
 \eeqnn
 By Corollary~\ref{Appendix.Coro01} and Proposition~\ref{Thm.Karamata}(1), we have as $t\to\infty$, 
 \beqnn
 \mathcal{I}_{|\mathcal{I}_R|^2}(t)
 \ar =\ar \frac{m^2 \cdot t}{(1-m)^2}  - \frac{ 2m}{(1-m)^3} \int_0^t {\it\Phi}(s)ds+ \int_0^t o({\it\Phi}(s))ds \cr
 \ar\sim \ar \frac{m^2\cdot t}{(1-m)^2} -  \frac{2m}{(1-m)^3} \cdot \frac{t {\it\Phi}(t)}{1-\alpha},\cr
 \ar\ar\cr 
 \mathcal{I}_R*\mathcal{I}_R(t) 
 \ar=\ar \frac{m^2\cdot t}{(1-m)^2} -\frac{2m}{(1-m)^3}\int_0^t   {\it\Phi}(s) ds + \frac{2m}{1-m} \int_0^t o({\it\Phi}(s))ds\cr
 \ar\ar -\frac{2}{(1-m)^2} \int_0^t  {\it\Phi}(t-s) o({\it\Phi}(s))  ds + \int_0^t  o({\it\Phi}(t-s) {\it\Phi}(s))  ds  \cr
 \ar\sim\ar \frac{m^2\cdot t}{(1-m)^2} -  \frac{2m}{(1-m)^3} \cdot \frac{t {\it\Phi}(t)}{1-\alpha}
 \eeqnn
 and 
 \beqnn
 |\mathcal{I}_R|^2*\mathcal{I}_R(t) 
 \ar=\ar \frac{m^3\cdot t}{(1-m)^3} - \int_0^t\frac{3 m^2{\it\Phi}(s)}{(1-m)^4} ds + \frac{ 1+m}{(1-m)^4}\int_0^t {\it\Phi}(t-s){\it\Phi}(s)ds \cr
 \ar\ar + \int_0^t \frac{m|{\it\Phi}(s)|^2}{(1-m)^5}  ds + \int_0^t o({\it\Phi}(t-s){\it\Phi}(s)) ds  + \int_0^t  o({\it\Phi}(s)) ds 
 \cr
 \ar\sim\ar \frac{m^3\cdot t}{(1-m)^3} -  \frac{3m^2}{(1-m)^4} \cdot \frac{t {\it\Phi}(t)}{1-\alpha} . 
 \eeqnn
 Plugging the above into (\ref{eqn.3007}) shows that claim (1.a) holds. To prove (3.a), we first deduce from Proposition~\ref{Prop.2Est.03} and Corollary~\ref{Coro.KaraRep} that
 \beqnn
 \mathcal{I}_R(t) = C_{\mathcal{I}_R} \cdot t^\alpha -1 +o(1)
 \quad\mbox{and}\quad 
 \mathcal{I}_R^2(t) \sim C_{\mathcal{I}_R^2} \cdot t^{\alpha+1} -t + o(t),
 \eeqnn
 as $t\to\infty$. 
 Similarly as in the previous argument, by Corollary~\ref{Appendix.Coro01} and Proposition~\ref{Thm.Karamata}(1) we have as $t\to\infty$, 
 \beqnn
 \mathcal{I}_{ |\mathcal{I}_R|^2 }(t)  
 \ar=\ar  \frac{|C_{\mathcal{I}_R}|^2}{2\alpha+1} \cdot t^{2\alpha+1} - \frac{2 C_{\mathcal{I}_R}}{\alpha+1} \cdot t^{\alpha+1} + t +o(t) ,\cr
 \mathcal{I}_R * \mathcal{I}_R (t)\ar=\ar |C_{\mathcal{I}_R}|^2 \mathrm{B}(1+\alpha,1+\alpha) \cdot t^{2\alpha+1}  - \frac{2 C_{\mathcal{I}_R}}{1+\alpha} \cdot t^{\alpha+1} + t +o(t),\cr
 |\mathcal{I}_R|^2 * \mathcal{I}_R(t) \ar=\ar |C_{\mathcal{I}_R}|^3 \mathrm{B}(1+2\alpha,1+\alpha) \cdot t^{3\alpha+1} + \frac{3 C_{\mathcal{I}_R}}{1+\alpha} \cdot t^{\alpha+1} - t\cr
 \ar\ar -  \big(2\cdot\mathrm{B}(1+\alpha,1+\alpha)+ 1/(1+2\alpha )\big) \cdot |C_{\mathcal{I}_R}|^2   \cdot t^{2\alpha+1} +o(t).
 \eeqnn
 Plugging these asymptotic results into (\ref{eqn.3007}), one can get claim (3.a) immediately with a simple calculation. 
 \qed

 \subsection{Examples}
 In contrast to the strongly critical case, the preceding first- and second-order approximations for subcritical and weakly critical Hawkes processes were established under mild conditions on the function ${\it\Phi}$ that are easily satisfied in many examples. In this section, we provide two concrete examples of second-order approximations of strongly critical Hawkes processes where all our assumptions can be verified. 
 
 
 For $\alpha\in(0,1)$ and $\kappa>0$, we denote by $E_{\alpha,\kappa} $ the \textit{Mittag-Leffler function} 
 \beqlb\label{MLF}
 E_{\alpha,\kappa}(t) := \sum_{k=0}^\infty \frac{t^k}{\Gamma(\kappa+k\alpha)},\quad t\in\mathbb{R}.
 \eeqlb
 It is customary to write $E_{\alpha}$ for $E_{\alpha,1}$. We refer to \cite{MathaiHaubold2008,HauboldMathaiSaxena2011}  for an in-depth analysis of Mittag-Leffler functions. For some constant $\beta>0$ we denote by
 \beqlb\label{MLdistribution}
 F^{\alpha,\beta}(t):= 1- E_\alpha(-\beta t^\alpha)
 \quad \mbox{and}\quad
 f^{\alpha,\beta}(t): = \beta t^{\alpha-1}  E_{\alpha,\alpha}(-\beta t^{\alpha}).
 \eeqlb
 the \textit{Mittag-Leffler distribution} and the \textit{Mittag-Leffler density function}  on $\mathbb{R}_+$, respectively.  The Laplace transform of Mittag-Leffler distribution admits the representation
 \beqlb\label{ML.LP}
 \hat{F}^{\alpha,\beta}(\lambda) = \int_0^\infty \lambda e^{-\lambda t} F^{\alpha,\beta}(t)dt = \int_0^\infty  e^{-\lambda t} f^{\alpha,\beta}(t)dt=  \frac{\beta}{\beta+\lambda^\alpha},\quad \lambda > 0.
 \eeqlb
 Moreover, for $\alpha\in(0,1)$ and $\beta>0$, the asymptotic expansion of Mittag-Leffler functions given in \cite[Chapter XVIII]{ErdelyiMagnusOberhettingerTricomi1955} and \cite[Section~6]{HauboldMathaiSaxena2011} shows that as $t\to\infty$,
 \beqlb\label{eqn.ML.2RV}
 f^{\alpha,\beta} (t) \ar\sim\ar  \frac{-  t^{-\alpha-1}}{\beta\Gamma(-\alpha)} + \frac{ t^{-2\alpha-1}}{\beta^2\Gamma(-2\alpha)} + O(t^{-3\alpha-1})
 \quad\mbox{and}\quad
 1- F^{\alpha,\beta} (t)  \sim \frac{ t^{-\alpha}}{\beta\Gamma(1-\alpha)} -\frac{t^{-2\alpha}}{\beta^2\Gamma(1-2\alpha)} + O(t^{-3\alpha}). \qquad\ 
 \eeqlb  
 
 \subsubsection{Fractional Hawkes processes}
 
 A strongly critical Hawkes process $N$ is said to be of \textit{Mittag-Leffler type} with index $(\alpha,\beta)$, also known as \textit{fractional type}, if 
 \beqnn
 \phi=f^{\alpha,\beta} \quad \mbox{and} \quad {\it\Phi} = 1-F^{\alpha,\beta}.
 \eeqnn
 It follows  from (\ref{ML.LP}) and (\ref{eqn.LapIR}) that
 \beqnn
 {\it\hat\Phi}(\lambda)
 = 1-\hat{F}^{\alpha,\beta}(\lambda)
 =\frac{\lambda^\alpha}{\beta+\lambda^\alpha},\quad
 \hat{\mathcal{I}}_R(\lambda)= \frac{\beta}{\lambda^\alpha} 
 \quad\mbox{and}\quad
 \hat{\mathcal{I}}_R^2(\lambda)=\frac{\beta}{\lambda^{\alpha+1}},\quad 
 \lambda > 0.
 \eeqnn 
 The one-to-one correspondence between functions and their Laplace transforms yields that
 \beqnn
 \mathcal{I}_R(t) = \frac{\beta \cdot  t^{\alpha} }{\Gamma(\alpha+1)}
 \quad \mbox{and}\quad
 \mathcal{I}_R^2(t) = \frac{\beta \cdot  t^{\alpha+1} }{\Gamma(\alpha+2)},\quad t\geq 0.
 \eeqnn
 It is straightforward to show that for $t\geq 0$,
 \beqnn
 \mathcal{I}_{|\mathcal{I}_R|^2}(t) = \frac{\beta^2 \cdot  t^{1+2\alpha} }{(1+2\alpha)|\Gamma(1+\alpha)|^2}, \quad 
 \mathcal{I}_R* \mathcal{I}_R(t) = \frac{\beta^2\cdot t^{1+2\alpha} }{\Gamma(2+2\alpha)},
 \quad 
 |\mathcal{I}_R|^2 * \mathcal{I}_R(t) = \frac{\beta^3\mathrm{B}(1+2\alpha,1+\alpha)}{|\Gamma(1+\alpha)|^3}\cdot t^{3\alpha+1}. 
 \eeqnn
 Plugging these equalities into (\ref{eqn.3007}), we obtain the following exact expressions for the mean and the variance of a fractional Hawkes process: 
 \beqnn
 \mathbf{E}\big[N(t) \big] = \frac{\mu_0\beta \cdot  t^{\alpha+1} }{\Gamma(\alpha+2)}  + \mu_0\cdot t
 \eeqnn
 and
 \beqnn
 {\rm Var}\big( N(t) \big) = \mu_0\Big(\frac{\mu_0\beta^3\mathrm{B}(1+2\alpha,1+\alpha)}{|\Gamma(1+\alpha)|^3}\cdot t^{3\alpha+1}
 + \frac{2\beta^2\cdot t^{1+2\alpha} }{\Gamma(2+2\alpha)}
 + \frac{\beta^2 \cdot  t^{1+2\alpha} }{(1+2\alpha)|\Gamma(1+\alpha)|^2}
 + \frac{3\beta \cdot  t^{\alpha+1} }{\Gamma(\alpha+2)} +t \Big).
 \eeqnn
 
 \subsubsection{Processes of mixed Mittag-Leffler type}
 
 For any four constants $0< \alpha_1\leq \alpha_2< 1$ and $\beta_1,\beta_2>0$,
 let $f^{\alpha_1,\beta_1}$ and $f^{\alpha_2,\beta_2}$ be two Mittag-Leffler density functions with parameters $( \alpha_1, \beta_1)$ and $( \alpha_2,\beta_2)$ respectively, and let 
 $F^{\alpha_1,\beta_1}$ and $F^{\alpha_2,\beta_2}$ be the corresponding distribution functions. 
 Their convolutions are defined by 
 \beqlb\label{Def.MixedML}
 f^{\alpha_1,\beta_1}_{\alpha_2,\beta_2}(t):= f^{\alpha_1,\beta_1} * f^{\alpha_2,\beta_2}(t)
 \quad \mbox{and}\quad
 F^{\alpha_1,\beta_1}_{\alpha_2,\beta_2}(t):= \int_0^t f^{\alpha_1,\beta_1}_{\alpha_2,\beta_2}(s)ds,\quad t>0 . 
 \eeqlb
 A  critical Hawkes process $N$ is said to be of \textit{mixed Mittag-Leffler type} with index $(\alpha_i,\beta_i)_{i=1,2} $  if 
 \[
 \phi= f^{\alpha_1,\beta_1}_{\alpha_2,\beta_2} \quad \mbox{and} \quad
 {\it\Phi}= 1-F^{\alpha_1,\beta_1}_{\alpha_2,\beta_2}. 
 \]
 
 An application of Proposition 3.12 in \cite{Xu2021a} or Lemma 4.9 in \cite[p.77]{FossKorshunovZachary2011}, together with (\ref{eqn.ML.2RV}), yields that as $t\to\infty$,
 \beqlb\label{ML.RV}
 f^{\alpha_1,\beta_1}_{\alpha_2,\beta_2}(t)
 \sim \frac{\mathcal{C}_\beta\alpha_1\cdot t^{-\alpha_1-1}}{ \Gamma(1-\alpha_1)}\in {\rm RV}^\infty_{-\alpha_1 -1} 
 \quad \mbox{and} \quad 
 1-F^{\alpha_1,\beta_1}_{\alpha_2,\beta_2}(t)
 \sim\frac{\mathcal{C}_\beta\cdot t^{-\alpha_1 }}{ \Gamma(1-\alpha_1)} \in {\rm RV}^\infty_{-\alpha_1 }.\label{ML.RV02}
 \eeqlb
 with $$\mathcal{C}_\beta:= \frac{1}{\beta_1}+\frac{1}{\beta_2}\cdot \mathbf{1}_{\{\alpha_1=\alpha_2\}}.$$
 Taking (\ref{ML.RV02}) back into Proposition~\ref{Thm.AsymR}(3), we have as $t\to\infty$, 
 \beqlb\label{ML.RV0202}
 \mathcal{I}_R(t)  \sim C_{\mathcal{I}_R} \cdot t^{\alpha_1} \in {\rm RV}_{\alpha_1}^\infty
 \quad \mbox{and}\quad
 \mathcal{I}^2_R(t) \sim C_{\mathcal{I}_R^2} \cdot t^{\alpha_1+1}  \in {\rm RV}_{\alpha_1+1}^\infty,
 \eeqlb 
 with 
 \[
 C_{\mathcal{I}_R}= \frac{1}{\mathcal{C}_\beta  \Gamma(\alpha_1+1)} \quad \mbox{and}  \quad
 C_{\mathcal{I}_R^2}= \frac{1}{\mathcal{C}_\beta  \Gamma(\alpha_1+2)}.
 \] 
 The next proposition analyzes the second-order regular variation of $\mathcal{I}_R$ and $\mathcal{I}^2_R$.

 \begin{proposition}\label{Prop.Example01}
 The function $\mathcal{I}_R(t)-  \frac{ t^{\alpha_1}}{\mathcal{C}_\beta\Gamma(\alpha_1+1)}$ is non-positive on $\mathbb{R}_+$ and as $t\to\infty$,
 \beqlb\label{eqn.70001}
 \mathcal{I}_R(t)-  \frac{ t^{\alpha_1}}{\mathcal{C}_\beta\Gamma(\alpha_1+1)} 
 \sim 
 \begin{cases}
 	\displaystyle{- \frac{\beta_1\beta_2}{(\beta_1+\beta_2)^2} ,} & \mbox{if $\alpha_1=\alpha_2$}; \vspace{5pt} \cr
 	\displaystyle{-  \frac{\beta_1^2}{\beta_2} \frac{t^{2\alpha_1-\alpha_2}}{\Gamma(1+2\alpha_1-\alpha_2)},} & \mbox{if $\alpha_1\neq \alpha_2$} .
 \end{cases}
 \eeqlb
 \end{proposition}
 \proof 
 To obtain the desired result, we provide an exact representation of the function $\mathcal{I}_R(t)-  \frac{ t^{\alpha_1}}{\mathcal{C}_\beta\Gamma(\alpha_1+1)}$. 
 By the convolution theorem of Laplace transform,
 \beqlb\label{eqn.4.01}
 {\it\hat\Phi}(\lambda)
 \ar=\ar 1- \int_0^\infty e^{-\lambda t} f_{\alpha_1,\beta_1}(t)dt \cdot \int_0^\infty e^{-\lambda t} f_{\alpha_2,\beta_2}(t)dt
 = 1- \frac{1}{1+\lambda^{\alpha_1}/\beta_1} \frac{1}{1+\lambda^{\alpha_2}/\beta_2} .
 \eeqlb
 Plugging this back into (\ref{eqn.LapIR}), we have $\hat{\mathcal{I}}_R(\lambda)= \lambda^{-\alpha_1} \cdot G(\lambda)$ where $G(\lambda)$ is a completely monotone function on $\mathbb{R}_+$ given by
 \beqnn
 G(\lambda)=  \frac{\beta_1\beta_2}{ \beta_2+  \beta_1  \lambda^{\alpha_2-\alpha_1} + \lambda^{\alpha_2} }.
 \eeqnn
 By Bernstein's theorem \cite[Theorem~1.4, p.3]{SchillingSongVondracek2012}, there exists a finite measure $\mathcal{V}_{\alpha}$ on $\mathbb{R}_+$ such that 
 \beqnn
 \int_0^\infty e^{-\lambda t} \mathcal{V}_{\alpha}(dt)= G(\lambda)
 \quad \mbox{and}\quad 
 \mathcal{V}_{\alpha}(\mathbb{R}_+)= \frac{1}{\mathcal{C}_\beta}. 
 \eeqnn 
 
 The function $G(0)-G(1/\cdot) $ belongs to $ {\rm RV}^\infty_{-\alpha_1}$ when $\alpha_1=\alpha_2$ and to $ {\rm RV}^\infty_{\alpha_1-\alpha_2}$  when $\alpha_1<\alpha_2$. 
 By Proposition~\ref{Thm.KaramataTauberian}(2), we have as $t\to\infty$,
 \beqlb\label{eqn.40001}
 \mathcal{V}_{\alpha}(t,\infty) \sim \frac{G(0)-G(1/t)}{\Gamma(1+\rho_\alpha)} \sim 
 \begin{cases}
 \displaystyle{\frac{\beta_1^2}{\beta_2} \frac{t^{\alpha_1-\alpha_2}}{\Gamma(1+\alpha_1-\alpha_2)},} & \mbox{if $\alpha_1\neq \alpha_2$}; \cr
 \displaystyle{\frac{\beta_1\beta_2}{(\beta_1+\beta_2)^2} \frac{t^{-\alpha_1}}{\Gamma(1-\alpha_1)},} & \mbox{if $\alpha_1=\alpha_2$} . 
 \end{cases}
 \eeqlb
 By (\ref{eqn.002}), the function $ \lambda^{-\alpha_1}$ is the Laplace transform of the function  
 $ g (t):= t^{\alpha_1}/\Gamma(\alpha_1+1)$ for $t\geq 0$.
 The one-to-one correspondence between functions and their Laplace transforms induces that 
 \beqnn
 \mathcal{I}_R(t)=  \int_0^t g(t-s)\mathcal{V}_{\alpha}(ds)
 = \int_0^t \int_0^s \frac{ (s-r)^{\alpha_1-1}}{\Gamma(\alpha_1)} \mathcal{V}_{\alpha}(dr)ds
 = \int_0^t   \frac{ (t-s)^{\alpha_1-1}\cdot \mathcal{V}_{\alpha}(0,s)}{\Gamma(\alpha_1)} ds,\quad t\geq 0. 
 \eeqnn 
 As a result, 
 \beqnn
 \mathcal{I}_R(t)- \frac{  t^{\alpha_1}}{\mathcal{C}_\beta\Gamma(\alpha_1+1)} 
 = \mathcal{I}_R(t)- \int_0^t \frac{(t-s)^{\alpha_1-1}}{\Gamma(\alpha_1)}\mathcal{V}_{\alpha}(0,\infty)ds
 = -\int_0^t \frac{(t-s)^{\alpha_1-1}}{\Gamma(\alpha_1)} \cdot \mathcal{V}_{\alpha}(s,\infty)ds,
 \eeqnn
 which is non-positive; hence  the first claim holds. 
 The second claim follows applying Corollary~\ref{Appendix.Coro01} and (\ref{eqn.40001}) to the convolution on the right-hand side of the second equality in the preceding equation. 
 \qed 
 
 
 \begin{proposition}\label{Prop.HP}
 We have $\mathcal{I}_R \in {\rm 2RV}^\infty_{\alpha_1, \rho_\alpha } ( A_{\tt HP})$ and $\mathcal{I}^2_R \in {\rm 2RV}^\infty_{\alpha_1+1, \rho_\alpha } ( \frac{\alpha_1+1}{\alpha_1+\rho_\alpha+1} \cdot A_{\tt HP})$ with 
 \beqnn 
 \rho_\alpha:= \begin{cases}
 	-\alpha_1, & \mbox{if } \alpha_1=\alpha_2;  \vspace{5pt}\\
 	\alpha_1-\alpha_2, & \mbox{if } \alpha_1<\alpha_2 
 \end{cases}
 \quad\mbox{and}\quad
 A_{\tt HP}(t):=
 \begin{cases}
 	\displaystyle{ \frac{\alpha_1\Gamma(1+\alpha_1)}{\beta_1+\beta_2} \cdot t^{-\alpha_1},} & \mbox{if } \alpha_1=\alpha_2;  \vspace{5pt}\\
 	\displaystyle{(\alpha_2-\alpha_1)\frac{\beta_1\Gamma(1+\alpha_1)}{\beta_2\Gamma(1+2\alpha_1-\alpha_2)} \cdot t^{\alpha_1-\alpha_2},} & \mbox{if } \alpha_1<\alpha_2 .
 \end{cases}
 \eeqnn
 \end{proposition}
 \proof To prove the second-order regular variation of $\mathcal{I}^2_R$ we recall the function $G(\lambda)$ defined in the proof of Proposition~\ref{Prop.Example01}. 
 By (\ref{eqn.LapIR}), 
 \beqlb
 \hat{\mathcal{I}}_R^2(1/\lambda)= \lambda \cdot \hat{\mathcal{I}}_R(1/\lambda)= \lambda^{\alpha_1+1}\cdot  G(1/\lambda) = \lambda^{\alpha_1+1}\cdot\frac{\beta_1\beta_2}{\beta_2 + \beta_1\lambda^{\alpha_1-\alpha_2} + \lambda^{-\lambda_2}},\quad \lambda >0.
 \eeqlb
 By the Taylor expansion of the function $1/(1+x)$ at $0$, we have as $\lambda \to \infty$,
 \beqnn
 \hat{\mathcal{I}}_R^2(1/\lambda) \sim
 \begin{cases}
 \displaystyle{\lambda^{\alpha_1+1}\cdot\frac{\beta_1\beta_2}{\beta_1+\beta_2}\Big(1+ \frac{\alpha_1}{\beta_1+\beta_2} \cdot \frac{\lambda^{-\alpha_1}}{-\alpha_1} + o(\lambda^{-\alpha_1})\Big),} & \mbox{if $\alpha_1=\alpha_2$}; \vspace{5pt}\\
 \displaystyle{\lambda^{\alpha_1+1}\cdot\beta_1 \Big(1+ (\alpha_2-\alpha_1) \frac{\beta_1}{\beta_2} \frac{\lambda^{\alpha_1-\alpha_2}}{\alpha_1-\alpha_2} +o(\lambda^{\alpha_1-\alpha_2})  \Big),} & \mbox{if $\alpha_1<\alpha_2$},\\ 
 \end{cases}
 \eeqnn
 From Corollary~\ref{Coro.KaraRep}, we see that $\hat{\mathcal{I}}_R^2(1/\cdot) \in {\rm 2RV}^\infty_{\alpha_1+1,\rho_\alpha}(\frac{\Gamma(1+\alpha_1+\rho_\alpha)}{\Gamma(1+\alpha_1)}\cdot A_{\tt HP})$. 
 Moreover, Proposition~\ref{Prop.Example01} shows that the function $\mathcal{I}_R^2(t)-  C_{\mathcal{I}_R^2} \cdot t^{\alpha_1+1} $ is non-increasing on $\mathbb{R}_+$ and hence $\mathcal{I}_R^2 \in \mathscr{M}^\infty_{\alpha_1+1,-}$.  
 By Theorem~\ref{Thm.2KaraTaub01}, we thus have that
 \[
 \mathcal{I}^2_R \in {\rm 2RV}^\infty_{\alpha_1+1, \rho_\alpha } ( \frac{\alpha_1+1}{\alpha_1+\rho_\alpha+1} A_{\tt HP}).
 \] 
 
 We  consider the second-order regular variation of $\mathcal{I}_R$. 
 By Corollary~\ref{Coro.KaraRep} and the preceding result, we have as $t\to\infty$, 
 \beqlb\label{eqn.40002}
 \mathcal{I}_R^2(t) = C_{\mathcal{I}_R^2}\cdot t^{\alpha_1+1}\cdot \big(1+ t^{\rho_\alpha}\cdot \ell(t)  \big),\quad t>0,
 \eeqlb
 where $\ell$ is a twice differentiable function on $(0,\infty)$ defined by 
 \beqlb\label{eqn.40003}
 \ell(t) = t^{-\rho_\alpha} \cdot \frac{\alpha_1+1}{\alpha_1+\rho_\alpha+1} \cdot \frac{A_{\tt HP}(t)}{\rho_\alpha}  \cdot \big(1+ o(1)\big)
 \to 
 \begin{cases}
 \displaystyle{ -\frac{\Gamma(2+\alpha_1)}{ \beta_1+\beta_2}  ,} & \mbox{if } \alpha_1=\alpha_2;  \vspace{5pt}\\
 \displaystyle{-\frac{\beta_1\Gamma(2+\alpha_1)}{\beta_2\Gamma(2+2\alpha_1-\alpha_2)}  ,} & \mbox{if } \alpha_1<\alpha_2 .
 \end{cases}
 \eeqlb
 Differentiating both sides of (\ref{eqn.40002}) and then using the equality $C_{\mathcal{I}_R}= \Gamma(1+\alpha_1)C_{\mathcal{I}_R^2}$, we see that  
 \beqlb
 \mathcal{I}_R(t) = C_{\mathcal{I}_R}\cdot t^{\alpha_1} + C_{\mathcal{I}_R^2}\cdot(\alpha_1+\rho_\alpha+1)\cdot t^{\alpha_1+\rho_\alpha}\cdot \ell(t)  +C_{\mathcal{I}_R^2} \cdot t^{\alpha_1+\rho_\alpha+1}\cdot \ell'(t) ,\quad t>0.
 \eeqlb
 A simple modification, together with the equality in (\ref{eqn.40003}), shows that 
 \beqlb
 \mathcal{I}_R(t)=C_{\mathcal{I}_R}\cdot t^{\alpha_1}\cdot
 \Big(1+ \Big(  \frac{\alpha_1+\rho_\alpha+1}{\alpha_1+1} + \frac{1}{\alpha_1+1} \cdot \frac{t\ell'(t)}{\ell(t)} \Big) \cdot \frac{\alpha_1+1}{\alpha_1+\rho_\alpha+1} \cdot \frac{A_{\tt HP}(t)}{\rho_\alpha} \cdot \big(1+ o(1)\big)  \Big) , \quad t>0.
 \eeqlb
 
 By Corollary~\ref{Coro.KaraRep}, the second-order regular variation of $ \mathcal{I}_R$ follows if we can prove $t\ell'(t)/\ell(t) \to 0$ as $t\to\infty$. 
 By the asymptotic results in (\ref{eqn.40003}), it holds if and only if $t\ell'(t) \to 0$ as $t\to\infty$. 
 Using (\ref{eqn.40002}) again, 
 \beqnn
 t\ell'(t) =  \frac{\mathcal{I}_R(t) -C_{\mathcal{I}_R}\cdot t^{\alpha_1} }{C_{\mathcal{I}_R^2} \cdot t^{\alpha_1+\rho_\alpha}} - (\alpha_1+\rho_\alpha+1)\cdot \ell(t), \quad t>0.
 \eeqnn
 Plugging the asymptotic results in (\ref{eqn.40003}) and (\ref{eqn.70001}) into the right side of this equality, we have as $t\to\infty$
 \beqnn
 t\ell'(t) \sim
 \begin{cases}
 \displaystyle{- \frac{\beta_1\beta_2}{C_{\mathcal{I}_R^2} \cdot(\beta_1+\beta_2)^2} + \frac{\Gamma(2+\alpha_1)}{ \beta_1+\beta_2} =0,} & \mbox{if $\alpha_1=\alpha_2$}; \vspace{5pt} \cr
 \displaystyle{ \frac{-\beta_1^2/C_{\mathcal{I}_R^2}+\beta_1\Gamma(2+\alpha_1)}{\beta_2\Gamma(1+2\alpha_1-\alpha_2)} =0,} & \mbox{if $\alpha_1\neq \alpha_2$} .
 \end{cases}
 \eeqnn 
 \qed 
 
 Similarly as in the proof of Theorem~\ref{Thm.2RV.HP}, the second-order approximations for $\mathbf{E}\big[N(t) \big]$ and $ {\rm Var} \big( N(t) \big) $ can be established by applying Proposition~\ref{Prop.HP}, Corollary~\ref{Coro.KaraRep} and Corollary~\ref{Appendix.Coro01} to their exact expressions given in (\ref{eqn.3007}). 
 We just give the respective formulate; the proof is straightforward and is hence omitted. 
 
 \begin{corollary}
 For a strongly critical Hawkes processes of mixed Mittag-Leffler type $N$ with parameter $(\alpha_i,\beta_i)_{i=1,2}$ with $0<\alpha_1\leq \alpha_2<1$ and $\beta_1,\beta_2>0$, we have as $t\to\infty$,
 \beqnn
 \mathbf{E}\big[N(t) \big] -\mu_0\cdot C_{\mathcal{I}^2_R} \cdot t^{\alpha_1+1} 
 \ar\sim\ar 
 \begin{cases}
 	\displaystyle{  \mu_0     \cdot \Big( 1- \frac{ \Gamma(2+\alpha_1)}{\beta_1+\beta_2} \Big)\cdot t ,} & \mbox{if } \alpha_1=\alpha_2;  \vspace{5pt}\\
 	\displaystyle{ -   \frac{\mu_0 \cdot\beta_1\Gamma(2+\alpha_1)}{\beta_2\Gamma(2+2\alpha_1-\alpha_2)} \cdot t^{1+2\alpha_1-\alpha_2},} & \mbox{if } \alpha_1<\alpha_2<2\alpha_1;  \vspace{5pt}\\
 	\displaystyle{  \mu_0 \cdot \Big(1-  \frac{\beta_1}{\beta_2 } \Gamma(2+\alpha_1)\Big) \cdot t ,} & \mbox{if } 2\alpha_1=\alpha_2;  \vspace{5pt}\\
 	\displaystyle{  \mu_0  \cdot t ,} & \mbox{if } 2\alpha_1<\alpha_2,\\
 \end{cases}\cr
 \ar\ar\cr
 \ar\ar\cr
 {\rm Var} \big( N(t) \big) - \mathcal{C}^{\rm Var}_1\cdot t^{3\alpha_1+1}
 \ar\sim\ar  
 \begin{cases}
 	\displaystyle{\mathcal{C}^{\rm Var}_2 \cdot \Big(1- \frac{\beta_1\beta_2}{(\beta_1+\beta_2)^2}  \Big)  \cdot t^{2\alpha_1+1},} & \mbox{if $\alpha_1=\alpha_2$}; \vspace{5pt} \cr
 	\displaystyle{ -\frac{ \mathcal{C}_2^{\rm Var} \cdot\beta_1^2/\beta_2 }{ \Gamma(1+2\alpha_1-\alpha_2)} \cdot t^{4\alpha_1-\alpha_2+1} ,} & \mbox{if $\alpha_1<\alpha_2<2\alpha_1$};  \vspace{5pt} \cr
 	\displaystyle{  \mathcal{C}_2^{\rm Var} \cdot \Big(1- \frac{\beta_1^2}{ \beta_2 }  \Big) \cdot t^{2\alpha_1+1},} & \mbox{if $2\alpha_1=\alpha_2$};  \vspace{5pt} \cr
 	\displaystyle{ \mathcal{C}_2^{\rm Var}  \cdot t^{2\alpha_1+1},} & \mbox{if $2\alpha_1< \alpha_2$};  \vspace{5pt} \cr
 \end{cases}
 \eeqnn 
 with 
 \[
 \mathcal{C}^{\rm Var}_1:= \mu_0\cdot |C_{\mathcal{I}_R}|^3 \cdot \mathrm{B}(1+2\alpha_1,1+\alpha_1) \quad \mbox{and} \quad 
 \mathcal{C}^{\rm Var}_2:= \mu_0\cdot |C_{\mathcal{I}_R}|^2\cdot \big(2 \cdot \mathrm{B}(1+\alpha_1,1+\alpha_1) +\mathrm{B}(1+2\alpha_1,1)\big).
 \]
 \end{corollary} 
 
 \appendix
 
 \renewcommand{\theequation}{A.\arabic{equation}}

  \section{Some properties of regular variation and $\Pi$-variation}\label{Appendix.RV}
  \setcounter{equation}{0}
  
  In this appendix, we provide several results on regular variation and $\Pi$-variation that are used in this paper. 
  The reader is encouraged to refer to the monographs \cite{BinghamGoldieTeugels1987, DeHaanFerreira2006, Resnick2007} for more detailed and extensive results. 
  We assume all functions in this appendix are locally integrable; since we are only interested in their behavior at infinity, we assume integrability on intervals including $0$ as well. 
  
  \begin{proposition}[Uniform convergence theorem]\label{Thm.UniConver}
  	If $F \in {\rm RV}_\alpha^\infty$ with $\alpha\in\mathbb{R}$, we have that (\ref{Def.RV}) holds uniformly in $x \in [a,b]$ if $\alpha=0$, $x\in(0,b]$ if $\alpha >0$ or $x\in [a,\infty)$ if $\alpha<0$ for each $0<a<b<\infty$.
  \end{proposition}
  
  \begin{proposition}[Potter's theorem]\label{Thm.PotterThm}
  	If $F\in {\rm RV}_\alpha^\infty$ with $\alpha\in\mathbb{R}$,     for any $\varepsilon,\delta >0$,  there exists a constant $t_0>0$ such that 
  	$| F(tx)/F(t)-x^\alpha | \leq \varepsilon (x^{\alpha+\delta} \vee x^{\alpha-\delta}) $ for any $x,tx\geq t_0$. 
  \end{proposition}

  \begin{corollary}\label{Appendix.Coro01}
  	For $i=1,2$, let $F_i\in {\rm RV}^\infty_{\alpha_i}$ with $\alpha_i>-1$. Then $F_1*F_2 \in {\rm RV}^\infty_{\alpha_1+\alpha_2+1}$ and as $t\to\infty$,
  	\beqnn
  	F_1*F_2(t) \sim \mathrm{B}(\alpha_1+1,\alpha_2+1) \cdot t \cdot F_1(t) \cdot F_2(t). 
  	\eeqnn
  \end{corollary}
  \proof For each $\epsilon\in(0,1/2)$, by the change of variables we have 
  \beqlb\label{eqn.6001}
  \frac{F_1*F_2(t)}{ t F_1(t)F_2(t)}  \ar=\ar  \int_0^1 \frac{F_1(t(1-s))F_2(ts)}{F_1(t)F_2(t)} ds  \cr
  \ar=\ar  \int_\epsilon^{1-\epsilon} \frac{F_1(t(1-s))F_2(ts)}{F_1(t)F_2(t)} ds  +\int_0^\epsilon \frac{F_1(t(1-s))F_2(ts)}{F_1(t)F_2(t)} ds +  \int_{1-\epsilon}^1 \frac{F_1(t(1-s))F_2(ts)}{F_1(t)F_2(t)} ds . \qquad
  \eeqlb
  By Proposition~\ref{Thm.UniConver},
  \beqlb\label{eqn.6002}
  \lim_{\epsilon \to 0+} \lim_{t\to\infty}\int_\epsilon^{1-\epsilon} \frac{F_1(t(1-s))F_2(ts)}{F_1(t)F_2(t)}  ds 
  =  \lim_{\epsilon \to 0+} \int_\epsilon^{1-\epsilon} (1-s)^{\alpha_1} s^{\alpha_2}ds 
  =\mathrm{B} (\alpha_1+1,\alpha_2+1) . 
  \eeqlb
  For $\delta\in(0,\alpha_2+1)$, by Proposition~\ref{Thm.PotterThm}, there exist two constants $C, t_0>0$ such that   for any $t\geq t_0$ and $ts\geq t_0$, 
  \beqlb\label{eqn.6003}
  |F_2(ts)/F(t)| \leq C \cdot ( s^{\alpha_2-\delta} \vee s^{\alpha_2+\delta} ). 
  \eeqlb
  For $t>t_0/\epsilon$, we have
  \beqlb\label{eqn.6004}
  \int_0^\epsilon F_1(t(1-s))F_2(ts)ds = \int_0^{t_0/t} F_1(t(1-s))F_2(ts)ds  + \int_{t_0/t}^\epsilon F_1(t(1-s))F_2(ts)ds .
  \eeqlb
  From Proposition~\ref{Thm.UniConver} and the local boundedness of $F_2$,  
  \beqnn
  \int_0^{t_0/t} F_1(t(1-s))F_2(ts)ds 
  \ar\leq \ar \sup_{r\in [0,t_0/t]}|F_2(tr)|\cdot \int_0^{t_0/t}  \frac{F_1(t(1-s))}{F_1(t)}  ds \cdot F_1(t) \cr
  \ar\leq\ar  C   \int_0^{t_0/t} (1-s)^{\alpha_1} ds \cdot F_1(t)
  \leq C \cdot  t^{ -1 }  \cdot F_2(t) = o(F_1(t)F_2(t)).
  \eeqnn
  Moreover, by using Proposition~\ref{Thm.UniConver} again and (\ref{eqn.6003}), 
  \beqnn
  \int_{t_0/t}^\epsilon \frac{F_1(t(1-s))F_2(ts)}{F_1(t)F_2(t)}ds 
  \leq C \int_{t_0/t}^\epsilon (1-s)^{\alpha_1} s^{\alpha_2-\delta} ds 
  \leq C \int_0^\epsilon  s^{\alpha_2-\delta} ds \leq C \epsilon^{\alpha_2-\delta+1},
  \eeqnn
  which goes to $0$ as $\epsilon \to 0+$. 
  Putting these two estimates back into (\ref{eqn.6004}), we have 
  \beqnn
  \lim_{\epsilon \to 0+}\limsup_{t\to\infty }  \int_0^{\epsilon} \frac{F_1(t(1-s))F_2(ts) }{F_1(t)F_2(t)}ds  =0. 
  \eeqnn
  Similarly, we also have 
  \beqnn
  \lim_{\epsilon \to 0+}\limsup_{t\to\infty }  \int_{1-\epsilon}^1 \frac{F_1(t(1-s))F_2(ts) }{F_1(t)F_2(t)}ds  =0. 
  \eeqnn
  The whole proof can be ended by taking these two results and  (\ref{eqn.6002}) back into (\ref{eqn.6001}).
  \qed 
  
  \begin{proposition}[Karamata's theorem; see page 25 in \cite{Resnick2007}]\label{Thm.Karamata}
  	The following hold.
  	\begin{enumerate}
  		\item[(1)] If $\alpha\geq -1$, then $F \in {\rm RV}_{\alpha}^\infty$ implies
  		\beqlb\label{KaraThm01}
  		\int_{0}^t F(s)ds \in {\rm RV}_{\alpha+1}^\infty
  		\quad\mbox{and}\quad
  		\lim_{t\to\infty} \frac{tF(t)}{\int_{0}^t F(s)ds} = \alpha+1.
  		\eeqlb
  		
  		\item[(2)]
  		If $\alpha<-1$ or if $\alpha=-1$ and $\int_{t_0}^\infty F(s)ds<\infty$ for some $t_0>0$, then $F\in {\rm RV}_{\alpha}^\infty$ implies
  		\beqlb\label{KaraThm02}
  		\int_t^\infty F(s)ds \in {\rm RV}_{\alpha+1}^\infty
  		\quad\mbox{and}\quad \lim_{t\to\infty} \frac{tF(t)}{\int_t^\infty  F(s)ds} = -\alpha-1.
  		\eeqlb
  		
  	\end{enumerate}
  \end{proposition}

  \begin{proposition}[Karamata's Tauberian theorem]\label{Thm.KaramataTauberian}
  	For $\alpha>-1$, the following hold. 
  	\begin{enumerate}
  		\item[(1)] If $F\in {\rm RV}^\infty_\alpha$, then $\hat{F}(1/\cdot) \in {\rm RV}^\infty_\alpha$ and $\hat{F}(1/\lambda) \sim \Gamma(1+\alpha) F(\lambda)$ as $\lambda \to \infty$. 
  		
  		\item[(2)] If $F$ is eventually monotone and $\hat{F}(1/\cdot) \in {\rm RV}^\infty_\alpha$, then $F\in {\rm RV}^\infty_\alpha$. 
  		
  	\end{enumerate}
  \end{proposition}
  
  \begin{proposition}[Corollary~8.1.7 in \cite{BinghamGoldieTeugels1987}]\label{Prop.A.7}
  	Let $F$ be a probability distribution on $\mathbb{R}_+$, let $\ell\in {\rm RV}_0^\infty$ and $\alpha\in[0,1]$. The following two statements are equivalent. 
  	\begin{enumerate}
  		\item[(1)] $\displaystyle{  \int_0^\infty \big(1- e^{-t / \lambda} \big) dF(t)  \sim   \lambda^{-\alpha} \ell(\lambda)} $ as $\lambda\to \infty$; 
  		
  		\item[(2)] one of the following holds as $t\to\infty$:
  		\begin{enumerate}
  			\item[(2.i)] 	$1-\displaystyle{F(t) \sim \frac{ t^{-\alpha}\ell(t)}{\Gamma(1-\alpha)} }$ if $\alpha\in[0,1)$;
  			
  			\item[(2.ii)] $\displaystyle{\int_0^t s  dF(s) \sim  \ell(t)}$ if $\alpha=1$. 
  		\end{enumerate}

  	\end{enumerate}

  \end{proposition}
  
  \begin{proposition}[Theorem~8.1.6 in \cite{BinghamGoldieTeugels1987}]\label{Prop.A.7.1}
  	Let $F$ be a probability distribution on $\mathbb{R}_+$ with finite mean and let $\ell\in {\rm RV}_0^\infty$ and $\alpha\in[1,2]$. The following two statements are equivalent. 
  	\begin{enumerate}
  		\item[(1)] $\displaystyle{ \int_0^\infty \big(e^{-t/\lambda}-1+ \frac t \lambda \big) dF(t)  \sim   \lambda^{-\alpha} \ell(\lambda)} $ as $\lambda\to \infty$; 
  		
  		\item[(2)] one of the following holds as $t\to\infty$:
  		\begin{enumerate}
  			 
  			\item[(2.i)] $\displaystyle{\int_t^\infty s  dF(s) \sim  \ell(t)}$ if $\alpha=1$;
  			
  			\item[(2.ii)] 	$1-\displaystyle{F(t) \sim \frac{ t^{-\alpha}\ell(t)}{\Gamma(1-\alpha)} }$ if $\alpha\in(1,2)$;
  			
  			\item[(2.iii)] $\displaystyle{\int_0^t s^2  dF(s) \sim  2\cdot \ell(t)}$ if $\alpha=2$.
  		\end{enumerate}

  	\end{enumerate}

  \end{proposition}

  \begin{proposition}[Theorem B.2.12 in \cite{DeHaanFerreira2006}]\label{Prop.A.5}
  	For $A\in {\rm RV}^\infty_0$ and $t_0\geq 0$,  assume that $\overline G(t):=t^{-1}\cdot \mathcal{I}^{t_0,\uparrow}_{\overline F,1}(t)$ is well defined for $t\geq t_0$, we have $\overline F\in \Pi^\infty(A)$ if and only if $\overline G\in {\rm RV}^\infty_0$ and then $A(t)\sim \overline G(t)$ as $t\to\infty$.  
  \end{proposition}
  
  \begin{proposition}[Corollary B.2.13 in \cite{DeHaanFerreira2006}]\label{Prop.A.6}
  	For $A\in {\rm RV}^\infty_0$, if $\overline F\in \Pi^\infty(A)$, then  $\overline F(t)\to \overline F(\infty)\in [0,\infty]$. 
  	Moreover, $\overline F\in {\rm RV}^\infty_{0}$ if $\overline F(\infty)=\infty$ and $\overline F(\infty)-\overline F(t) \in {\rm RV}^\infty_0$ if $\overline F(\infty)<\infty$. 
  \end{proposition}
  
  \begin{proposition}[Wiener-Pitt theorem; see Lemma 2.33 in \cite{GelukdeHaan1987}]\label{Wiener-Pitt.Thm}
  	Let $g$ be a bounded, slowly increasing (or decreasing) function on $\mathbb{R}_+$ and $k_0 $ be an integrable kernel on $\mathbb{R}_+$ with 
  	$\int_0^\infty k_0(s)s^{-\mathtt{i} z}ds \neq 0$ for all $z\in \mathbb{R}$.
  	For some $c \in \mathbb{R}$, 
  	we have $g(t)\to c$ as $t\to\infty$ if 
  	\beqnn
  	\lim_{\lambda \to\infty} k_0\overset{\rm M}* g(\lambda) =  c\cdot \int_0^\infty k_0(s)ds. 
  	\eeqnn

  \end{proposition}

 \bibliographystyle{plain}

 \bibliography{Reference}

  	\end{document}